\documentclass[11pt]{article}
\usepackage[T1]{fontenc}
\usepackage{lmodern}
\usepackage[margin=27mm]{geometry}
\usepackage{amsmath,amssymb,amsthm,mathtools}
\usepackage{microtype}
\usepackage{xcolor}
\usepackage{tocloft}
\usepackage{newtxtext}
\usepackage{newtxmath}
\usepackage[colorlinks=true,linkcolor=blue,citecolor=blue,urlcolor=blue]{hyperref}

\numberwithin{equation}{section}

\newtheorem{theorem}{Theorem}[section]
\newtheorem{lemma}[theorem]{Lemma}
\newtheorem{proposition}[theorem]{Proposition}
\newtheorem{corollary}[theorem]{Corollary}

\theoremstyle{remark}
\newtheorem{remark}[theorem]{Remark}

\newcommand{\R}{\mathbb R}

\newcommand{\C}{\mathbb C}

\newcommand{\dd}{\,\mathrm d}
\newcommand{\ii}{\mathrm i}
\DeclareMathOperator{\im}{Im}
\DeclareMathOperator{\re}{Re}
\DeclareMathOperator{\supp}{supp}

\newcommand{\D}{\langle D\rangle}
\newcommand{\norm}[2]{\|#1\|_{#2}}

\title{A new proof of global well-posedness for the mass-critical
	nonlinear Schr\"odinger equation and its extension to Hamiltonian systems}
\author{Jiajun Wang}
\date{}
\newcommand{\Sclass}{\mathcal S}
\begin{document}
\maketitle
\begin{abstract}
We give a new proof of Dodson's theorem~\cite{Dodson} on global
well-posedness and scattering for the defocusing cubic Schr\"odinger
equation on $\R^2$ with arbitrary initial data in $L^2(\R^2)$.
The proof combines a modified energy with the bilinear virial identity
of Planchon--Vega~\cite{PV}. For initial data in $H^s(\R^2)$ with
$s>0$, we prove global well-posedness by estimating the high-frequency
Duhamel integral on time intervals whose lengths grow at most
polynomially. For general $L^2(\R^2)$ data,
the scaling argument gives no such polynomial bound. We instead apply
the local mass and momentum identities to the forced high-frequency
wave to obtain a Duhamel estimate uniform in the length of the time
interval. The proof is direct, uses no argument by contradiction, and requires
neither profile decompositions nor inverse Strichartz estimates.
It provides a new framework that may also apply to other dispersive equations.

The main extension is an arbitrary-data theorem for general Hamiltonian couplings. We prove global well-posedness and scattering for arbitrary
$L^2(\R^2;\C^m)$ data for finite Hamiltonian systems with a nonnegative real quartic potential of bidegree $(2,2)$.
This extends the $H^s$, $s>4/7$, result of Ma--Song--Zhao~\cite{MSZ}
to the mass-critical regularity and to more general quartic interactions.

For every $s<0$, we also construct an initial datum
$u_0\in H^s(\R^2)\setminus L^2(\R^2)$ and a distributional
solution satisfying
$\lim_{t\uparrow1}\|u(t)\|_{H^s}=\infty$.
The construction combines the choice of three resonant frequencies
in \cite[Example~2.3]{CDS2012}, quadratic phases, spatial
separation, and a Strichartz perturbation argument.
It applies to both the focusing and defocusing scalar cubic
equations and extends through an invariant one-dimensional
subspace to every nonzero real Hamiltonian potential of
bidegree $(2,2)$, without a positivity assumption.
\end{abstract}
\tableofcontents

\vspace{10pt}
\section{Introduction and statement of the main results}
\vspace{7pt}
\subsection{2D mass-critical NLS}
We consider the 2D defocusing cubic nonlinear Schr\"odinger equation
\begin{equation}\label{eq:nls}
	\left\{
	\begin{array}{ll}
		\ii\partial_tu+\Delta u=|u|^2u,
		&(t,x)\in\R\times\R^2,\\
		u(0,x)=u_0(x),
		&x\in\R^2.
	\end{array}
	\right.
\end{equation}
For sufficiently regular, decaying solutions, the conserved mass
and energy are
\begin{equation*}
	M(u):=\norm{u}{L^2}^2,\qquad
	E(u):=\frac12\norm{\nabla u}{L^2}^2
	+\frac14\norm{u}{L^4}^4.
\end{equation*}
Equation~\eqref{eq:nls} is invariant under the scaling
\begin{equation*}
	u_\lambda(t,x):=\lambda^{-1}u(t/\lambda^2,x/\lambda),
	\qquad
	\norm{u_\lambda(0)}{\dot H^s}
	=\lambda^{-s}\norm{u_0}{\dot H^s},
	\qquad \lambda>0.
\end{equation*}
In particular,
\begin{equation*}
	\norm{u_\lambda(0)}{L^2}=\norm{u_0}{L^2},
	\qquad M(u_\lambda(0))=M(u_0).
\end{equation*}
Thus the scaling preserves the mass, and the critical Sobolev
regularity is $s=0$. This is why \eqref{eq:nls} is called
mass-critical.

We first establish global well-posedness at positive Sobolev
regularity $s>0$. This introduces the modified-energy argument before
the additional difficulties at the critical regularity $s=0$ are addressed.

\begin{theorem}[Global well-posedness at positive regularity]\label{thm:positive}
	For every $s>0$ and $u_0\in H^s(\R^2)$,
	\eqref{eq:nls} has a unique global solution
	\begin{equation*}
		u\in C(\R;H^s(\R^2))
		\cap L^4_{\mathrm{loc}}(\R;W^{s,4}(\R^2)).
	\end{equation*}
	The solution depends continuously on the initial datum in
	$C(I;H^s(\R^2))\cap L^4(I;W^{s,4}(\R^2))$
	for every compact time interval $I\subset\R$.
\end{theorem}

We then extend the argument to the mass-critical regularity
and establish scattering in both time directions.

\begin{theorem}[Global well-posedness and scattering at critical regularity]\label{thm:endpoint}
	For every $u_0\in L^2(\R^2)$,
	\eqref{eq:nls}  has a unique global solution
	\begin{equation*}
		u\in C(\R;L^2(\R^2))
		\cap L^4(\R\times\R^2).
	\end{equation*}
	The solution depends continuously on the initial datum in
	the uniform-in-time $L^2(\R^2)$ norm and the $L^4(\R\times\R^2)$ norm.
	
	Moreover, $u$ scatters in $L^2(\R^2)$ in both time directions:
	there exist unique $u_+,u_-\in L^2(\R^2)$ such that
	\begin{equation*}
		\lim_{t\to+\infty}
		\norm{u(t)-e^{\ii t\Delta}u_+}{L^2}=0,
		\qquad
		\lim_{t\to-\infty}
		\norm{u(t)-e^{\ii t\Delta}u_-}{L^2}=0.
	\end{equation*}
\end{theorem}

\begin{remark}\label{persist}
	Theorem~\ref{thm:positive} also follows from
	Theorem~\ref{thm:endpoint} by persistence of regularity. Indeed, for $u_0\in H^s\subset L^2$ and $I=\R$,
	the solution $u$ has finite $L^4(I\times\R^2)$ norm.
	Subdivide $I$ into finitely many intervals $J$ on which this norm
	is sufficiently small. For $t_0\in J$, Strichartz estimates and
	the fractional product estimate give
	\begin{equation}\label{a}
		\norm{u}{L^\infty(J;H^s(\R^2))}
		+\norm{u}{L^4(J;W^{s,4}(\R^2))}
		\lesssim_s
		\norm{u(t_0)}{H^s}
		+\norm{u}{L^4(J\times\R^2)}^2
		\norm{u}{L^4(J;W^{s,4}(\R^2))}.
	\end{equation}
	Absorption and iteration over these intervals propagate
	$H^s$ regularity throughout $I$; local $H^s$ well-posedness
	also gives uniqueness and continuous dependence.
	
	We nevertheless give a proof of
	Theorem~\ref{thm:positive} first. The positive-regularity estimates
	introduce the main argument in a simpler setting and identify
	precisely which steps require modification at $s=0$.
	This order makes the passage to the critical case more transparent.
\end{remark}
\begin{remark}
	The proof of Theorem~\ref{thm:positive} establishes
	global well-posedness and
	$u\in L^4_{\mathrm{loc}}(\R;W^{s,4}(\R^2))$,
	but does not establish a global spacetime bound or scattering.
	Once Theorem~\ref{thm:endpoint} is available, its bound
	$\norm{u}{L^4(\R\times\R^2)}<\infty$, together with the same estimates as in (\ref{a}), yields
	\begin{equation*}
		u\in L^\infty(\R;H^s(\R^2))
		\cap L^4(\R;W^{s,4}(\R^2)).
	\end{equation*}
	Consequently, $u$ scatters in $H^s(\R^2)$ in both time directions:
	there exist unique $u_\pm\in H^s(\R^2)$ such that
	\begin{equation*}
		\norm{u(t)-e^{\ii t\Delta}u_\pm}{H^s}
		\longrightarrow0
		\qquad(t\to\pm\infty).
	\end{equation*}
	These stronger conclusions therefore follow after the critical
	theory is established.
\end{remark}

\vspace{7pt}
\subsection{Hamiltonian systems}
The proof of Theorem~\ref{thm:endpoint} (with some modifications) also applies to a broader
class of Hamiltonian Schr\"odinger systems. We now replace the
scalar cubic nonlinearity by the gradient of a real quartic
potential and establish global well-posedness and scattering
for arbitrary $L^2$ initial data under the assumption that
the potential is nonnegative. 

Let $m<\infty$, let $W:\C^m\to\R$ be a polynomial of bidegree $(2,2)$, and set
\begin{equation}\label{sys-eq:system}
	\left\{
	\begin{array}{ll}
		\ii\partial_tu_j+\Delta u_j=F_j(u),
		&(t,x)\in\R\times\R^2,\\
		u_j(0,x)=u_{0,j}(x),
		&x\in\R^2,
	\end{array}
	\right.
	\qquad j=1,\ldots,m,
\end{equation}
where $F_j(z):=\partial_{\bar z_j}W(z)$.
Bidegree $(2,2)$ is equivalent here to real quartic homogeneity together with $W(e^{\ii\vartheta}z)=W(z)$. Our assumption is
\begin{equation}\label{sys-eq:nonnegative}
 W(z)\ge0\quad\text{for every }z\in\C^m.
\end{equation}
No convexity or strict positivity is assumed.  This condition is necessary for an arbitrary-data theorem: Proposition~\ref{sys-prop:negative} gives smooth blow-up data when $W$ takes a negative value.

\begin{theorem}[Global well-posedness and scattering]\label{sys-thm:main}
	Under \eqref{sys-eq:nonnegative}, for every $s\ge0$ and
	$u_0\in H^s(\R^2;\C^m)$, \eqref{sys-eq:system} has a unique
	global solution
	\begin{equation*}
		u\in C(\R;H^s(\R^2;\C^m))
		\cap L^\infty(\R;H^s(\R^2;\C^m))
		\cap L^4(\R;W^{s,4}(\R^2;\C^m)).
	\end{equation*}
	The solution depends continuously on the initial data in the above space.
	
	Moreover, $u$ scatters in $H^s(\R^2;\C^m)$ in both time
	directions: there exist unique $u_+,u_-\in H^s(\R^2;\C^m)$
	such that
	\begin{equation*}
		\lim_{t\to+\infty}
		\norm{u(t)-e^{\ii t\Delta}u_+}{H^s}=0,
		\qquad
		\lim_{t\to-\infty}
		\norm{u(t)-e^{\ii t\Delta}u_-}{H^s}=0.
	\end{equation*}

\end{theorem}

\begin{corollary}\label{sys-cor:matrix}
Let $\beta$ be a fixed real symmetric $m\times m$ matrix satisfying the copositivity condition: 
\begin{equation}\label{sys-eq:copositive}
 \sum_{j,k}\beta_{jk}r_jr_k\ge0\qquad(r_1,\ldots,r_m\ge0).
\end{equation}
Then Theorem~\ref{sys-thm:main} applies to
\begin{equation*}
 \ii\partial_tu_j+\Delta u_j=
 \left(\sum_k\beta_{jk}|u_k|^2\right)u_j.
\end{equation*}
In particular, the conclusion applies whenever $\beta$ is a real
symmetric matrix that is either entrywise nonnegative or positive
semidefinite.
\end{corollary}

\begin{remark}\label{rem:scope}
	For the case $m=2$, the copositivity condition
	\eqref{sys-eq:copositive} is equivalent to
	\begin{equation*}
		\beta_{11}\ge0,\qquad \beta_{22}\ge0,\qquad
		\beta_{12}\ge-\sqrt{\beta_{11}\beta_{22}}.
	\end{equation*}
\end{remark}

\vspace{7pt}
\subsection{Finite-time blow-up below the mass space}
The next statements concern one solution for each $s<0$, rather than a family of norm-inflation examples. We also allow the focusing scalar sign:
\begin{equation}\label{neg-eq:scalar}
 \ii\partial_tu+\Delta u=\mu|u|^2u,\qquad \mu\in\{-1,1\}.
\end{equation}
Here $\Sclass$ denotes the Schwartz space, and $\Sclass'$ its dual.

\begin{theorem}\label{neg-thm:scalar}
For every $s<0$, $\varepsilon>0$, and either sign in \eqref{neg-eq:scalar}, there exists $u_0\in H^s(\R^2)\setminus L^2(\R^2)$ with $\|u_0\|_{H^s}<\varepsilon$ and a distributional solution $u$ on $[0,1)$ such that
\begin{equation}\label{neg-eq:blowup}
 u\in C([0,1);H^s),\qquad
 \lim_{t\uparrow1}\|u(t)\|_{H^s}=\infty.
\end{equation}
The solution also belongs to $L^4([0,1];L^4_{\mathrm{loc}}(\R^2))$
and extends to an element of
$C([0,1];L^2_{\mathrm{loc}}(\R^2))
\cap C([0,1];\Sclass'(\R^2))$,
with $u(1)\notin H^s(\R^2)$.
 Then the cubic nonlinearity in (\ref{neg-eq:scalar}) is understood as the ordinary product. 
\end{theorem}

The above construction directly extends to the Hamiltonian systems (\ref{sys-eq:system}) through an invariant one-dimensional subspace.
\begin{theorem}\label{neg-thm:system}
Let $m<\infty$ and let $W:\C^m\to\R$ be a nonzero real polynomial of bidegree $(2,2)$. For every $s<0$ and $\varepsilon>0$, system \eqref{sys-eq:system} has an initial datum $u_0\in H^s(\R^2;\C^m)\setminus L^2(\R^2;\C^m)$, with $\|u_0\|_{H^s}<\varepsilon$, and a distributional solution $u$ satisfying all the conclusions of Theorem~\ref{neg-thm:scalar}, with vector-valued spaces. 
\end{theorem}

\begin{corollary}\label{neg-cor:matrix}
Let $\beta$ be any nonzero real symmetric $m\times m$ matrix. For every $s<0$, the system
\begin{equation*}
 \ii\partial_tu_j+\Delta u_j=
 \left(\sum_k\beta_{jk}|u_k|^2\right)u_j
\end{equation*}
has the blow-up asserted in Theorem~\ref{neg-thm:system}. Copositivity is not needed for this blow-up statement.
\end{corollary}

\vspace{7pt}
\subsection{Earlier results and the contribution of the paper}
For the scalar equation \eqref{eq:nls}, the $I$-method of
Colliander--Keel--Staffilani--Takaoka--Tao~\cite{CKSTT}
gave global well-posedness in $H^s(\R^2)$ for $s>4/7$.
Their resonant decomposition and corrected
energy~\cite{CKSTT2} improved this range to $s>1/2$.
Dodson~\cite[Theorem 1.3]{Dodson} subsequently established
global well-posedness and scattering for every
$L^2(\R^2)$ initial datum.
Thus the conclusions of
Theorems~\ref{thm:positive} and~\ref{thm:endpoint}
are already known.
Our contribution for \eqref{eq:nls} is an independent
modified-energy proof which avoids inverse Strichartz estimates and
profile decompositions and extends to the more general
nonlinearities in \eqref{sys-eq:system}.

For the two-component case ($m=2$) of \eqref{sys-eq:system},
Ma--Song--Zhao~\cite[Theorem 1]{MSZ} proved global
well-posedness in $H^s(\R^2;\C^2)$ for $s>4/7$ when
\begin{equation*}
	F_1(z)=(\mu_1|z_1|^2+\beta|z_2|^2)z_1,\qquad
	F_2(z)=(\mu_2|z_2|^2+\beta|z_1|^2)z_2,
	\qquad
	\mu_1,\mu_2>0,\quad
	\beta>-\sqrt{\mu_1\mu_2}.
\end{equation*}
Corollary~\ref{sys-cor:matrix} extends this conclusion
to $s=0$, includes the boundary
$\beta=-\sqrt{\mu_1\mu_2}$, and allows any finite
number of components.
Negative cross-coupling was already permitted
in~\cite{MSZ}; the improvement concerns the regularity
and the range of admissible systems.

Yang--Zhao~\cite[Theorem 1.1]{YZ} established global
well-posedness and scattering in $L_x^2h_j^1$
for an infinite resonant system.
Its one-dimensional discrete resonance relation
reduces the nonlinearity to
\begin{equation*}
	F_j(z)=\left(2\sum_k|z_k|^2-|z_j|^2\right)z_j.
\end{equation*}
Every finite set of components is invariant.
The resulting finite system is a special case of
\eqref{sys-eq:system}, with
\begin{equation*}
	W(z)=\left(\sum_j|z_j|^2\right)^2
	-\frac12\sum_j|z_j|^4\ge0.
\end{equation*}
Theorem~\ref{sys-thm:main} recovers global well-posedness
and scattering for these finite systems and extends
these conclusions to every potential satisfying
\eqref{sys-eq:nonnegative}.

Below the mass space $L^2$, Christ--Colliander--Tao~\cite{CCT}
established failure of uniform continuity for the
scalar cubic Schr\"odinger equation, including
\eqref{eq:nls}.
Oh~\cite[Theorem 1.1]{Oh2017} proved norm inflation
at every prescribed $H^s(\R^2)$ initial datum
for every $s<0$, for both \eqref{eq:nls}
and its focusing counterpart.
These results concern families of solutions.
In contrast, Theorem~\ref{neg-thm:scalar} constructs
a single initial datum in $H^s(\R^2)\setminus L^2(\R^2)$
whose solution satisfies the finite-time blow-up
\eqref{neg-eq:blowup}.
Theorem~\ref{neg-thm:system} extends this construction
to \eqref{sys-eq:system} for every nonzero real
potential of bidegree $(2,2)$, without assuming
\eqref{sys-eq:nonnegative}.

\vspace{7pt}
\subsection{Road map of the proofs}

\textit{Global well-posedness and scattering.}
Let $P_N$ be a projection to frequencies $|\xi|\lesssim  N$
and set $Q_N:=\operatorname{Id}-P_N$. Writing $F(u):=|u|^2u$, we
consider the modified energy
\begin{equation*}
	\mathcal E_N(u)
	=\frac12\norm{\nabla P_Nu}{L^2(\R^2)}^2
	+\frac14\norm{P_Nu}{L^4(\R^2)}^4
	+\frac{N^2}{2}\norm{Q_Nu}{L^2(\R^2)}^2 .
\end{equation*}
Differentiating this modified energy in time gives
\begin{equation*}
	\frac{\dd}{\dd t}\mathcal E_N(u)
	=\im\int_{\R^2}
	\overline{-\Delta P_Nu+F(P_Nu)}
	\bigl(P_NF(u)-F(P_Nu)\bigr)
	+N^2\overline{Q_Nu}\,Q_NF(u)\dd x .
\end{equation*}
The difference $P_NF(u)-F(P_Nu)$ vanishes when all three
input frequencies are at most $N/3$. The last term already
contains $Q_Nu$. Thus each integral has a factor at frequency
$\gtrsim N$, and frequency conservation supplies a second one.
Splitting the remaining factors at $L=N^{1-\varepsilon}$,
Lemma~\ref{source} pairs a low factor with a high factor in
$L^2_{t,x}$ and estimates factors above $L$ in $L^4_{t,x}$.
This is why \eqref{bootstrap} includes, alongside the energy,
the two quantities
\begin{equation*}
	Y(Q_{aN}u)
	=\sup_{d\in\R^2}
	\norm{(P_Lu)(t,x+d)(Q_{aN}u)(t,x)}{L^2(J\times\R^2)},
	\qquad
	\norm{Q_Lu}{L^4(J\times\R^2)} .
\end{equation*}

To improve the bound for $Y(Q_{aN}u)$, the key observation
is that the bilinear virial identity controls a derivative
of the product. If $h$ has frequency $H\gtrsim N$, then
$(P_Lu)\bar h$ has frequency at least $H/2-2L\gtrsim H$.
We can therefore remove the derivative at a gain of $H^{-1}$.

To improve the $L^4$ bound for $Q_Lu$, \eqref{w} already
controls $Q_Nu$, so it remains to estimate
$(P_Nu)_{\ge r}=\mathbf1_{\{|D|\ge r\}}P_Nu$ for
$L/2<r\le L$. Fix $\rho<r$. There are two cases.

 \emph{First case}: if one of the three
factors in $(P_Nu)\overline{(P_Nu)}(P_Nu)$ has frequency
below $\rho$, the relation $\xi_1-\xi_2+\xi_3-\xi_4=0$,
with $|\xi_4|\ge r$, allows us to form two pairs with
frequency differences $\gtrsim r-\rho$. Plancherel then
bounds each product by $(r-\rho)^{-1}$ times the $L^2$ norm
of its gradient. The bilinear virial estimates control both gradients and yield
the square-summable errors in \eqref{separated}, which is easy to control.
\emph{Second case}: if all three
factors have frequency
above $\rho$, their product is bounded by
$\norm{(P_Nu)_{\ge\rho}}{L^4(I\times\R^2)}^3$. 

Partition $J$ into subintervals $I$,
so that $\norm{(P_Nu)_{\ge r}}{L^4(I\times\R^2)}
\sim L^{-1}$. On each interval $I$
where the \emph{second case} dominates, we have
\begin{equation*}
	\norm{(P_Nu)_{\ge r}}{L^4(I\times\R^2)}^4
	\sim L^{-4}
	\lesssim L^{-8/3}
	\norm{(P_Nu)_{\ge\rho}}{L^4(I\times\R^2)}^4 .
\end{equation*}
Summing this bound and choosing $r,\rho$ as in
\eqref{eq:positive-two-radii} bounds the contribution of
these intervals by
$ L^{-8/3}\norm{(P_Nu)_{\ge r}}{L^4(J\times\R^2)}^4$.
For large $L$ this term is absorbed in \eqref{count},
giving \eqref{tail-improved} and closing the bootstrap.
The quartic term of $\mathcal E_N$ and \eqref{w} now control
$P_Nu$ and $Q_Nu$, respectively, in $L^4(J\times\R^2)$ on
each finite interval $J$. Pseudoconformal invariance then
gives scattering; see also~\cite[Section~3]{TaoPC}.

\textit{Finite-time $H^s$-norm blow-up.}
For $s<0$, a high-frequency packet can have large mass and
small $H^s$ norm. To make its contribution large at $t=1$,
Proposition~\ref{neg-prop:profile} uses three packets that
meet near $t=1$. In the rescaled variables, the three packets have central frequencies
$Ne_1$, $N(e_1+e_2)$, and $Ne_2$, with $e_1=(1,0)$ and $e_2=(0,1)$,
as in the configuration of \cite[Example~2.3]{CDS2012}.
The corresponding spatial and temporal oscillations cancel in the
cubic interaction, since
\begin{equation*}
	e_1-(e_1+e_2)+e_2=0,\qquad
	|e_1|^2-|e_1+e_2|^2+|e_2|^2=0 .
\end{equation*}
Before the packets meet, we evolve them separately. Their high
frequencies keep the $H^{-\sigma}$ norm of $Z_n$ small. Spatial separation makes every mixed
cubic term small, and the quadratic phases make the packets
spread backwards in time, giving the small $L^2_tL^\infty_x$
norm needed for the correction.

When the packets meet, we take a finite Duhamel iterate of their sum.
In the first iterate, the cubic interaction
$(b\otimes a)\overline{(b\otimes b)}(a\otimes b)$ contributes
$\bigl(\int_{\R}|a|^2|b|^2\dd x\bigr)^2$ when tested against
$a\otimes a$. The overlap of the packets gives a positive lower
bound for its time integral, and the choice of $\theta$ ensures
that the full pairing retains this lower bound. The later iterates preserve the resulting contribution,
which after scaling gives the large terminal pairing used in the
blow-up construction.

Finally, following the author's construction for the periodic
Zakharov system~\cite[Sections~14--15]{WangZakharov}, we add
one translated profile at a time and correct the sum to an
exact solution before choosing the next profile. The corrections
are summable, and the translations preserve each profile's
large terminal pairing while making its pairings with the
other test functions small. The limit therefore belongs to
$C([0,T];H^s(\R^2))$ for every $T<1$, while its distributional
value at $t=1$ lies outside $H^s(\R^2)$. This gives
\eqref{neg-eq:blowup}.

\vspace{7pt}
\subsection{Organization of the proofs}
Section~\ref{sec:tools} collects the analytic estimates and establishes the local theory.
Section~\ref{sec:energy} estimates the variation of the modified energy and derives
the key mass--momentum estimates. Section~\ref{sec:positive} proves
Theorem~\ref{thm:positive} on the polynomial time scale given by positive Sobolev regularity.
Section~\ref{sec:endpoint} replaces the time restriction by an estimate
on arbitrary intervals. Section~\ref{sec:global} then proves Theorem~\ref{thm:endpoint},
including scattering. Section~\ref{sec:vector} proves Theorem~\ref{sys-thm:main}
by deriving the corresponding virial identity and checking the nonlinear estimates.
Section~\ref{sec:sharpness} discusses the necessity
of the sign condition and the algebraic reduction
to scalar equations.
Sections~\ref{sec:negative-profile}--\ref{sec:negative-limit}
prove Theorem~\ref{neg-thm:scalar}.
Section~\ref{sec:negative-profile} constructs the individual approximate profiles.
Section~\ref{sec:negative-correction} shows how to add one such profile
to a smooth exact solution and correct the resulting error.
Section~\ref{sec:negative-limit} iterates this step and passes to the limit
to obtain one initial datum with finite-time $H^s$ norm blow-up.
Section~\ref{sec:negative-systems} transfers
this construction to \eqref{sys-eq:system},
proving Theorem~\ref{neg-thm:system}.

\vspace{10pt}
\noindent
\textbf{Notation.}
\begin{itemize}
	\item By $u\in C([0,T];B)$ (or $L^p([0,T];B)$) for a
	Banach space $B$, we mean that $u$ is a continuous
	(or Bochner $L^p$) map from $[0,T]$ to $B$;
	see \cite[page 301]{evans2022partial}.
	
	\item By $A\lesssim B$ (resp. $A\sim B$), we mean that
	there exists a positive constant $C$ such that $A\le CB$
	(resp. $C^{-1}B\le A\le CB$).
	If the constant $C$ depends on $p$, then we write
	$A\lesssim_p B$ (resp. $A\sim_p B$).

	\item By $A\ll B$, we mean that $A/B$ is sufficiently small.
	
		\item By $\norm f{L^p}$ and $\langle f,g\rangle$, we mean
	$\norm f{L^p(\R^2)}$ and
	$\int_{\R^2}f(x)\overline{g(x)}\dd x$, respectively.
	Unless otherwise indicated, all spatial norms are taken on $\R^2$.
	
	\item By $|x|$ and $\langle x\rangle$, we mean the Euclidean norm
	of $x$ and $(1+|x|^2)^{1/2}$, respectively. For
	$u=(u_1,\ldots,u_m)\in\C^m$, we likewise write
	$|u|=(\sum_{j=1}^m|u_j|^2)^{1/2}$.
	
	\item By $\widehat f$, we mean the spatial Fourier
	transform of $f$.

\end{itemize}

\vspace{10pt}
\noindent
\textbf{Disclosure of AI assistance.}
The author used ChatGPT-6 Astra during the preparation and
revision of this paper to assist with developing several
arguments (especially in the endpoint case and the construction of blow-up), checking calculations, and revising the exposition
and \LaTeX{} code.
The author takes full responsibility for the mathematical
claims, proofs, and references presented in the paper.

\vspace{10pt}
\section{Analytic tools and local solutions}\label{sec:tools}

Fix a real radial smooth function $\chi$, with $0\le\chi\le1$, equal to one on $|\xi|\le1$ and zero on $|\xi|\ge2$. Set
\begin{equation*}
 P_L:=\chi(D/L),\qquad Q_L:=1-P_L.
\end{equation*}
\begin{equation*}
 u_1:=P_1u,\qquad u_L:=(P_L-P_{L/2})u, \quad \forall L\ge 2.
\end{equation*}
Thus $u=\sum_Lu_L$, and $\supp\widehat u_L\subset\{L/2\le|\xi|\le2L\}$ for $L\ge2$. The cutoffs $P_L,Q_L$ also make sense for non-dyadic $L$.

\vspace{7pt}
\subsection{Linear estimates}
The following is the two-dimensional Strichartz estimate from Keel--Tao~\cite[Theorem 1.2]{KT}. 

\begin{proposition}\label{prop:strichartz}
	We call a pair $(q,r)$ Schr\"odinger admissible if
	\begin{equation}\label{eq:admissible}
		\frac{2}{q}+\frac{2}{r}=1,\qquad
		2\le q,r\le\infty,\qquad (q,r)\ne(2,\infty).
	\end{equation}
	Let $S(t)=e^{\ii t\Delta}$.
	For every admissible pair $(q,r)$
	in the sense of \eqref{eq:admissible},
	\begin{equation}\label{eq:strichartz}
		\norm{S(t)f}{L^q(\R;L^r(\R^2))}
		\lesssim_{q,r}\norm{f}{L^2(\R^2)}.
	\end{equation}
	If $(\widetilde q,\widetilde r)$ is another
	admissible pair, $J\subset\R$ is an interval,
	and $a\in J\cap\R$, then
	\begin{equation}\label{eq:inhomogeneous}
		\left\|\int_a^t S(t-\tau)G(\tau)\dd\tau
		\right\|_{L^q(J;L^r(\R^2))}
		\lesssim_{q,r,\widetilde q,\widetilde r}
		\norm{G}{L^{\widetilde q'}(J;L^{\widetilde r'}(\R^2))}.
	\end{equation}
	In particular, the right side $\norm G{L^{4/3}(J\times \R^2)}$ controls the output in $L^4(J\times\R^2)$ and $L^\infty(J;L^2(\R^2))$.
\end{proposition}

We next recall the estimate from
Christ--Kiselev~\cite[Theorem 1.2]{CK}, together with
its Banach-valued extension discussed there.

\begin{proposition}\label{prop:ck}
	Let $X,Y$ be Banach spaces, let $1\le p<q\le\infty$,
	and let $K:\R^2\to\mathcal L(X,Y)$ be locally integrable.
	Suppose that
	\begin{equation*}
		\left\|\int_\R K(t,\tau)G(\tau)\dd\tau
		\right\|_{L^q(\R;Y)}
		\le A\norm{G}{L^p(\R;X)}.
	\end{equation*}
	Then we have 
	\begin{equation*}
		\left\|
		\sup_{a\in\R}
		\left\|\int_{-\infty}^a
		K(t,\tau)G(\tau)\dd\tau\right\|_Y
		\right\|_{L^q(\R)}
		\lesssim_{p,q} A
		\norm{G}{L^p(\R;X)}.
	\end{equation*}
	In particular,
	\begin{equation*}
		\left\|\int_{-\infty}^t
		K(t,\tau)G(\tau)\dd\tau
		\right\|_{L^q(\R;Y)}
		\lesssim_{p,q} A\norm{G}{L^p(\R;X)}.
	\end{equation*}
	The same conclusions hold on an interval $J\subset\R$,
	with all integrations restricted to $J$,
	and with $\tau<t$ (retarded) replaced by $\tau>t$ (advanced).
\end{proposition}

Combining Proposition~\ref{prop:ck} with the dual
of \eqref{eq:strichartz}, we transfer estimates
for the free evolution to the corresponding
retarded and advanced Duhamel integrals.
\begin{lemma}\label{lem:transfer}
	Let $J\subset\R$ be an interval, let $Y$ be a Banach
	space, and let $M(D)$ be an $L^2(\R^2)$-bounded
	Fourier multiplier.
	Suppose that $\mathcal A(t)$ is a measurable family
	of linear maps for which the expressions below
	are well-defined, and assume that
	\begin{equation}\label{eq:free-test}
		\norm{\mathcal A(t)S(t)M(D)f}{L^r(J;Y)}
		\le A\norm{f}{L^2(\R^2)}.
	\end{equation}
	Let $p'<r\le\infty$ and $(p,q)$ admissible, then we have 
	\begin{equation}\label{eq:transfer-retarded}
		\left\|
		\mathcal A(t)
		\int_{J\cap(-\infty,t)}
		S(t-\tau)M(D)G(\tau)\dd\tau
		\right\|_{L^r(J;Y)}
		\lesssim_{p,q,r}
		A\norm{G}{L^{p'}(J;L^{q'}(\R^2))}.
	\end{equation}
	The same estimate holds with $J\cap(-\infty,t)$
	replaced by $J\cap(t,\infty)$.
	In particular, taking $(p,q)=(4,4)$ gives
	\begin{equation*}
		\left\|
		\mathcal A(t)
		\int_{J\cap(-\infty,t)}
		S(t-\tau)M(D)G(\tau)\dd\tau
		\right\|_{L^r(J;Y)}
		\lesssim_r
		A\norm{G}{L^{4/3}(J\times\R^2)},
		\qquad r>4/3.
	\end{equation*}
\end{lemma}

\begin{proof}
	We first estimate the integral over all of $J$.
	The dual of \eqref{eq:strichartz}, applied to
	the admissible pair $(p,q)$, gives
	\begin{equation*}
		\left\|\int_J S(-\tau)G(\tau)\dd\tau
		\right\|_{L^2(\R^2)}
		\lesssim_{p,q}
		\norm{G}{L^{p'}(J;L^{q'}(\R^2))}.
	\end{equation*}
	Since $M(D)$ commutes with $S(t)$,
	\begin{equation*}
		\int_J S(t-\tau)M(D)G(\tau)\dd\tau
		=
		S(t)M(D)\int_J S(-\tau)G(\tau)\dd\tau.
	\end{equation*}
	Applying \eqref{eq:free-test} therefore yields
	\begin{equation*}
		\left\|
		\mathcal A(t)\int_J
		S(t-\tau)M(D)G(\tau)\dd\tau
		\right\|_{L^r(J;Y)}
		\lesssim_{p,q}
		A\norm{G}{L^{p'}(J;L^{q'}(\R^2))}.
	\end{equation*}
	Since $p'<r$, Proposition~\ref{prop:ck}, with
	input space $L^{q'}(\R^2)$ and time exponents
	$p'$ and $r$, gives \eqref{eq:transfer-retarded}
	and the corresponding estimate with $\tau>t$.
	
\end{proof}

\vspace{7pt}
\subsection{Multipliers and local theory}
We use the classical fractional Leibniz rule
(Kato--Ponce inequality), in the form stated in
Grafakos--Oh~\cite[Theorem 1]{GO}:
if $s>0$, $1<r,p_j,q_j<\infty$, and
$1/r=1/p_j+1/q_j$ for $j=1,2$, then
\begin{equation}\label{eq:leibniz}
	\norm{\D^s(fg)}{L^r(\R^2)}
	\lesssim_{s,r,p_1,q_1,p_2,q_2}
	\norm{\D^sf}{L^{p_1}(\R^2)}
	\norm{g}{L^{q_1}(\R^2)}
	+
	\norm{f}{L^{p_2}(\R^2)}
	\norm{\D^sg}{L^{q_2}(\R^2)}.
\end{equation}
Applying \eqref{eq:leibniz} twice and then
H\"older's inequality in time gives
\begin{equation}\label{eq:cubic-leibniz}
	\norm{\D^s(|u|^2u)}{L^{4/3}(J\times\R^2)}
	\lesssim_s
	\norm{u}{L^4(J\times\R^2)}^2
	\norm{\D^su}{L^4(J\times\R^2)}.
\end{equation}

We also recall Bernstein and the two-dimensional interpolation inequality in the following forms:
\begin{equation*}
 \supp\widehat v\subset\{ |\xi|\le L\}\quad\Longrightarrow\quad
 \norm v{L^q}\lesssim_{p,q}L^{2/p-2/q}\norm v{L^p},\quad1\le p\le q\le\infty;
\end{equation*}
\begin{equation}\label{eq:sobolev-gn}
  \norm v{L^4}^4\lesssim\norm v2^2\norm{\nabla v}2^2.
\end{equation}

The product estimates below will also be used with
smooth frequency projections.
We therefore record how to insert these projections
using convolution kernels.

Let $L>0$, and let $m_L$ be a smooth multiplier
supported in $\{|\xi|\lesssim L\}$ such that
\begin{equation*}
	|\partial^\alpha m_L(\xi)|
	\lesssim_\alpha L^{-|\alpha|},
	\qquad |\alpha|\le4,
\end{equation*}
with constants independent of $L$.
We claim that its convolution kernel $k_L$ satisfies
\begin{equation*}
	\norm{k_L}{L^1(\R^2)}\lesssim1.
\end{equation*}
Indeed, rescaling $\xi=L\eta$ and integrating by
parts with $(1-\Delta_\eta)^2$ gives
\begin{equation*}
	|k_L(x)|
	\lesssim L^2(1+L^2|x|^2)^{-2}.
\end{equation*}
Integration in $x$ proves the claim.
The same conclusion holds for smooth annular
multipliers satisfying these derivative bounds,
and for the normalized derivative multipliers
$L^{-1}\xi_jm_L(\xi)$, $j=1,2$.

To apply this kernel bound to products, we use
Minkowski's inequality: for $k\in L^1(\R^2)$,
\begin{equation}\label{eq:translation-kernel}
	\norm{(k*v)w}{L^2(J\times\R^2)}
	\le
	\norm{k}{L^1(\R^2)}
	\sup_{y\in\R^2}
	\norm{v(t,x-y)w(t,x)}{L^2(J\times\R^2)}.
\end{equation}
Thus a bound for the translated products on the
right-hand side of \eqref{eq:translation-kernel},
uniform in $y$, also controls $(m_L(D)v)w$
with a constant independent of $L$.
When several factors are convolved, we apply
\eqref{eq:translation-kernel} successively,
using a product estimate uniform under spatial
translations of those factors.

We first apply \eqref{eq:inhomogeneous} to the $L^2$ local theory
of \eqref{eq:nls}, following the standard argument in
\cite[Section~3.3, Proposition~3.17]{Tao}. We then use
\eqref{eq:cubic-leibniz} to show that a uniform $L^4$ spacetime
bound for smooth solutions on a given interval suffices to pass
to $H^s$ data on that interval.

\begin{lemma}\label{lem:local}
	Let $J\subset\R$ be an interval containing $a$, and let
	$f\in L^2(\R^2)$. There is an absolute $\varepsilon_*>0$
	such that, if
	\begin{equation}\label{eq:local-small}
		\norm{S(t-a)f}{L^4(J\times\R^2)}\le\varepsilon_*,
	\end{equation}
	then \eqref{eq:nls} with $u(a)=f$ has a unique solution in
	$C(J;L^2(\R^2))\cap L^4(J\times\R^2)$. 
	
	Let now $s\ge0$, let $J$ be a compact interval containing $a$,
	and let $u,\widetilde u$ be smooth solutions on $J$ with
	$u(a)=f$ and $\widetilde u(a)=\widetilde f$. If
	\begin{equation*}
		\max\{\norm f{H^s(\R^2)},\norm{\widetilde f}{H^s(\R^2)}\}
		\le R,\qquad
		\max\{\norm u{L^4(J\times\R^2)},
		\norm{\widetilde u}{L^4(J\times\R^2)}\}\le A,
	\end{equation*}
	then
	\begin{equation}\label{eq:positive-stability}
		\norm{u-\widetilde u}{C(J;H^s(\R^2))}
		+\norm{u-\widetilde u}{L^4(J;W^{s,4}(\R^2))}
		\lesssim_{s,R,A}\norm{f-\widetilde f}{H^s(\R^2)} .
	\end{equation}
	Consequently, for any family of smooth data satisfying these
	bounds, the solution map extends to its $H^s(\R^2)$ closure.
	For each datum in the closure, the extension gives the unique
	solution in $C(J;H^s(\R^2))\cap L^4(J;W^{s,4}(\R^2))$.
\end{lemma}

\begin{proof}
	\textit{Local solutions in $L^2$.}
	Consider the solution map
	\begin{equation}\label{eq:duhamel}
		\Phi(u)(t):=S(t-a)f
		-\ii\int_a^tS(t-\tau)(|u|^2u)(\tau)\dd\tau.
	\end{equation}
	By \eqref{eq:inhomogeneous}, the nonlinear part
	of \eqref{eq:duhamel} in
	$L^\infty(J;L^2(\R^2))\cap L^4(J\times\R^2)$
	is bounded by
	$\lesssim\norm{u}{L^4(J\times\R^2)}^3$.
	Moreover,
	\begin{equation*}
		\norm{|u|^2u-|v|^2v}{L^{4/3}(J\times\R^2)}
		\lesssim
		\bigl(
		\norm{u}{L^4(J\times\R^2)}^2
		+
		\norm{v}{L^4(J\times\R^2)}^2
		\bigr)
		\norm{u-v}{L^4(J\times\R^2)}.
	\end{equation*}
	Thus $\Phi$ is a contraction on the closed ball
	\begin{equation*}
		\norm{u}{L^4(J\times\R^2)}\le2\varepsilon_*,
		\qquad
		\norm{u-S(t-a)f}{L^\infty(J;L^2(\R^2))}
		\le\varepsilon_*,
	\end{equation*}
	equipped with the sum of the
	$L^\infty(J;L^2(\R^2))$ and
	$L^4(J\times\R^2)$ distances,
	provided $\varepsilon_*$ is sufficiently small.
	The Duhamel integral is continuous in
	$L^2(\R^2)$, so the fixed point belongs to
	$C(J;L^2(\R^2))$.
	The difference estimate and subdivision into
	intervals with small $L^4$ norms give uniqueness
	in the stated class and continuous dependence.
	
   \textit{Stability and smooth approximation in $H^s$.}
   To prove \eqref{eq:positive-stability}, we first use the assumed
   $L^4(J\times\R^2)$ bound to choose a partition on which both solutions
   are small in $L^4(I;W^{s,4}(\R^2))$.
   
   For $v\in\{u,\widetilde u\}$, partition $J$ into intervals $I$ on which
   $\norm v{L^4(I\times\R^2)}\le \eta\ll 1$.  Let $c_I$ be the left endpoint of $I$. 
   By \eqref{eq:strichartz}, \eqref{eq:inhomogeneous}, and
   \eqref{eq:cubic-leibniz},
   \begin{equation*}
   	\norm v{C(I;H^s(\R^2))}+\norm v{L^4(I;W^{s,4}(\R^2))}
   	\lesssim_s \norm{v(c_I)}{H^s(\R^2)}
   	+\eta^2\norm v{L^4(I;W^{s,4}(\R^2))}.
   \end{equation*}
   Absorbing the last term, iterating from
   $\norm{v(a)}{H^s(\R^2)}\le R$, and summing over
   $I$, we can obtain
   \begin{equation*}
   	\norm u{L^4(J;W^{s,4}(\R^2))}
   	+\norm{\widetilde u}{L^4(J;W^{s,4}(\R^2))}
   	\lesssim_{s,A} R .
   \end{equation*}
  
  We now choose a common partition for $u$ and $\widetilde u$ into intervals $K$,
  such that
  \begin{equation*}
  	\norm u{L^4(K;W^{s,4}(\R^2))}^4
  	+\norm{\widetilde u}{L^4(K;W^{s,4}(\R^2))}^4
  	\le\delta^4.
  \end{equation*}
  
  Let $c_K$ be the left endpoint of $K$. On $K$, the
  difference $w=u-\widetilde u$ satisfies
  \begin{equation*}
  	w(t)=S(t-c_K)w(c_K)
  	-\ii\int_{c_K}^t S(t-\tau)
  	\bigl(|u|^2u-|\widetilde u|^2\widetilde u\bigr)(\tau)\dd\tau .
  \end{equation*}
  For the free term, unitarity and \eqref{eq:strichartz} give
  \begin{equation*}
  	\norm{S(t-c_K)w(c_K)}{C(K;H^s(\R^2))}
  	+\norm{S(t-c_K)w(c_K)}{L^4(K;W^{s,4}(\R^2))}
  	\lesssim \norm{w(c_K)}{H^s(\R^2)} .
  \end{equation*}
  For the Duhamel integral, \eqref{eq:inhomogeneous} gives, for
  $(q,r)\in\{(\infty,2),(4,4)\}$,
  \begin{equation*}
  	\left\|\int_{c_K}^t S(t-\tau)
  	\bigl(|u|^2u-|\widetilde u|^2\widetilde u\bigr)(\tau)\dd\tau
  	\right\|_{L^q(K;W^{s,r}(\R^2))}
  	\lesssim
  	\norm{|u|^2u-|\widetilde u|^2\widetilde u}
  	{L^{4/3}(K;W^{s,4/3}(\R^2))}.
  \end{equation*}
  \begin{equation*}
  	\lesssim
  	\bigl(\norm u{L^4(K;W^{s,4}(\R^2))}^2
  	+\norm{\widetilde u}{L^4(K;W^{s,4}(\R^2))}^2\bigr)
  	\norm w{L^4(K;W^{s,4}(\R^2))}
  	\lesssim \delta^2\norm w{L^4(K;W^{s,4}(\R^2))}.
  \end{equation*}
  Choose $\delta$
   small so that the coefficient of
  $\norm w{L^4(K;W^{s,4}(\R^2))}$ is at most $1/2$. Absorbing this term gives
  \begin{equation*}
  	\norm w{C(K;H^s(\R^2))}+\norm w{L^4(K;W^{s,4}(\R^2))}
  	\lesssim_s \norm{w(c_K)}{H^s(\R^2)} .
  \end{equation*}
  Iterating it again proves
  \eqref{eq:positive-stability}.

\end{proof}

\vspace{10pt}
\section{The modified energy and a mass--momentum identity}\label{sec:energy}
Write $F(z):=|z|^2z$ and
$M:=1+\norm{u(0)}{L^2(\R^2)}$. Fix $0<\varepsilon<1/24$ and a
sufficiently small absolute constant $a>0$. For $N$ sufficiently
large, set $L=N^{1-\varepsilon}$ and let $J=[0,t_1]$. We use the
modified energy
\begin{equation}\label{energy}
	\mathcal E_N(u)
	:=
	\frac12\norm{\nabla P_Nu}{L^2(\R^2)}^2
	+\frac14\norm{P_Nu}{L^4(\R^2)}^4
	+\frac{N^2}{2}\norm{Q_Nu}{L^2(\R^2)}^2 .
\end{equation}
\begin{remark}
We now explain the choice of the modified energy in \eqref{energy}.
	The first two terms of \eqref{energy} give the low-frequency energy,
	\begin{equation*}
		\frac12\norm{\nabla P_Nu}{L^2}^2
		+\frac14\norm{P_Nu}{L^4}^4
		=E(P_Nu).
	\end{equation*}
	We therefore expect the energy cancellation to give better control
	of the low-frequency part of the solution.
	
	For the first and third terms, recall the $I$-operator
	of the $I$-method~\cite{CKSTT} with $s=0$:
	\begin{equation*}
		\widehat{I_Nu}(\xi)=m_N(\xi)\widehat u(\xi),\qquad m_N(\xi)=
		\left\{
		\begin{array}{ll}
			1, & |\xi|\le N,\\[2pt]
			\displaystyle\frac{N}{|\xi|}, & |\xi|\ge 2N.
		\end{array}
		\right.
	\end{equation*}
	Since $m_N(\xi)\sim1$ for $N\le|\xi|\le2N$, Plancherel gives
	\begin{equation*}
		\frac12\norm{\nabla P_Nu}{L^2}^2
		+\frac{N^2}{2}\norm{Q_Nu}{L^2}^2
		\sim \frac12\norm{I_Nu}{\dot H^1}^2 .
	\end{equation*}
\end{remark}

\vspace{7pt}
\subsection{Framework of the energy bootstrapping}
We first estimate the change in $\mathcal E_N(u)$ on $J=[0,t_1]$.
Applying $P_N$ and $Q_N$ to \eqref{eq:nls} gives
\begin{equation*}
	\ii\partial_tP_Nu=-\Delta P_Nu+F(P_Nu)+R,
	\qquad R:=P_NF(u)-F(P_Nu),
\end{equation*}
and
\begin{equation*}
	\ii\partial_tQ_Nu=-\Delta Q_Nu+Q_NF(u).
\end{equation*}
 Hence
$\partial_t\mathcal E_N(u)=J_1+J_2+J_3$, where
\begin{equation}\label{energy-kin}
	J_1:=\im\int_{\R^2}\overline{-\Delta P_Nu}\,R\dd x ,
\end{equation}
\begin{equation}\label{energy-pot}
	J_2:=\im\int_{\R^2}\overline{F(P_Nu)}\,R\dd x ,
\end{equation}
\begin{equation}\label{energy-tail}
	J_3:=N^2\im\int_{\R^2}\overline{Q_Nu}\,Q_NF(u)\dd x .
\end{equation}

To see which estimates are needed, consider the frequencies in
\eqref{energy-kin}--\eqref{energy-tail}. In \eqref{energy-kin}
and \eqref{energy-pot}, the multiplier of $R$ vanishes if its
three input frequencies all have magnitude at most $N/3$.
In \eqref{energy-tail}, the factor $Q_Nu$ already has frequency
at least $N$. Since all the frequencies sum to zero
in each spatial integral, there is another frequency
$\gtrsim N$. Thus each quartic integral contains two high-frequency
factors. For the sextic integral \eqref{energy-pot}, we shall
retain two such factors and place two factors of frequency
$\lesssim N$ in $L^\infty(J\times\R^2)$. This leaves four factors
to estimate in each case. For more details, we refer to the proof of Lemma \ref{source}.

If two of the four factors are $P_Lu$, we pair each of them
with a high-frequency factor in $L^2(J\times\R^2)$. Expanding
the convolution kernels of the projections $P_N$ and $Q_N$ introduces spatial
translations (see (\ref{eq:translation-kernel})), so for any $h$ we put
\begin{equation}\label{Y}
	Y(h):=\sup_{y\in\R^2}
	\norm{(P_Lu)(t,x+y)h(t,x)}{L^2(J\times\R^2)} .
\end{equation}
If only one factor is $P_Lu$, we use one such $L^2$ product
and put the other two factors in $L^4$. If all four factors
come from $Q_Lu$, we put each in $L^4$. Thus the energy estimate
requires control of both $Y(Q_{aN}u)$ and
$\norm{Q_Lu}{L^4(J\times\R^2)}$.

Suppose $\mathcal E_N(u(0))\le 1/4$. We use a bootstrap argument to prove that,
for every $J=[0,t_1]$ in the time range under consideration,
\begin{equation}\label{bootstrap}
	\sup_{t\in J}\mathcal E_N(u(t))\le 1,\qquad
	\norm{Q_Lu}{L^4(J\times\R^2)}\le N^{-2/3},\qquad
	Y(Q_{aN}u)\le N^{-1-\varepsilon/4}.
\end{equation}
At $t_1=0$, the first quantity is at most $1/4$ and the other two
vanish. Moreover, all three quantities depend continuously on $t_1$. 

It therefore suffices to rule out a first time $T>0$ at which one
of the three quantities in \eqref{bootstrap} reaches its stated
upper bound. Suppose such a time exists and take $J=[0,T]$.
By the choice of $T$, all three inequalities in \eqref{bootstrap}
hold on this particular interval, with equality in at least one.
We may therefore apply Lemma~\ref{source} on $[0,T]$ to obtain
\begin{equation*}
	\int_0^T |\partial_t\mathcal E_N(u(t))|\dd t
	\lesssim_M N^{-\varepsilon/2}.
\end{equation*}
For $N$ sufficiently large, this gives
\begin{equation*}
	\sup_{0\le t\le T}\mathcal E_N(u(t))
	\le \mathcal E_N(u(0))
	+\int_0^T|\partial_t\mathcal E_N(u(t))|\dd t
	\le \frac14+\frac14=\frac12 .
\end{equation*}
On the same interval, Proposition~\ref{tail} gives
\begin{equation*}
	\norm{Q_Lu}{L^4([0,T]\times\R^2)}
	\lesssim_{M}N^{-3/4+2\varepsilon}<N^{-2/3}.
\end{equation*}
For the third quantity, we verify \eqref{operator-small} and
apply Proposition~\ref{high} to obtain
\begin{equation*}
	Y(Q_{aN}u)\lesssim_{M,a}N^{-1-\varepsilon/2}
	<N^{-1-\varepsilon/4}.
\end{equation*}
 Thus none of the three quantities can reach
its upper bound at $T$. This rules out a first such time and proves
\eqref{bootstrap} throughout that time range. We begin with the
energy variation estimate used above, where (\ref{eq:translation-kernel}) will be implicitly and repeatedly used.

\begin{lemma}\label{source}
	Assume \eqref{bootstrap} and put
	\begin{equation*}
		D_4:=
		Y(Q_{aN}u)^2
		+Y(Q_{aN}u)\norm{Q_Lu}{L^4(J\times\R^2)}^2
		+\norm{Q_Lu}{L^4(J\times\R^2)}^4 .
	\end{equation*}
	Then
	\begin{equation}\label{variation}
		\int_J|\partial_t\mathcal E_N(u(t))|\dd t
		\lesssim_M N^2D_4
		\lesssim N^{-\varepsilon/2}.
	\end{equation}
	Moreover, if $q\in L^\infty(J\times\R^2)$ has spatial Fourier
	support in $\{|\xi|\le L\}$ for almost every $t\in J$, then
	\begin{equation}\label{qs1}
		\int_J\left|\int_{\R^2}q\,\overline{P_Nu}\,R\dd x\right|\dd t
		\lesssim \norm q{L^\infty(J\times\R^2)}D_4
	\end{equation}
	and, for $j=1,2$,
	\begin{equation}\label{qs2}
		\int_J\left(
		\left|\int_{\R^2}q\,\partial_j\overline{P_Nu}\,R\dd x\right|
		+\left|\int_{\R^2}q\,\overline{P_Nu}\,\partial_jR\dd x\right|
		\right)\dd t
		\lesssim N\norm q{L^\infty(J\times\R^2)}D_4 .
	\end{equation}
	The implicit constants may also depend on the fixed cutoffs.
\end{lemma}

\begin{proof}
	\textit{Why two high-frequency factors are present.}
	The first high frequency follows from the definition of $R$.
	Indeed, its multiplier is
	\begin{equation*}
		\chi\big((\xi_1-\xi_2+\xi_3)/N\big)
		-\chi(\xi_1/N)\chi(\xi_2/N)\chi(\xi_3/N).
	\end{equation*}
	If $|\xi_1|,|\xi_2|,|\xi_3|\le N/3$, then
	$|\xi_1-\xi_2+\xi_3|\le N$, and all four cutoffs equal one.
	Thus a nonzero contribution from $R$ requires an input of
	frequency greater than $N/3$. In \eqref{energy-tail}, the factor
	$\overline{Q_Nu}$ already has frequency at least $N$.
	
	The second high frequency follows from integration in $x$.
	To see this, decompose each input as
	\begin{equation}\label{eq:source-decomposition}
		u=P_Lu+P_{2aN}Q_Lu+Q_{2aN}u .
	\end{equation}
	 The first two summands in
	\eqref{eq:source-decomposition} have Fourier support in
	$\{|\xi|\le4aN\}$. Fix $t\in J$. In the spatial integrals in \eqref{qs1}
	and \eqref{qs2}, let $\eta$ be the frequency of $q(t,\cdot)$
	and let $\xi_1,\ldots,\xi_4$ be the frequencies of the four
	inputs before conjugation. Integration in $x$ requires
	\begin{equation*}
		\eta+\xi_1-\xi_2+\xi_3-\xi_4=0 .
	\end{equation*}
	Since $|\eta|\le L$, we obtain
	$|\xi_1-\xi_2+\xi_3-\xi_4|\le L$.
	The spatial integrals in \eqref{energy-kin}--\eqref{energy-tail}
	have weight $1$, so the signed sum is exactly zero: in
	\eqref{energy-kin} and \eqref{energy-tail} it is the sum of four
	input frequencies, and in \eqref{energy-pot} it is the sum of six.
	If only one input had frequency greater than $4aN$, its frequency
	$\xi$ would therefore satisfy
	\begin{equation*}
		|\xi|\le L+5(4aN)<N/3
	\end{equation*}
	for $a$ sufficiently small and $N$ sufficiently large. This
	contradicts the frequency restriction just proved. Hence every
	nonzero integral in \eqref{energy-kin}--\eqref{energy-tail} and
	\eqref{qs1}--\eqref{qs2} contains at least two factors from
	$Q_{2aN}u$. We use this restriction before taking absolute values.
	
	\textit{Estimating four factors.}
	 We next bound the product of $P_Lu$ and $Q_{2aN}u$ by
	$Y(Q_{aN}u)$. By definition,
	\begin{equation*}
		Q_{2aN}u
		=Q_{aN}u-P_{2aN}Q_{aN}u .
	\end{equation*}
	
	To estimate both terms, observe that for every $y,z\in\R^2$,
	translation in $x$ and \eqref{Y} give
	\begin{equation}\label{eq:translated-high-low-pair}
		\norm{(P_Lu)(t,x+y)(Q_{aN}u)(t,x+z)}{L^2(J\times\R^2)}
		=
		\norm{(P_Lu)(t,x+y-z)(Q_{aN}u)(t,x)}{L^2(J\times\R^2)}
		\le Y(Q_{aN}u).
	\end{equation}
	Let $k_{2aN}$ be the convolution kernel of $P_{2aN}$. By
	Minkowski's inequality and \eqref{eq:translated-high-low-pair}
	with $z$ replaced by $z-\eta$,
	\begin{equation*}
		\norm{(P_Lu)(t,x+y)(P_{2aN}Q_{aN}u)(t,x+z)}{L^2(J\times\R^2)}
		\le \int_{\R^2}|k_{2aN}(\eta)|\,Y(Q_{aN}u)\dd\eta .
	\end{equation*}
	Since $\norm{k_{2aN}}{L^1(\R^2)}=\norm{k_1}{L^1(\R^2)}$, the
	triangle inequality now yields
	\begin{equation*}
		\sup_{y,z\in\R^2}
		\norm{(P_Lu)(t,x+y)(Q_{2aN}u)(t,x+z)}{L^2(J\times\R^2)}
		\le \bigl(1+\norm{k_1}{L^1(\R^2)}\bigr)Y(Q_{aN}u).
	\end{equation*}
	
    We now apply these estimates to \eqref{qs1}. Expand the four inputs
    using \eqref{eq:source-decomposition}, discard the terms that vanish
    by the frequency restriction above, and then expand the projection
    kernels. For fixed spatial translations $y_1,\ldots,y_4$, denote
    the four translated factors by $b_1,\ldots,b_4$, in any order and
    up to complex conjugation. At least two come from $Q_{2aN}u$, so
    we may take
    \begin{equation*}
    	b_i(t,x)=(Q_{2aN}u)(t,x+y_i),\qquad i=1,2.
    \end{equation*}
    For each choice of $b_3,b_4$, we have
    \begin{equation*}
    	\int_J\left|\int_{\R^2}q\,b_1b_2b_3b_4\dd x\right|\dd t
    	\le \norm q{L^\infty(J\times\R^2)}
    	\norm{b_1b_2b_3b_4}{L^1(J\times\R^2)} .
    \end{equation*}
    
    If $b_3,b_4$ both come from $P_Lu$, pair $b_3$ with $b_1$ and
    $b_4$ with $b_2$. The translated product estimate gives
    \begin{equation*}
    	\norm{b_1b_2b_3b_4}{L^1(J\times\R^2)}
    	\le \norm{b_1b_3}{L^2(J\times\R^2)}
    	\norm{b_2b_4}{L^2(J\times\R^2)}
    	\lesssim Y(Q_{aN}u)^2 .
    \end{equation*}
  If $b_3$ is the only factor from $P_Lu$, pair $b_3$ with $b_1$
  as above. Since $b_4$ comes from either $P_{2aN}Q_Lu$ or $Q_{2aN}u$,
  the $L^4$ boundedness of $P_{2aN}$ and $Q_{2aN}$ gives
  \begin{equation*}
  	\norm{b_2}{L^4(J\times\R^2)}+\norm{b_4}{L^4(J\times\R^2)}
  	\lesssim \norm{Q_Lu}{L^4(J\times\R^2)} .
  \end{equation*}
  By H\"older's inequality,
  \begin{equation*}
  	\norm{b_1b_2b_3b_4}{L^1(J\times\R^2)}
  	\le \norm{b_1b_3}{L^2(J\times\R^2)}
  	\norm{b_2}{L^4(J\times\R^2)}
  	\norm{b_4}{L^4(J\times\R^2)}
  	\lesssim Y(Q_{aN}u)\norm{Q_Lu}{L^4(J\times\R^2)}^2 .
  \end{equation*}
    If neither $b_3$ nor $b_4$ comes from $P_Lu$, then
    \begin{equation*}
    	\norm{b_1b_2b_3b_4}{L^1(J\times\R^2)}
    	\le \prod_{i=1}^4\norm{b_i}{L^4(J\times\R^2)}
    	\lesssim \norm{Q_Lu}{L^4(J\times\R^2)}^4 .
    \end{equation*}
    All the bounds above are uniform in $y_1,\ldots,y_4$. Integrating
    against the projection kernel of $P_N$ and summing gives
    \begin{equation*}
    	\int_J\left|\int_{\R^2}q\,\overline{P_Nu}\,R\dd x\right|\dd t
    	\lesssim \norm q{L^\infty(J\times\R^2)}D_4 ,
    \end{equation*}
    which proves \eqref{qs1}. Similarly, we can obtain \eqref{qs2}.
	
	\textit{Estimating the energy derivative.}
	For \eqref{energy-kin}, the kernel of $-N^{-2}\Delta P_N$ has
	uniformly bounded $L^1$ norm. Following the same arguments as in \eqref{qs1}-\eqref{qs2}, 
	\begin{equation*}
		\int_J|J_1(t)|\dd t\lesssim N^2D_4 .
	\end{equation*}
	
	For \eqref{energy-pot}, we claim that both factors estimated in
	$L^\infty$ can be chosen from $\overline{F(P_Nu)}$. To see this,
	apply the decomposition \eqref{eq:source-decomposition} to all six inputs. Also recall that, in every
	nonzero summand, at least one of the three inputs to $R$ comes
	from $Q_{2aN}u$.
	
	Choose one such $Q_{2aN}u$ factor inside $R$. By the frequency
	restriction proved above, we can choose a second $Q_{2aN}u$ factor
	among the other five inputs. At most one of the two chosen factors
	belongs to $\overline{F(P_Nu)}$. Thus at least two of its three
	factors are still available for the $L^\infty$ estimate.
	
	Each of these two factors has the form $P_Nb$ or its conjugate,
	where $b$ is one of the three summands in
	\eqref{eq:source-decomposition}. Since $P_Nb$ has Fourier support
	in $|\xi|\le2N$, Bernstein's inequality and mass conservation give
	\begin{equation*}
		\norm{P_Nb}{L^\infty(J\times\R^2)}
		\lesssim N\sup_{t\in J}\norm{b(t)}{L^2(\R^2)}
		\lesssim N\norm{u(0)}{L^2(\R^2)} .
	\end{equation*}
	After placing these two factors in $L^\infty$, the four remaining
	factors still contain both chosen $Q_{2aN}u$ factors. The quartic
	product bound proved before therefore gives
	\begin{equation*}
		\int_J |J_2(t)|\dd t
		\lesssim
		\bigl(N\norm{u(0)}{L^2(\R^2)}\bigr)^2 D_4
		\lesssim_M N^2D_4 .
	\end{equation*}
	
	For \eqref{energy-tail}, the four factors already include two
	from $Q_{2aN}u$. Then one can similarly derive
	\begin{equation*}
		\int_J|J_3(t)|\dd t\lesssim N^2D_4 .
	\end{equation*}
	Adding the three estimates and substituting \eqref{bootstrap}, we obtain \eqref{variation}.
	
\end{proof}

\vspace{7pt}
\subsection{The mass--momentum identity}

Lemma~\ref{source} controls the variation of $\mathcal E_N(u)$.
For the two spacetime bounds in \eqref{bootstrap}, we shall need
two kinds of product estimates. In the proof of the bound for
$Y(Q_{aN}u)$, Proposition~\ref{high} pairs $P_Lu$ with the free
term in the high-frequency Duhamel formula and, later, with a
component of $P_Nu$ at frequencies comparable to $aN$. These
products are estimated in \eqref{freehigh} and \eqref{shell},
respectively. In the proof of the bound for
$\norm{Q_Lu}{L^4(J\times\R^2)}$, Proposition~\ref{tail} uses
Lemma~\ref{quartet} to estimate the Duhamel terms at frequencies
between $L$ and $N$. That lemma requires spatial gradients of
two products: one formed from two components of $P_Nu$, and one
formed from a component of $P_Nu$ and a test wave. We establish
the corresponding estimates in \eqref{vv} and \eqref{vf}.

The connection between these estimates is that, for separated
frequencies, a spatial gradient also controls the product itself.
To see this for the free term, let $h(t)=S(t-t_0)f_H$, where
$H\gtrsim N$ and $\supp\widehat f_H\subset\{H/2\le|\xi|\le2H\}$.
For fixed $t$ and $z\in\R^2$, the spatial Fourier support of
$(P_Lu)(t,x)\overline{h(t,x-z)}$ lies in
$|\xi|\ge H/2-2L\gtrsim H$ for $N$ sufficiently large. Thus
Plancherel's identity gives
\begin{equation}\label{as}
	\norm{(P_Lu)(t,x)\overline{h(t,x-z)}}{L^2_x(\R^2)}
	\lesssim \frac{1}{H}
	\norm{\nabla_x\bigl[(P_Lu)(t,x)\overline{h(t,x-z)}\bigr]}
	{L^2_x(\R^2)} .
\end{equation}

This is the spatial derivative that appears in the bilinear virial
identity of Planchon--Vega~\cite[Theorem~2.3]{PV}. Indeed, they consider
the interaction of the mass densities of two Schr\"odinger waves
against a convex weight depending on their separation:
\begin{equation}\label{eq:radial-interaction}
	\iint_{\R^2\times\R^2}
	A(x-y)|U(t,x)|^2|V(t,y)|^2\dd x\dd y .
\end{equation}
To see how this leads to the right-hand side of (\ref{as}), write
$\rho_U:=|U|^2$ and $p_U:=\im(\bar U\nabla U)$, and likewise for
$V$. If $U$ and $V$ solve the free Schr\"odinger equation,
differentiating \eqref{eq:radial-interaction} and integrating by
parts gives
\begin{equation*}
	\frac{\dd}{\dd t}\iint A(x-y)\rho_U(t,x)\rho_V(t,y)\dd x\dd y
	=2\iint\nabla A(x-y)\cdot
	\bigl[p_U(t,x)\rho_V(t,y)-\rho_U(t,x)p_V(t,y)\bigr]\dd x\dd y .
\end{equation*}
Differentiating it again, we have
\begin{equation*}
	\frac{\dd^2}{\dd t^2}
	\iint_{\R^2\times\R^2} A(x-y)\rho_U(t,x)\rho_V(t,y)\dd x\dd y
	=
	4\sum_{j,k=1}^2\iint_{\R^2\times\R^2}
	\partial_{jk}A(x-y)\,
	\re\bigl(\overline{B_j(t,x,y)}B_k(t,x,y)\bigr)\dd x\dd y ,
\end{equation*}
where
\begin{equation*}
	B(t,x,y)
	:=(\nabla_x+\nabla_y)\bigl[U(t,x)\overline{V(t,y)}\bigr]
	=\overline{V(t,y)}\,\nabla U(t,x)
	+U(t,x)\,\overline{\nabla V(t,y)} .
\end{equation*}
Since $A$ is convex, the quadratic form on the right-hand side
is nonnegative. For fixed $z\in\R^2$, setting $y=x-z$ and applying
the chain rule gives
\begin{equation*}
	B(t,x,x-z)
	=\nabla_x\bigl[U(t,x)\overline{V(t,x-z)}\bigr],
\end{equation*}
which has the same form as the gradient on the right-hand side
of \eqref{as}.

In our application the first wave is $U=P_Nu$. We use this field
before projecting to $P_Lu$ because its equation has forcing
$R=P_NF(u)-F(P_Nu)$, for which Lemma~\ref{source} gives the
weighted estimates \eqref{qs1}--\eqref{qs2}. The second wave is
either another copy of $P_Nu$ or a solution of a linear
Schr\"odinger equation, possibly with forcing. We therefore
choose the mass--momentum expression directly:
\begin{equation}\label{eq:radial-action}
	\mathcal M(t):=\iint\nabla A(x-y)\cdot
	\bigl[p_U(t,x)\rho_V(t,y)-\rho_U(t,x)p_V(t,y)\bigr]\dd x\dd y .
\end{equation}

Fix a nonnegative radial Schwartz function $\varphi$ on $\R^2$
such that $\int_{\R^2}\varphi=1$, $\varphi(0)>0$, and
$\supp\widehat\varphi\subset\{|\xi|\le1\}$. In \eqref{eq:radial-action}, take
\begin{equation}\label{weight}
	A(z):=\int_{\R^2}\sqrt{L^{-2}+|z-y|^2}\,L^2\varphi(Ly)\dd y .
\end{equation}

We first differentiate \eqref{eq:radial-action}, keeping the
forcing terms for both waves. This gives the form of the bilinear
virial identity of Planchon--Vega~\cite[Theorem~2.3]{PV} that we
will use for both free and forced test waves.

\begin{lemma}\label{lem:forced-radial-identity}
	Let $U,V$ be smooth, spatially decaying solutions on $J$ of
	\begin{equation*}
		(\ii\partial_t+\Delta)U=\sigma_U|U|^2U+R_U,
		\qquad
		(\ii\partial_t+\Delta)V=\sigma_V|V|^2V+R_V,
		\qquad \sigma_U,\sigma_V\in\{0,1\}.
	\end{equation*}
	Let $\mathcal M$ be given by \eqref{eq:radial-action}, with $A$
	as in \eqref{weight}, and put
	\begin{equation*}
		B(t,x,y)=\overline{V(t,y)}\,\nabla U(t,x)
		+U(t,x)\nabla\overline{V(t,y)} .
	\end{equation*}
	Then we can obtain
	\begin{equation*}
		\frac{d\mathcal M}{dt}=I_1+I_2+I_3,
	\end{equation*}
	\begin{equation}\label{principal}
		I_1
		:=2\sum_{j,k=1}^2\iint \partial_{jk}A\,
		\re(\overline{B_j}B_k)\dd x\dd y
		+\frac12\iint\Delta A\,
		\bigl(\sigma_U\rho_U^2\rho_V+\sigma_V\rho_U\rho_V^2\bigr)
		\dd x\dd y ,
	\end{equation}
	\begin{equation}\label{actualsource}
		I_2:=\iint\nabla A\cdot
		\bigl[\rho_V\re(\overline{R_U}\nabla U-\bar U\nabla R_U)
		-2p_V\im(\bar U R_U)\bigr]\dd x\dd y ,
	\end{equation}
	\begin{equation}\label{testsource}
		I_3:=\iint\nabla A\cdot
		\bigl[2p_U\im(\bar V R_V)
		-\rho_U\re(\overline{R_V}\nabla V-\bar V\nabla R_V)\bigr]
		\dd x\dd y,
	\end{equation}
	 where the $U$ factors are
	evaluated at $(t,x)$, the $V$ factors at $(t,y)$, and the
	derivatives of $A$ at $x-y$.
	
In particular, uniformly in $L>0$,
\begin{equation}\label{eq:radial-positive}
	I_1(t)\gtrsim_\varphi
	\iint L(1+L|x-y|)^{-3}|B(t,x,y)|^2\dd x\dd y .
\end{equation}
\end{lemma}

\begin{proof}
	We first record the properties of \eqref{weight} used in the
	proof. The gradient of $\sqrt{L^{-2}+|z|^2}$ has length at most
	one. Since $L^2\varphi(L\,\cdot)$ is nonnegative and has integral
	one, convolution gives
	\begin{equation*}
		\norm{\nabla A}{L^\infty(\R^2)}\le1,
		\qquad
		\supp\widehat{\nabla A}\subset\{|\xi|\le L\}.
	\end{equation*}
	For each unit vector $e\in\R^2$, we have
	\begin{equation*}
		\partial_e^2\sqrt{L^{-2}+|z|^2}
		=\frac{L^{-2}+|z|^2-(e\cdot z)^2}{(L^{-2}+|z|^2)^{3/2}}
		\ge L(1+L|z|)^{-3}.
	\end{equation*}
	Since $\varphi\ge0$ and $1+L|z-y|\le 2(1+L|z|)$ for $|y|\le L^{-1}$,
	\eqref{weight} gives
	\begin{equation*}
		\partial_e^2 A(z)
		\ge \int_{|y|\le L^{-1}}
		\frac{L}{(1+L|z-y|)^3}\,L^2\varphi(Ly)\dd y
		\ge \frac{L}{8(1+L|z|)^3}\int_{|w|\le1}\varphi(w)\dd w .
	\end{equation*}
	 Hence
	\begin{equation}\label{hessian}
		D^2A(z)\gtrsim_\varphi L(1+L|z|)^{-3}\operatorname{Id}.
	\end{equation}
	
	We next compute the local equations for $U$. Substituting
	$U_t=\ii\Delta U-\ii\sigma_U\rho_U U-\ii R_U$ into
	$\partial_t\rho_U=2\re(\bar U U_t)$ and using
	$\im(\bar U\Delta U)=\operatorname{div}p_U$ gives
	\begin{equation}\label{masslaw}
		\partial_t\rho_U+2\operatorname{div}p_U
		=2\im(\bar U R_U).
	\end{equation}
	For the momentum, differentiating
	$p_{U,j}=\im(\bar U\partial_jU)$ gives
	\begin{equation*}
		\partial_t p_{U,j}
		=\re(\bar U\partial_j\Delta U-\Delta\bar U\,\partial_jU)
		+\sigma_U\re\bigl(\rho_U\bar U\partial_jU
		-\bar U\partial_j(\rho_U U)\bigr)
		+\re(\overline{R_U}\,\partial_jU-\bar U\,\partial_jR_U).
	\end{equation*}
	The first term equals
	$\sum_{k=1}^2\partial_k
	\re(\bar U\partial_{jk}U-\partial_k\bar U\,\partial_jU)$,
	and the product rule gives
	\begin{equation*}
		\re(\bar U\partial_{jk}U-\partial_k\bar U\,\partial_jU)
		=\tfrac12\partial_{jk}\rho_U
		-2\re(\partial_j\bar U\,\partial_kU).
	\end{equation*}
	The second term equals $-\sigma_U\rho_U\partial_j\rho_U$.
	Consequently,
	\begin{equation}\label{momlaw}
		\partial_t p_{U,j}
		=\tfrac12\partial_j\Delta\rho_U
		-2\sum_{k=1}^2\partial_k\re(\partial_j\bar U\,\partial_kU)
		-\tfrac{\sigma_U}{2}\partial_j(\rho_U^2)
		+\re(\overline{R_U}\,\partial_jU-\bar U\,\partial_jR_U).
	\end{equation}
	The equations for $\rho_V,p_V$ follow by substituting
	$(V,\sigma_V,R_V)$ for $(U,\sigma_U,R_U)$.
	
	Differentiating \eqref{eq:radial-action} gives
	\begin{equation*}
		\frac{\dd\mathcal M}{\dd t}
		=\iint\nabla A\cdot
		\bigl[\rho_V\partial_t p_U+p_U\partial_t\rho_V
		-p_V\partial_t\rho_U-\rho_U\partial_t p_V\bigr]
		\dd x\dd y .
	\end{equation*}
	Substituting \eqref{masslaw} and \eqref{momlaw} for both fields
	and grouping the terms gives
	\begin{equation*}
		\frac{\dd\mathcal M}{\dd t}=J_1+J_2+J_3+J_4+J_5+J_6,
	\end{equation*}
	where
	\begin{equation}\label{eq:radial-expansion-density}
		J_1:=\frac12\iint\nabla A\cdot
		\bigl[\rho_V\nabla_x\Delta_x\rho_U
		-\rho_U\nabla_y\Delta_y\rho_V\bigr]\dd x\dd y ,
	\end{equation}
	\begin{equation}\label{eq:radial-expansion-stress}
		J_2:=-2\sum_{j,k=1}^2\iint\partial_j A
		\bigl[\rho_V\partial_{x_k}\re(\partial_j\bar U\,\partial_kU)
		-\rho_U\partial_{y_k}\re(\partial_j\bar V\,\partial_kV)\bigr]
		\dd x\dd y ,
	\end{equation}
	\begin{equation}\label{eq:radial-expansion-current}
		J_3:=2\iint\nabla A\cdot
		\bigl[p_V\operatorname{div}_x p_U
		-p_U\operatorname{div}_y p_V\bigr]\dd x\dd y ,
	\end{equation}
	\begin{equation}\label{eq:radial-expansion-pressure}
		J_4:=\frac12\iint\nabla A\cdot
		\bigl[-\sigma_U\rho_V\nabla_x(\rho_U^2)
		+\sigma_V\rho_U\nabla_y(\rho_V^2)\bigr]\dd x\dd y ,
	\end{equation}
	\begin{equation}\label{eq:radial-expansion-source-u}
		J_5:=\iint\nabla A\cdot
		\bigl[\rho_V\re(\overline{R_U}\nabla U-\bar U\nabla R_U)
		-2p_V\im(\bar U R_U)\bigr]\dd x\dd y ,
	\end{equation}
	and
	\begin{equation}\label{eq:radial-expansion-source-v}
		J_6:=\iint\nabla A\cdot
		\bigl[2p_U\im(\bar V R_V)
		-\rho_U\re(\overline{R_V}\nabla V-\bar V\nabla R_V)\bigr]
		\dd x\dd y .
	\end{equation}
	
	We first show that $J_1+J_2+J_3$ equals the first term in
	\eqref{principal}. In the integrations by parts below, we use
	\begin{equation*}
		\partial_{x_k}(\partial_j A(x-y))=\partial_{jk}A(x-y),
		\qquad
		\partial_{y_k}(\partial_j A(x-y))=-\partial_{jk}A(x-y).
	\end{equation*}
	For \eqref{eq:radial-expansion-density}, write
	$\partial_j\Delta\rho=\sum_{k=1}^2\partial_k\partial_{jk}\rho$
	and integrate by parts in $x_k$ and $y_k$, respectively. This gives
	\begin{equation*}
		J_1=-\frac12\sum_{j,k=1}^2\iint\partial_{jk}A
		\bigl[(\partial_{jk}\rho_U)\rho_V
		+\rho_U\partial_{jk}\rho_V\bigr]\dd x\dd y .
	\end{equation*}
	Now set $z=x-y$. At fixed $z$, the coefficient $\partial_{jk}A(z)$
	is constant in $x$. Integrating the first term by parts in $x_k$
	and the second in $x_j$ gives
	\begin{equation*}
		-\frac12\int_{\R^2}
		\bigl[(\partial_{jk}\rho_U)(t,x)\rho_V(t,x-z)
		+\rho_U(t,x)(\partial_{jk}\rho_V)(t,x-z)\bigr]\dd x
		=\int_{\R^2}(\partial_j\rho_U)(t,x)
		(\partial_k\rho_V)(t,x-z)\dd x .
	\end{equation*}
	Returning to $(x,y)$, we obtain
	\begin{equation*}
		J_1=\sum_{j,k=1}^2\iint\partial_{jk}A\,
		(\partial_j\rho_U)(\partial_k\rho_V)\dd x\dd y .
	\end{equation*}
	
	For \eqref{eq:radial-expansion-stress}, integration by parts in
	$x_k$ and $y_k$ gives
	\begin{equation*}
		J_2=2\sum_{j,k=1}^2\iint\partial_{jk}A
		\bigl[\rho_V\re(\partial_j\bar U\,\partial_kU)
		+\rho_U\re(\partial_j\bar V\,\partial_kV)\bigr]\dd x\dd y .
	\end{equation*}
	For \eqref{eq:radial-expansion-current}, integrating the two
	divergences by parts gives
	\begin{equation*}
		J_3=-2\sum_{j,k=1}^2\iint\partial_{jk}A
		\bigl[p_{U,k}p_{V,j}+p_{U,j}p_{V,k}\bigr]\dd x\dd y
		=-4\sum_{j,k=1}^2\iint\partial_{jk}A\,
		p_{U,j}p_{V,k}\dd x\dd y.
	\end{equation*}

	Since $\partial_j\rho_U=2\re(\bar U\partial_jU)$ and
	$p_{U,j}=\im(\bar U\partial_jU)$, we have
	\begin{equation*}
		(\partial_j\rho_U)(\partial_k\rho_V)-4p_{U,j}p_{V,k}
		=4\re\bigl[(\bar U\partial_jU)(\bar V\partial_kV)\bigr].
	\end{equation*}
	These are precisely the mixed terms obtained by expanding
	$\re(\overline{B_j}B_k)$ and summing against the symmetric matrix
	$D^2A$. Indeed, at each $(t,x,y)$,
	\begin{equation*}
		\sum_{j,k=1}^2\partial_{jk}A\,\re(\overline{B_j}B_k)
		=\sum_{j,k=1}^2\partial_{jk}A
		\bigl[\rho_V\re(\partial_j\bar U\,\partial_kU)
		+\rho_U\re(\partial_j\bar V\,\partial_kV)
		+2\re\bigl((\bar U\partial_jU)(\bar V\partial_kV)\bigr)\bigr].
	\end{equation*}
	The first two terms on the right give $J_2$ after multiplication
	by two and integration, and the last gives $J_1+J_3$. Hence
	\begin{equation*}
		J_1+J_2+J_3
		=2\sum_{j,k=1}^2\iint\partial_{jk}A\,
		\re(\overline{B_j}B_k)\dd x\dd y .
	\end{equation*}
	
	For \eqref{eq:radial-expansion-pressure}, integration by parts
	in $x$ and $y$, respectively, gives
	\begin{equation*}
		J_4=\frac12\iint\Delta A\,
		\bigl(\sigma_U\rho_U^2\rho_V+\sigma_V\rho_U\rho_V^2\bigr)
		\dd x\dd y .
	\end{equation*}
	Thus $J_1+J_2+J_3+J_4=I_1$ by \eqref{principal}. Comparing
	\eqref{eq:radial-expansion-source-u} with \eqref{actualsource}
	and \eqref{eq:radial-expansion-source-v} with \eqref{testsource}
	gives $J_5=I_2$ and $J_6=I_3$, proving the identity for
	$\dd\mathcal M/\dd t$.
	
	Finally, applying \eqref{hessian} to $\re B$ and $\im B$ gives
	\begin{equation*}
		\sum_{j,k=1}^2\partial_{jk}A(x-y)\re(\overline{B_j}B_k)
		\gtrsim_\varphi L(1+L|x-y|)^{-3}|B|^2 .
	\end{equation*}
	Since $\Delta A\ge0$ and $\sigma_U,\sigma_V\ge0$, the second term
	in \eqref{principal} is nonnegative. This proves
	\eqref{eq:radial-positive}.

\end{proof}

\begin{corollary}\label{cor:radial-weighted-estimates}
	Assume \eqref{bootstrap} and let $U=P_Nu$. Let $B$ be formed
	from $U,V$ as in Lemma~\ref{lem:forced-radial-identity}.
	If $V=S(t-t_0)f$, where $t_0\in\R$,
	$\supp\widehat f\subset\{|\xi|\le 2H\}$, and $H\ge1$, then
	\begin{equation}\label{commonB}
		\int_J\iint L(1+L|x-y|)^{-3}|B(t,x,y)|^2
		\dd x\dd y\dd t
		\lesssim_M H\norm f{L^2(\R^2)}^2 .
	\end{equation}
	If $V=P_Nu$, then
	\begin{equation}\label{eq:radial-two-projected}
		\int_J\iint L(1+L|x-y|)^{-3}|B(t,x,y)|^2
		\dd x\dd y\dd t
		\lesssim_M 1 .
	\end{equation}
	Both estimates hold uniformly under spatial
	translations of $U$ and $V$.
\end{corollary}

\begin{proof}
	By Lemma~\ref{lem:forced-radial-identity} and
	\eqref{eq:radial-positive}, integration over $J=[0,t_1]$ gives
	\begin{equation}\label{eq:radial-integrated}
		\int_J\iint L(1+L|x-y|)^{-3}|B(t,x,y)|^2
		\dd x\dd y\dd t
		\lesssim
		|\mathcal M(t_1)|+|\mathcal M(0)|
		+\int_J\bigl(|I_2(t)|+|I_3(t)|\bigr)\dd t .
	\end{equation}
	
	Since $\norm{\nabla A}{L^\infty(\R^2)}\le1$,
	\eqref{eq:radial-action} gives
	\begin{equation*}
		|\mathcal M(t)|
		\le \norm{p_U(t)}{L^1(\R^2)}\norm{\rho_V(t)}{L^1(\R^2)}
		+\norm{\rho_U(t)}{L^1(\R^2)}\norm{p_V(t)}{L^1(\R^2)} .
	\end{equation*}
	Mass conservation bounds $\norm{P_Nu(t)}{L^2(\R^2)}$ in
	terms of $M$, and the first bound in \eqref{bootstrap} bounds
	$\norm{\nabla P_Nu(t)}{L^2(\R^2)}$.
	Since $\rho_U=|U|^2$ and $|p_U|\le |U||\nabla U|$, it follows
	that
	\begin{equation}\label{eq:radial-action-mass-momentum}
		|\mathcal M(t)|
		\lesssim_M \norm{\rho_V(t)}{L^1(\R^2)}
		+\norm{p_V(t)}{L^1(\R^2)} .
	\end{equation}
	
For \eqref{actualsource}, 
$R_U=R=P_NF(u)-F(P_Nu)$. Integrating first in $y$
and using $\re(\bar R\,\partial_jU)=\re(R\,\partial_j\bar U)$,
we obtain $I_2(t)=J_1(t)+J_2(t)$, where
\begin{equation}\label{eq:radial-I2-first}
	J_1(t):=\sum_{j=1}^2\re\int_{\R^2}
	\bigl[(\partial_jA)*\rho_V\bigr](t,x)
	\bigl[R\,\partial_j\bar U-\bar U\,\partial_jR\bigr](t,x)\dd x
\end{equation}
and
\begin{equation}\label{eq:radial-I2-second}
	J_2(t):=-2\sum_{j=1}^2\im\int_{\R^2}
	\bigl[(\partial_jA)*p_{V,j}\bigr](t,x)
	(\bar U R)(t,x)\dd x .
\end{equation}

For \eqref{eq:radial-I2-first}, take
$q=(\partial_jA)*\rho_V$ in \eqref{qs2}. By \eqref{weight},
$\supp\widehat{\partial_jA}\subset\{|\xi|\le L\}$, so $q$ has
the required spatial Fourier support. Since
$\norm{\partial_jA}{L^\infty(\R^2)}\le1$, we also have
\begin{equation*}
	\norm{q(t)}{L^\infty(\R^2)}
	\le \norm{\rho_V(t)}{L^1(\R^2)} .
\end{equation*}
Thus \eqref{qs2} bounds the two integrals containing
$R\,\partial_j\bar U$ and $\bar U\,\partial_jR$, respectively,
and summing over $j$ gives
\begin{equation*}
	\int_J|J_1(t)|\dd t
	\lesssim ND_4\sup_{t\in J}\norm{\rho_V(t)}{L^1(\R^2)} .
\end{equation*}

For \eqref{eq:radial-I2-second}, take
$q=(\partial_jA)*p_{V,j}$ in \eqref{qs1}. Again $q$ has spatial
Fourier support in $\{|\xi|\le L\}$, and
$\norm{q(t)}{L^\infty(\R^2)}\le\norm{p_{V,j}(t)}{L^1(\R^2)}$.
Applying \eqref{qs1} to the integral containing $\bar U R$ and
summing over $j$ gives
\begin{equation*}
	\int_J|J_2(t)|\dd t
	\lesssim D_4\sup_{t\in J}\norm{p_V(t)}{L^1(\R^2)} .
\end{equation*}
Adding the two bounds yields
\begin{equation}\label{eq:radial-source-bound}
	\int_J|I_2(t)|\dd t
	\lesssim D_4\left(
	N\sup_{t\in J}\norm{\rho_V(t)}{L^1(\R^2)}
	+\sup_{t\in J}\norm{p_V(t)}{L^1(\R^2)}\right).
\end{equation}

Moreover, estimate \eqref{eq:radial-source-bound} is uniform under translations of $U$ and $V$, which preserve spatial Fourier support and $L^\infty(\R^2)$ norm.
	
For \eqref{commonB}, $V=S(t-t_0)f$ and $\supp\widehat f\subset\{|\xi|\le 2H\}$. 
Cauchy--Schwarz and Plancherel give
\begin{equation*}
	\norm{\rho_V(t)}{L^1(\R^2)}=\norm{V(t)}{L^2(\R^2)}^2,
	\qquad
	\norm{p_V(t)}{L^1(\R^2)}
	\le \norm{V(t)}{L^2(\R^2)}\norm{\nabla V(t)}{L^2(\R^2)}
	\le 2H\norm{V(t)}{L^2(\R^2)}^2 .
\end{equation*}
Substituting these bounds into \eqref{eq:radial-action-mass-momentum}
and \eqref{eq:radial-source-bound}, respectively, gives
\begin{equation}\label{actionbound}
	\sup_{t\in J}|\mathcal M(t)|
	\lesssim_M H\norm{V}{L^\infty(J;L^2(\R^2))}^2,
	\qquad
	\int_J|I_2(t)|\dd t
	\lesssim_M (N+H)D_4\norm{V}{L^\infty(J;L^2(\R^2))}^2 .
\end{equation}
Unitarity gives
$\norm{V}{L^\infty(J;L^2(\R^2))}=\norm f{L^2(\R^2)}$,
and $R_V=0$ gives $I_3=0$ by \eqref{testsource}.
Since $N^2D_4\lesssim_M N^{-\varepsilon/2}$ by \eqref{variation}
and $N,H\ge1$, we have $(N+H)D_4\lesssim_M H$.
Substituting \eqref{actionbound} into \eqref{eq:radial-integrated}
therefore proves \eqref{commonB}.

For \eqref{eq:radial-two-projected}, $U=V=P_Nu$. Since
$\sup_{t\in J}\norm{P_Nu(t)}{H^1(\R^2)}\lesssim_M1$, we have
\begin{equation*}
	\norm{\rho_V(t)}{L^1(\R^2)}+\norm{p_V(t)}{L^1(\R^2)}
	\le \norm{V(t)}{L^2(\R^2)}^2
	+\norm{V(t)}{L^2(\R^2)}\norm{\nabla V(t)}{L^2(\R^2)}
	\lesssim_M1 .
\end{equation*}
Thus \eqref{eq:radial-action-mass-momentum} gives
$\sup_{t\in J}|\mathcal M(t)|\lesssim_M1$, and
\eqref{eq:radial-source-bound} gives
$\int_J|I_2(t)|\dd t\lesssim_M(1+N)D_4$.
Since $U=V$ and $R_U=R_V$, interchanging $x$ and $y$ in
\eqref{testsource} and using $\nabla A(y-x)=-\nabla A(x-y)$
gives $I_3(t)=I_2(t)$. Estimate \eqref{eq:radial-two-projected}
now follows from \eqref{eq:radial-integrated} and
$(1+N)D_4\lesssim_M1$.

Both estimates hold uniformly under independent spatial translations
of $U$ and $V$, since \eqref{eq:radial-action-mass-momentum} and
\eqref{eq:radial-source-bound} are uniform under these translations,
and interchanging the two fields gives the same bound for $I_3$.

\end{proof}

To pass from the weighted integrals in
Corollary~\ref{cor:radial-weighted-estimates} to a product at
$x=y$, we use the Fourier support in the first variable.
Suppose that a function $D(x,y)$ has Fourier support in
$\{|\xi|\le K\}$ in the $x$ variable for each fixed $y$.
Then $P_KD(\cdot,y)=D(\cdot,y)$. Writing $k_K$ for the
convolution kernel of $P_K$, we obtain by Cauchy--Schwarz
\begin{equation*}
	|D(y,y)|^2
	=\left|\int_{\R^2}k_K(y-x)D(x,y)\dd x\right|^2
	\le \norm{k_K}{L^1(\R^2)}
	\int_{\R^2}|k_K(y-x)|\,|D(x,y)|^2\dd x .
\end{equation*}
Since $k_K(w)=K^2k_1(Kw)$ and $k_1$ is Schwartz, it follows that
\begin{equation}\label{diag}
	|D(y,y)|^2
	\lesssim K^2\int_{\R^2}(1+K|x-y|)^{-4}|D(x,y)|^2\dd x .
\end{equation}
For $K\ge L$, the weight in \eqref{diag} satisfies
\begin{equation}\label{eq:radial-diagonal-cost}
	K^2(1+K|w|)^{-4}
	\le \frac{K^2}{L}\,L(1+L|w|)^{-3}.
\end{equation}

We now apply \eqref{diag} to the products in
Corollary~\ref{cor:radial-weighted-estimates}.

\begin{lemma}\label{lem:radial-product-estimates}
	Assume \eqref{bootstrap}, let $L=N^{1-\varepsilon}$, and fix
	$c>0$. For $N$ sufficiently large depending on $a,c,\varepsilon$,
	the following estimates hold uniformly in $J$.
	
	If $H\ge cN$ and $f\in L^2(\R^2)$ has Fourier support in
	$\{H/2\le|\xi|\le2H\}$, then, for every $t_0\in\R$,
	\begin{equation}\label{freehigh}
		Y\bigl(S(t-t_0)f\bigr)
		\lesssim_M (L/H)^{1/2}\norm f{L^2(\R^2)} .
	\end{equation}
	For all $y,z\in\R^2$,
	\begin{equation}\label{vv}
		\norm{\nabla_x\bigl[(P_Nu)(t,x+y)
			\overline{(P_Nu)(t,x+z)}\bigr]}{L^2(J\times\R^2)}
		\lesssim_M NL^{-1/2}.
	\end{equation}
	If $f\in L^2(\R^2)$ has Fourier support in $\{|\xi|\le4N\}$,
	then, for all $y,z\in\R^2$,
	\begin{equation}\label{vf}
		\norm{\nabla_x\bigl[(P_Nu)(t,x+y)
			\overline{(S(t)f)(x+z)}\bigr]}{L^2(J\times\R^2)}
		\lesssim_M N^{3/2}L^{-1/2}\norm f{L^2(\R^2)} .
	\end{equation}
	Finally, set $b=aN/4$. Then
	\begin{equation}\label{shell}
		Y\bigl((P_{aN}-P_b)P_Nu\bigr)
		\lesssim_M \frac{\sqrt L}{b-2L}
		\sim \frac{\sqrt L}{N}.
	\end{equation}
\end{lemma}

\begin{proof}
	For \eqref{freehigh}, set $V(t)=S(t-t_0)f$. We apply
	\eqref{commonB} with the pair $(P_Nu,V)$ and then project
	the first factor to $P_Lu$. To make this step explicit, set
	$z=x-y$. For the pair $(P_Nu,V)$ in
	Corollary~\ref{cor:radial-weighted-estimates}, the chain rule gives
	\begin{equation*}
		B(t,x,x-z)
		=\nabla_x\bigl[(P_Nu)(t,x)\overline{V(t,x-z)}\bigr].
	\end{equation*}
	Since $N$ is sufficiently large,
	$P_LP_Nu=P_Lu$. Writing $k_L$ for the kernel of $P_L$, we have
	\begin{equation*}
		\nabla_x\bigl[(P_Lu)(t,x)\overline{V(t,x-z)}\bigr]
		=\int_{\R^2}k_L(s)\,
		\nabla_x\bigl[(P_Nu)(t,x-s)\overline{V(t,x-z)}\bigr]\dd s .
	\end{equation*}
	For each $s$, \eqref{commonB} applies to the translate
	$(P_Nu)(t,x-s)$ with the same bound. Minkowski's inequality
	in $L^2$ with measure $L(1+L|z|)^{-3}\dd x\dd z\dd t$
	therefore gives
	\begin{equation*}
		\int_J\iint L(1+L|z|)^{-3}
		\bigl|\nabla_x[(P_Lu)(t,x)\overline{V(t,x-z)}]\bigr|^2
		\dd x\dd z\dd t
		\lesssim_M \norm{k_L}{L^1(\R^2)}^2 H\norm f{L^2(\R^2)}^2 .
	\end{equation*}
	Here $\norm{k_L}{L^1(\R^2)}=\norm{k_1}{L^1(\R^2)}$.
	
	For fixed $t,z$, the product inside the gradient has Fourier
	support in $\{|\xi|\ge H/2-2L\}$. Since $2L\le cN/4\le H/4$,
	Plancherel gives
	\begin{equation*}
		\int_{\R^2}|(P_Lu)(t,x)V(t,x-z)|^2\dd x
		\le \frac{16}{H^2}\int_{\R^2}
		\bigl|\nabla_x[(P_Lu)(t,x)\overline{V(t,x-z)}]\bigr|^2\dd x .
	\end{equation*}
	Multiplying by $L(1+L|z|)^{-3}$, integrating in $z,t$, and
	changing variables back to $y=x-z$, we obtain
	\begin{equation*}
		\int_J\iint L(1+L|x-y|)^{-3}|(P_Lu)(t,x)V(t,y)|^2
		\dd x\dd y\dd t
		\lesssim_M H^{-1}\norm f{L^2(\R^2)}^2 .
	\end{equation*}
	For fixed $t,y$, the product $(P_Lu)(t,x)V(t,y)$ has
	$x$-Fourier support in $\{|\xi|\le2L\}$. Applying \eqref{diag}
	with $K=2L$ and $D(x,y)=(P_Lu)(t,x)V(t,y)$, \eqref{eq:radial-diagonal-cost} gives
	\begin{equation*}
		\norm{(P_Lu)V}{L^2(J\times\R^2)}^2
		\lesssim L\int_J\iint L(1+L|x-y|)^{-3}
		|(P_Lu)(t,x)V(t,y)|^2\dd x\dd y\dd t
		\lesssim_M \frac{L}{H}\norm f{L^2(\R^2)}^2 .
	\end{equation*}
	Replacing $(P_Lu)(t,x)$ by $(P_Lu)(t,x+d)$ replaces $x-s$
	by $x+d-s$ in the convolution above, so the bound is uniform
	in $d$. Taking square roots and the supremum in translations
	proves \eqref{freehigh}.
	
	For \eqref{vv} and \eqref{vf}, we use $B$ directly with the
	pairs in Corollary~\ref{cor:radial-weighted-estimates}. In both
	cases its first factor is $P_Nu$, so each component of $B(t,x,y)$
	has $x$-Fourier support in $\{|\xi|\le2N\}$ for fixed $t,y$.
	On the diagonal,
	\begin{equation*}
		B(t,x,x)=\nabla_x\bigl[(P_Nu)(t,x)\overline{V(t,x)}\bigr].
	\end{equation*}
	Applying \eqref{diag} componentwise with $K=2N$ and using
	\eqref{eq:radial-diagonal-cost}, we obtain
	\begin{equation*}
		\norm{\nabla_x[(P_Nu)\bar V]}{L^2(J\times\R^2)}^2
		\lesssim \frac{N^2}{L}
		\int_J\iint L(1+L|x-y|)^{-3}|B(t,x,y)|^2
		\dd x\dd y\dd t .
	\end{equation*}
	For \eqref{vv}, take $V=P_Nu$. By
	\eqref{eq:radial-two-projected}, the integral is bounded by
	$\lesssim_M 1$.
	For \eqref{vf}, take $V=S(t)f$. Since
	$\supp\widehat f\subset\{|\xi|\le4N\}$, \eqref{commonB}
	with $H=2N$ bounds the integral by
	$\lesssim_M N\norm f{L^2(\R^2)}^2$. Taking square roots gives
	\eqref{vv} and \eqref{vf}, including the stated translations
	by the uniformity in Corollary~\ref{cor:radial-weighted-estimates}.
	
	For \eqref{shell}, start with the pair $(P_Nu,P_Nu)$ in
	\eqref{eq:radial-two-projected}. We apply $P_L$ to the first
	factor and $P_{aN}-P_b$ to the second. The corresponding
	kernels are $k_L$ and $k_{aN}-k_b$, with
	\begin{equation*}
		\norm{k_L}{L^1(\R^2)}=\norm{k_1}{L^1(\R^2)},\qquad
		\norm{k_{aN}-k_b}{L^1(\R^2)}\le2\norm{k_1}{L^1(\R^2)} .
	\end{equation*}
	Expanding both convolutions gives, for each pair of translation
	parameters, the gradient of a product of two 
	translated copies of $P_Nu$. Thus Minkowski's inequality in
	the two convolution variables and \eqref{eq:radial-two-projected}
	give
	\begin{equation*}
		\int_J\iint L(1+L|z|)^{-3}
		\bigl|\nabla_x\bigl[(P_Lu)(t,x)
		\overline{((P_{aN}-P_b)P_Nu)(t,x-z)}\bigr]\bigr|^2
		\dd x\dd z\dd t
		\lesssim_M1 .
	\end{equation*}
	The second factor has Fourier support in
	$\{b\le|\xi|\le2aN\}$. Hence, for fixed $t,z$, the product
	inside the gradient has Fourier support in
	$\{|\xi|\ge b-2L\}$. Since $b=aN/4$, we may take $N$
	sufficiently large that $2L\le b/2$. Applying Plancherel
	at fixed $z$ and then changing variables to $y=x-z$ gives
	\begin{equation*}
		\int_J\iint L(1+L|x-y|)^{-3}
		\bigl|(P_Lu)(t,x)((P_{aN}-P_b)P_Nu)(t,y)\bigr|^2
		\dd x\dd y\dd t
		\lesssim_M (b-2L)^{-2}.
	\end{equation*}
	For fixed $t,y$, only $P_Lu$ depends on $x$. Applying
	\eqref{diag} with $K=2L$ and using \eqref{eq:radial-diagonal-cost}
	therefore yields
	\begin{equation*}
		\norm{(P_Lu)((P_{aN}-P_b)P_Nu)}{L^2(J\times\R^2)}^2
		\lesssim_M \frac{L}{(b-2L)^2}.
	\end{equation*}
	The bound is uniform under translation of $P_Lu$ by the
	translation uniformity of \eqref{eq:radial-two-projected}.
	Taking square roots and the supremum, and using
	$b-2L\ge b/2=aN/8$, proves \eqref{shell}.
	
\end{proof}

\begin{remark}
	The choice of $A$ in \eqref{weight} serves two purposes.
	First, the weights paired with $R$ in \eqref{actualsource}
	must satisfy the Fourier support and $L^\infty$ hypotheses
	of \eqref{qs1} and \eqref{qs2}. Convolution with
	$L^2\varphi(L\,\cdot)$ gives
	$\supp\widehat{\nabla A}\subset\{|\xi|\le L\}$ and
	preserves the gradient bound for $\sqrt{L^{-2}+|z|^2}$.
	This gives \eqref{eq:radial-source-bound}.
	
	Second, for each translation required from $Y(h)$, the desired
	product is the value of $(P_Lu)(t,x+d)h(t,y)$ at $x=y$.
	To see what this requires of $A$, fix $t,y,d$. Since
	$(P_Lu)(t,x+d)$ has $x$-Fourier support in $\{|\xi|\le 2L\}$,
	\eqref{diag} with $K=2L$ gives
	\begin{equation*}
		|(P_Lu)(t,y+d)h(t,y)|^2
		\lesssim \int_{\R^2}\frac{L^2}{(1+L|x-y|)^4}
		|(P_Lu)(t,x+d)h(t,y)|^2\dd x .
	\end{equation*}
	The choice of $A$ in \eqref{weight} gives the Hessian lower bound
	\eqref{hessian} in every direction. Thus, at fixed $z=x-y$,
	the virial estimate controls the full gradient
	$\nabla_x[(P_Lu)(t,x+d)\overline{h(t,x-z)}]$ with weight
	$L(1+L|z|)^{-3}$ in the squared norm. Since this weight is
	constant in $x$ for fixed $z$, the frequency separation in the
	proof of \eqref{freehigh} allows us to remove the gradient by
	Plancherel while retaining the same weight. It suffices for the
	diagonal estimate above because
	\begin{equation*}
		\frac{L^2}{(1+L|z|)^4}
		\le L\,\frac{L}{(1+L|z|)^3}.
	\end{equation*}
	This comparison is the reason for choosing $A$ with the lower
	bound \eqref{hessian}, and it gives the factor $L$ in the squared
	estimate \eqref{freehigh}.
\end{remark}

\vspace{7pt}
\subsection{Duhamel estimates}

Lemma~\ref{lem:radial-product-estimates} bounds the products in
\eqref{freehigh} and \eqref{vf} when one factor is a free
Schr\"odinger solution. We now use Lemma~\ref{lem:transfer} to
replace that factor by a Duhamel integral. We then
use \eqref{freehigh} by duality to treat the case in which the
low-frequency factor occurs inside the time integral.

Assume \eqref{bootstrap}, and let $I\subset J$ be an interval.
For a spatial Fourier multiplier $m(D)$ of $L^2(\R^2)$ operator
norm at most one, write
\begin{equation*}
	\mathcal T_m F(t)
	:=\int_{I\cap\{\tau<t\}}S(t-\tau)m(D)F(\tau)\dd\tau .
\end{equation*}
Each estimate below also holds with $\tau>t$ in place of $\tau<t$.

Suppose $H\ge cN$ and $\supp m\subset\{H/2\le|\xi|\le2H\}$.
Then
\begin{equation}\label{asd}
	\norm{(P_Lu)(t,x+d)\,\mathcal T_mF(t,x)}{L^2(I\times\R^2)}
	\lesssim_M (L/H)^{1/2}\norm F{L^{4/3}(I\times\R^2)},
	\quad \forall d\in\R^2 .
\end{equation}
Indeed, $m(D)f$ has the Fourier support required in
\eqref{freehigh}, and $\norm{m(D)f}{L^2(\R^2)}\le\norm f{L^2(\R^2)}$.
Thus \eqref{freehigh} gives
\begin{equation*}
	\norm{(P_Lu)(t,x+d)(S(t)m(D)f)(x)}{L^2(I\times\R^2)}
	\lesssim_M (L/H)^{1/2}\norm f{L^2(\R^2)} .
\end{equation*}
Applying Lemma~\ref{lem:transfer} with admissible pair $(4,4)$
 proves \eqref{asd}.

If $\supp m\subset\{|\xi|\le4N\}$, the same choice of exponents
in Lemma~\ref{lem:transfer}, now applied to \eqref{vf}, gives
for every $y,z\in\R^2$
\begin{equation}\label{eq:duhamel-gradient-product}
	\norm{\nabla_x\bigl[(P_Nu)(t,x+y)
		\overline{(\mathcal T_mF)(t,x+z)}\bigr]}{L^2(I\times\R^2)}
	\lesssim_M N^{3/2}L^{-1/2}\norm F{L^{4/3}(I\times\R^2)} .
\end{equation}

For the $L^4$ estimate, \eqref{eq:strichartz} gives
\begin{equation*}
	\norm{S(t)m(D)f}{L^4(I\times\R^2)}
	\lesssim \norm{m(D)f}{L^2(\R^2)}
	\le \norm f{L^2(\R^2)} .
\end{equation*}
Applying Lemma~\ref{lem:transfer} again yields
\begin{equation}\label{eq:duhamel-multiplier-l4}
	\norm{\mathcal T_mF}{L^4(I\times\R^2)}
	\lesssim \norm F{L^{4/3}(I\times\R^2)} .
\end{equation}

We next place the low-frequency factor inside the time integral.
Let $\ell=P_Lu$, let $\ell'$ be any spatial translate of $\ell$
or $\bar\ell$, and suppose $H\ge cN$ and
$\supp m\subset\{H/2\le|\xi|\le2H\}$. Then, for every $t_0\in I$,
\begin{equation}\label{low-input}
	\left\|\int_{t_0}^t S(t-\tau)m(D)(\ell'G)(\tau)\dd\tau
	\right\|_{L^4(I\times\R^2)}
	\lesssim_{M,c}(L/H)^{1/2}\norm G{L^2(I\times\R^2)} .
\end{equation}

In fact, we first integrate over all of $I$. For $f\in L^2(\R^2)$,
unitarity and Cauchy--Schwarz give
\begin{equation*}
	\left|\left\langle f,
	\int_I S(-\tau)m(D)(\ell'G)(\tau)\dd\tau
	\right\rangle_{L^2(\R^2)}\right|
	\le
	\norm{\overline{\ell'}\,S(\tau)m(D)^*f}{L^2(I\times\R^2)}
	\norm G{L^2(I\times\R^2)} .
\end{equation*}
The symbol of $m(D)^*$ is $\bar m$, so $m(D)^*f$ has Fourier
support in $\{H/2\le|\xi|\le2H\}$ and
$\norm{m(D)^*f}{L^2(\R^2)}\le\norm f{L^2(\R^2)}$.
Since $|\ell'|$ is a translate of $|P_Lu|$, \eqref{freehigh}
bounds the first factor on the right by
$\lesssim_M (L/H)^{1/2}\norm f{L^2(\R^2)}$. Taking the supremum
over $\norm f{L^2(\R^2)}=1$ gives
\begin{equation*}
	\left\|\int_I S(-\tau)m(D)(\ell'G)(\tau)\dd\tau
	\right\|_{L^2(\R^2)}
	\lesssim_M (L/H)^{1/2}\norm G{L^2(I\times\R^2)} .
\end{equation*}
Applying \eqref{eq:strichartz} with the admissible pair $(4,4)$
 and using $S(t)S(-\tau)=S(t-\tau)$ yield
\begin{equation*}
	\left\|\int_I S(t-\tau)m(D)(\ell'G)(\tau)\dd\tau
	\right\|_{L^4(I\times\R^2)}
	\lesssim_M (L/H)^{1/2}\norm G{L^2(I\times\R^2)} .
\end{equation*}
Moreover,
Proposition~\ref{prop:ck} gives the same bound with the integral
restricted to $\tau<t$ or to $\tau>t$. For $t\ge t_0$, apply
the first bound with $G$ restricted to $I\cap[t_0,\infty)$.
For $t<t_0$, apply the second with $G$ restricted to
$I\cap(-\infty,t_0]$.
Combining the two bounds proves \eqref{low-input}.

\vspace{10pt}
\section{Global well-posedness at positive regularity}\label{sec:positive}

\vspace{7pt}
\subsection{Time truncation of the Duhamel integral and the logarithmic loss}
Estimate \eqref{asd} places a translate of $P_Lu$ outside the
Duhamel integral, and \eqref{low-input} places one inside. We now
need both factors in the same estimate. 

\begin{lemma}\label{lem:polynomial-time-restriction}
	Assume \eqref{bootstrap}. Let $H\ge cN$, and let $m(D)$ be a
	spatial Fourier multiplier with $|m|\le1$ and
	$\supp m\subset\{H/2\le|\xi|\le2H\}$. Let $b_1,b_2$ be
	 spatial translates of $\ell=P_Lu$ or $\bar\ell$.
	Then
	\begin{equation}\label{poslog}
		\left\|b_1(t)\int_{J\cap\{\tau<t\}}
		S(t-\tau)m(D)(b_2G)(\tau)\dd\tau\right\|_{L^2(J\times\R^2)}
		\lesssim_M \frac LH
		\bigl[1+\log(2+|J|HL)\bigr]\norm G{L^2(J\times\R^2)} .
	\end{equation}
	The same estimate holds with $\tau>t$.
\end{lemma}

\begin{proof}
	We first integrate over all of $J$. Since $m(D)^*f$ has Fourier
	support in $\{H/2\le|\xi|\le2H\}$, \eqref{freehigh} gives
	\begin{equation*}
		\left|\left\langle f,
		\int_J S(-\tau)m(D)(b_2G)(\tau)\dd\tau
		\right\rangle_{L^2(\R^2)}\right|
		\le \norm{\overline{b_2}S(\tau)m(D)^*f}{L^2(J\times\R^2)}
		\norm G{L^2(J\times\R^2)}
	\end{equation*}
	\begin{equation*}
			\lesssim_M (L/H)^{1/2}\norm f{L^2(\R^2)}
		\norm G{L^2(J\times\R^2)} .
	\end{equation*}
	The integral paired with $f$ has the same Fourier support.
	Taking the supremum over $f$ and applying \eqref{freehigh}
	to this integral, now with the factor $b_1$, gives
	\begin{equation}\label{eq:unrestricted-two-low}
		\left\|b_1(t)\int_J S(t-\tau)m(D)(b_2G)(\tau)\dd\tau
		\right\|_{L^2(J\times\R^2)}
		\lesssim_M (L/H)\norm G{L^2(J\times\R^2)} .
	\end{equation}
	
Let $k$ be the least nonnegative integer such that
$2^{-k}|J|\le(HL)^{-1}$. For $0\le j\le k$, let $\mathcal D_j$
be the partition of $J$ obtained by $j$ successive bisections,
with $\mathcal D_0=\{J\}$. For each $I\in\mathcal D_j$, $j<k$,
write $I_-$ and $I_+$ for the left and right halves of $I$.

In the coordinates $(\tau,t)$, bisecting $I$ gives the disjoint
decomposition, up to endpoints,
\begin{equation*}
	\{(\tau,t)\in I^2:\tau<t\}
	=(I_-\times I_+)
	\sqcup\{(\tau,t)\in I_-^2:\tau<t\}
	\sqcup\{(\tau,t)\in I_+^2:\tau<t\}.
\end{equation*}
Starting with $I=J$ and iterating this identity $k$ times gives
\begin{equation}\label{eq:positive-time-triangle}
	\{(\tau,t)\in J^2:\tau<t\}
	=\bigsqcup_{j=0}^{k-1}\ \bigsqcup_{I\in\mathcal D_j}
	(I_-\times I_+)
	\;\sqcup\!
	\bigsqcup_{I\in\mathcal D_k}\{(\tau,t)\in I^2:\tau<t\}.
\end{equation}

Fix $j<k$. On $I_-\times I_+$, the condition $\tau<t$ holds
for every $\tau\in I_-$ and $t\in I_+$. Applying
\eqref{eq:unrestricted-two-low} to $\mathbf1_{I_-}G$ and
restricting $t$ to $I_+$ gives
\begin{equation}\label{eq:positive-time-generation}
	\left\|\sum_{I\in\mathcal D_j}\mathbf1_{I_+}(t)b_1(t)
	\int_{I_-}S(t-\tau)m(D)(b_2G)(\tau)\dd\tau
	\right\|_{L^2(J\times\R^2)}^2
	\lesssim_M \left(\frac LH\right)^2
	\sum_{I\in\mathcal D_j}\norm G{L^2(I_-\times\R^2)}^2
	\le \left(\frac LH\right)^2\norm G{L^2(J\times\R^2)}^2 .
\end{equation}

For the last union in \eqref{eq:positive-time-triangle}, fix
$I\in\mathcal D_k$. Since $b_1,b_2$ are translates of $P_Lu$
or its conjugate, Bernstein and mass conservation give
\begin{equation*}
	\norm{b_i}{L^\infty(J\times\R^2)}
	=\norm{P_Lu}{L^\infty(J\times\R^2)}
	\lesssim L\norm u{L^\infty(J;L^2(\R^2))}
	\lesssim_M L,\qquad i=1,2 .
\end{equation*}
Using these two bounds,
$\norm{S(t-\tau)m(D)}{L^2(\R^2)\to L^2(\R^2)}\le1$ and Minkowski's inequality, we obtain
for $t\in I$
\begin{equation*}
	\left\|b_1(t)\int_{I\cap\{\tau<t\}}S(t-\tau)m(D)
	(b_2G)(\tau)\dd\tau\right\|_{L^2(\R^2)}
	\lesssim_M L^2\int_I\norm{G(\tau)}{L^2(\R^2)}\dd\tau
	\le L^2|I|^{1/2}\norm G{L^2(I\times\R^2)} .
\end{equation*}
Taking the $L^2(I)$ norm in $t$ gives another factor $|I|^{1/2}$.
Since $L^2|I|=L^2 2^{-k}|J|\le L/H$, squaring and summing over
$I\in\mathcal D_k$ gives
\begin{equation}\label{eq:positive-time-last}
	\sum_{I\in\mathcal D_k}
	\left\|b_1(t)\int_{I\cap\{\tau<t\}}S(t-\tau)m(D)
	(b_2G)(\tau)\dd\tau\right\|_{L^2(I\times\R^2)}^2
	\lesssim_M \left(\frac LH\right)^2
	\sum_{I\in\mathcal D_k}\norm G{L^2(I\times\R^2)}^2
	=\left(\frac LH\right)^2\norm G{L^2(J\times\R^2)}^2 .
\end{equation}

Taking square roots in \eqref{eq:positive-time-generation}
and \eqref{eq:positive-time-last}, then applying the triangle
inequality to the decomposition \eqref{eq:positive-time-triangle},
yields
\begin{equation*}
	\left\|b_1(t)\int_{J\cap\{\tau<t\}}S(t-\tau)m(D)
	(b_2G)(\tau)\dd\tau\right\|_{L^2(J\times\R^2)}
	\lesssim_M (k+1)\frac LH\norm G{L^2(J\times\R^2)} .
\end{equation*}
If $k\ge1$, its minimality gives $2^{k-1}<|J|HL\le2^k$.
Thus $k+1\lesssim1+\log(2+|J|HL)$, also for $k=0$, which
proves \eqref{poslog}. For $\tau>t$, the rectangles in
\eqref{eq:positive-time-triangle} are $I_+\times I_-$, and the
bound on each $I\in\mathcal D_k$ still follows from
$L^2|I|\le L/H$.

\end{proof}

\vspace{7pt}
\subsection{Improving the bound for $Y(Q_{aN}u)$}

We next improve the bound for $Y(Q_{aN}u)$ in \eqref{bootstrap}.
Writing $u=P_Lu+Q_Lu$ and expanding $F(u)=|u|^2u$, we obtain
\begin{equation}\label{eq:positive-cubic-expansion}
	F(u)=F(P_Lu)
	+2|P_Lu|^2Q_Lu+(P_Lu)^2\overline{Q_Lu}
	+2P_Lu\,|Q_Lu|^2+\overline{P_Lu}(Q_Lu)^2
	+F(Q_Lu).
\end{equation}
Let $b=aN/4$ as in \eqref{shell}. The two terms linear in $Q_Lu$
will be iterated, so we write
\begin{equation*}
	Kz(t):=-\ii\int_0^t S(t-\tau)Q_b
	\bigl(2|P_Lu|^2z+(P_Lu)^2\bar z\bigr)(\tau)\dd\tau .
\end{equation*}
Since $\supp\widehat{F(P_Lu)}\subset\{|\xi|\le6L\}$ and
$L/N=N^{-\varepsilon}$, we have $Q_bF(P_Lu)=0$ for $N$
sufficiently large. Applying $Q_b$ to \eqref{eq:duhamel} and
using \eqref{eq:positive-cubic-expansion} therefore gives
\begin{equation}\label{duhamelh}
	Q_bu(t)-K(Q_Lu)(t)
	=S(t)Q_bu(0)-\ii\int_0^t S(t-\tau)Q_b
	\bigl(2P_Lu\,|Q_Lu|^2+\overline{P_Lu}(Q_Lu)^2
	+F(Q_Lu)\bigr)(\tau)\dd\tau .
\end{equation}

We state the hypothesis on the time integral separately so that
the same proposition can be used in Section~\ref{sec:endpoint}.
For $P=\operatorname{Id}$ or $Q_{aN}$, the dyadic decomposition is
\begin{equation*}
	PQ_b=\sum_{H\ge b/2}(P_H-P_{H/2})PQ_b .
\end{equation*}
Indeed, the sum telescopes, and the summands with $2H<b$ vanish.
The symbol of each summand is bounded in absolute value by one
and supported in $\{H/2\le|\xi|\le2H\}$. Let $d_H$ be the least
constant such that
\begin{equation}\label{eq:positive-annular-bound}
	\left\|b_1(t)\int_0^t S(t-\tau)(P_H-P_{H/2})PQ_b
	(b_2G)(\tau)\dd\tau\right\|_{L^2(J\times\R^2)}
	\le d_H\norm G{L^2(J\times\R^2)}
\end{equation}
for both choices of $P$ and all spatial translates
$b_1,b_2$ of $P_Lu$ or $\overline{P_Lu}$. The hypothesis is
\begin{equation}\label{operator-small}
	q:=\sum_{H\ge b/2}d_H\le N^{-\varepsilon/2}.
\end{equation}
\begin{remark}\label{J}
	If $|J|\le N^p$ for fixed $p>0$, then \eqref{poslog} applies on
	every annulus since $H\ge b/2=aN/8$. At $H\sim 2^jN$,
	\begin{equation*}
		\frac LH\bigl[1+\log(2+|J|HL)\bigr]
		\lesssim_{p}2^{-j}\frac LN(1+\log N+j).
	\end{equation*}
	Consequently,
	\begin{equation*}
		\sum_{H\ge b/2}d_H
		\lesssim_{M,p}\frac LN\sum_{j\ge0}2^{-j}(1+\log N+j)
		\lesssim N^{-\varepsilon}(1+\log N)
		\le N^{-\varepsilon/2}
	\end{equation*}
	for $N$ sufficiently large depending on $M,a,p,\varepsilon$.
\end{remark}

\begin{proposition}\label{high}
	Under \eqref{bootstrap} and \eqref{operator-small}, we have
	\begin{equation}\label{Ximprove}
		Y(Q_{aN}u)\lesssim_M N^{-1-\varepsilon/2}.
	\end{equation}
\end{proposition}

\begin{proof}
	To estimate $Q_{aN}u$, we use $Q_{aN}u=Q_{aN}Q_bu$, which
	follows from $2b<aN$. The right-hand side of \eqref{duhamelh}
	gives an explicit formula for $Q_bu-K(Q_Lu)$. To recover
	$Q_bu$ from this difference, write
	$Q_Lu=Q_bu+(P_b-P_L)u$. By the linearity of $K$,
	\begin{equation*}
		Q_bu=[Q_bu-K(Q_Lu)]+K((P_b-P_L)u)+K(Q_bu).
	\end{equation*}
	Fix an integer $k$ with $k\varepsilon>4$. Substituting this
	identity into the last term $k-1$ times and applying $Q_{aN}$
	gives
	\begin{equation*}
		Q_{aN}u
		=\sum_{j=0}^{k-1}Q_{aN}K^j[Q_bu-K(Q_Lu)]
		+\sum_{j=1}^{k}Q_{aN}K^j((P_b-P_L)u)
		+Q_{aN}K^k(Q_bu).
	\end{equation*}
	
	The second sum vanishes for $N$ sufficiently large. Indeed,
	if $z(t)$ has Fourier support in $\{|\xi|\le R\}$ for every
	$t\in J$, then
	\begin{equation*}
		\supp\widehat{\,2|P_Lu|^2z+(P_Lu)^2\bar z\,}(t)
		\subset\{|\xi|\le R+4L\}.
	\end{equation*}
	The multipliers $Q_b$ and $S(t-\tau)$ in the definition of
	$K$ preserve this inclusion. Since $L<b$ for large $N$,
	$(P_b-P_L)u$ has frequencies at most $2b$. Induction gives
	\begin{equation*}
		\supp\widehat{K^j((P_b-P_L)u)}(t)
		\subset\{|\xi|\le2b+4jL\}.
	\end{equation*}
	Now $2b=aN/2$ and $L/N=N^{-\varepsilon}$, so we may choose
	$N$ large enough that $2b+4kL<aN$. Hence, to prove
	\eqref{Ximprove}, it suffices to bound the two terms on the
	right of
	\begin{equation}\label{eq:positive-high-reduction}
		Y(Q_{aN}u)
		\le \sum_{j=0}^{k-1}Y\bigl(Q_{aN}K^j[Q_bu-K(Q_Lu)]\bigr)
		+Y\bigl(Q_{aN}K^k(Q_bu)\bigr).
	\end{equation}
	
	To control the powers of $K$ in \eqref{eq:positive-high-reduction},
	we claim that, for $P=\operatorname{Id}$ or $Q_{aN}$,
	\begin{equation}\label{K}
		Y(PKz)\le 3qY(z).
	\end{equation}
	Fix $d\in\R^2$. In \eqref{eq:positive-annular-bound}, take
	$b_1(t,x)=(P_Lu)(t,x+d)$ and $b_2=P_Lu$. For the two terms
	in $Kz$, the respective choices of $G$ are determined by
	\begin{equation*}
		|P_Lu|^2z=(P_Lu)(\overline{P_Lu}\,z),\qquad
		(P_Lu)^2\bar z=(P_Lu)((P_Lu)\bar z).
	\end{equation*}
	Both choices satisfy
	\begin{equation*}
		\norm{\overline{P_Lu}\,z}{L^2(J\times\R^2)}
		=\norm{(P_Lu)\bar z}{L^2(J\times\R^2)}
		\le Y(z).
	\end{equation*}
	 Thus the two terms on
	annulus $H$ are bounded by $2d_HY(z)$ and $d_HY(z)$,
	respectively. Summing over $H$ and taking the supremum over
	$d$ proves \eqref{K}. In particular, for $j\ge1$,
	\begin{equation*}
		Y(Q_{aN}K^jz)\le 3qY(K^{j-1}z)\le(3q)^jY(z).
	\end{equation*}
	
     For the last term in \eqref{eq:positive-high-reduction}, we claim that
     $Y(Q_{aN}K^k(Q_bu))\lesssim_{M}N^{-1-\varepsilon/2}$. Since $2aN<N$, the
	multiplier $P_N$ equals one on the support of $P_{aN}-P_b$,
	so
	\begin{equation*}
		Q_bu=Q_{aN}u+(P_{aN}-P_b)P_Nu .
	\end{equation*}
	Applying \eqref{bootstrap} and \eqref{shell} gives
	\begin{equation*}
		Y(Q_bu)\lesssim_{M,a}
		N^{-1-\varepsilon/4}+\frac{\sqrt L}{N}.
	\end{equation*}
	By \eqref{operator-small} and the choice $k\varepsilon>4$,
	we may increase $N$ so that $3q\le1/2$ and
	$(3q)^k\le 3^kN^{-k\varepsilon/2}\le N^{-2}$. Hence
	\begin{equation*}
		Y(Q_{aN}K^k(Q_bu))
		\le (3q)^kY(Q_bu)
		\lesssim_{M}
		N^{-2}\left(N^{-1-\varepsilon/4}+\frac{\sqrt L}{N}\right)
		\lesssim N^{-1-\varepsilon/2}.
	\end{equation*}
	
	For the first term in \eqref{eq:positive-high-reduction},
	\eqref{K} gives
	\begin{equation*}
		\sum_{j=0}^{k-1}Y\bigl(Q_{aN}K^j[Q_bu-K(Q_Lu)]\bigr)
		\le Y\bigl(Q_{aN}[Q_bu-K(Q_Lu)]\bigr)
		+\sum_{j=1}^{k-1}(3q)^jY\bigl(Q_bu-K(Q_Lu)\bigr).
	\end{equation*}
	Since $3q\le1/2$, the sum of the coefficients $(3q)^j$
	for $j\ge1$ is at most one. It therefore suffices to prove
	\begin{equation}\label{eq:positive-z-star}
		Y\bigl(P[Q_bu-K(Q_Lu)]\bigr)
		\lesssim_{M}N^{-1-\varepsilon/2},
		\qquad P=\operatorname{Id},\,Q_{aN}.
	\end{equation}
	We estimate the right-hand side of \eqref{duhamelh} term by term.
	
	For $S(t)Q_bu(0)$, split $u(0)=P_Nu(0)+Q_Nu(0)$.
	Since $Q_b$ vanishes on $|\xi|\le b$, Plancherel,
	\eqref{energy}, and \eqref{bootstrap} give
	\begin{equation*}
		\norm{Q_bu(0)}{L^2(\R^2)}
		\le b^{-1}\norm{\nabla P_Nu(0)}{L^2(\R^2)}
		+\norm{Q_Nu(0)}{L^2(\R^2)}
		\lesssim N^{-1}.
	\end{equation*}
	Apply \eqref{freehigh} to $(P_H-P_{H/2})PQ_bu(0)$.
	Its $L^2$ norm is at most $\norm{Q_bu(0)}{L^2(\R^2)}$,
	and $\sum_{H\ge b/2}(L/H)^{1/2}\lesssim (L/N)^{1/2}$.
	Thus
	\begin{equation*}
		Y\bigl(S(t)PQ_bu(0)\bigr)
		\lesssim (L/N)^{1/2}N^{-1}.
	\end{equation*}
	
	For $2P_Lu\,|Q_Lu|^2$ and $\overline{P_Lu}(Q_Lu)^2$,
	apply \eqref{eq:positive-annular-bound} with
	$(b_2,G)=(P_Lu,|Q_Lu|^2)$ and
	$(b_2,G)=(\overline{P_Lu},(Q_Lu)^2)$, respectively.
	The $L^2(J\times\R^2)$ norm of $G$ is
	$\norm{Q_Lu}{L^4(J\times\R^2)}^2$ in both cases. Summing over
	$H$ gives
	\begin{equation*}
		Y\left(\int_0^t S(t-\tau)PQ_b
		\bigl(2P_Lu\,|Q_Lu|^2+\overline{P_Lu}(Q_Lu)^2\bigr)(\tau)
		\dd\tau\right)
		\le 3q\norm{Q_Lu}{L^4(J\times\R^2)}^2 .
	\end{equation*}
	For $F(Q_Lu)$, apply \eqref{asd} with
	$m(D)=(P_H-P_{H/2})PQ_b$ and sum over $H$. Since
	$\norm{F(Q_Lu)}{L^{4/3}(J\times\R^2)}
	=\norm{Q_Lu}{L^4(J\times\R^2)}^3$, we obtain
	\begin{equation*}
		Y\left(\int_0^t S(t-\tau)PQ_bF(Q_Lu)(\tau)\dd\tau\right)
		\lesssim_{M}(L/N)^{1/2}\norm{Q_Lu}{L^4(J\times\R^2)}^3 .
	\end{equation*}

	Finally, substituting $L/N=N^{-\varepsilon}$,
	$q\le N^{-\varepsilon/2}$, and
	$\norm{Q_Lu}{L^4(J\times\R^2)}\le N^{-2/3}$ (from \eqref{bootstrap}) into the three
	bounds gives
	\begin{equation*}
		Y\bigl(P[Q_bu-K(Q_Lu)]\bigr)
		\lesssim_{M}N^{-\varepsilon/2}
		\bigl(N^{-1}+N^{-4/3}+N^{-2}\bigr)
		\lesssim_{M,a}N^{-1-\varepsilon/2}.
	\end{equation*}
	This proves \eqref{eq:positive-z-star} and hence
	\eqref{Ximprove} by \eqref{eq:positive-high-reduction}.
	
\end{proof}

Next we derive $L^4(J\times\R^2)$ bound for $Q_Nu$, which will be used in the proof
of \eqref{tail-improved} to control the desired $L^4(J\times\R^2)$ bound for $Q_Lu$ in \eqref{bootstrap}.

\begin{lemma}\label{lem:positive-above-N}
	Under \eqref{bootstrap},
	\begin{equation}\label{w}
		\norm{Q_Nu}{L^4(J\times\R^2)}\lesssim_M N^{-1}.
	\end{equation}
\end{lemma}
\begin{proof}
	Apply $Q_N$ to \eqref{eq:duhamel} and use
	\eqref{eq:positive-cubic-expansion}. For the free solution,
	\eqref{eq:strichartz} and \eqref{bootstrap} give
	\begin{equation*}
		\norm{S(t)Q_Nu(0)}{L^4(J\times\R^2)}
		\lesssim \norm{Q_Nu(0)}{L^2(\R^2)}
		\lesssim N^{-1}.
	\end{equation*}
	Since $6L<N$, we have $Q_NF(P_Lu)=0$. For $F(Q_Lu)$,
	\eqref{eq:duhamel-multiplier-l4} with $m(D)=Q_N$ gives
	\begin{equation*}
		\left\|\int_0^t S(t-\tau)Q_N F(Q_Lu)(\tau)\dd\tau
		\right\|_{L^4(J\times\R^2)}
		\lesssim \norm{Q_Lu}{L^4(J\times\R^2)}^3 .
	\end{equation*}
	
	For the terms containing both $P_Lu$ and $Q_Lu$, use
	$Q_N=\sum_{H\ge N/2}(P_H-P_{H/2})Q_N$. For
	$2P_Lu\,|Q_Lu|^2+\overline{P_Lu}(Q_Lu)^2$, apply
	\eqref{low-input} on each annulus with the two choices
	\begin{equation*}
		(\ell',G)=(P_Lu,|Q_Lu|^2)
		\qquad\text{and}\qquad
		(\ell',G)=(\overline{P_Lu},(Q_Lu)^2).
	\end{equation*}
	Summing $(L/H)^{1/2}$ over $H\ge N/2$ bounds these terms by
	$\lesssim_M(L/N)^{1/2}\norm{Q_Lu}{L^4(J\times\R^2)}^2$.
	
	For the terms linear in $Q_Lu$, write
	$Q_Lu=Q_{aN}u+(P_{aN}-P_L)u$. The products containing
	$(P_{aN}-P_L)u$ have frequencies at most $2aN+4L<N$.
	Consequently,
	\begin{equation*}
		Q_N\bigl(2|P_Lu|^2Q_Lu+(P_Lu)^2\overline{Q_Lu}\bigr)
		=Q_N\bigl(2|P_Lu|^2Q_{aN}u+(P_Lu)^2\overline{Q_{aN}u}\bigr).
	\end{equation*}
	Apply \eqref{low-input} to both terms on the right with
	$\ell'=P_Lu$. The remaining factors satisfy
	\begin{equation*}
		\norm{\overline{P_Lu}\,Q_{aN}u}{L^2(J\times\R^2)}
		=\norm{(P_Lu)\overline{Q_{aN}u}}{L^2(J\times\R^2)}
		\le Y(Q_{aN}u).
	\end{equation*}
	After summation in $H$, the resulting bound is
	$\lesssim_M(L/N)^{1/2}Y(Q_{aN}u)$. Combining the estimates and
	using \eqref{bootstrap}, we conclude that
	\begin{equation*}
		\norm{Q_Nu}{L^4(J\times\R^2)}
		\lesssim_{M}
		N^{-1}+N^{-2}
		+N^{-\varepsilon/2}
		\bigl(N^{-4/3}+N^{-1-\varepsilon/4}\bigr)
		\lesssim N^{-1}.
	\end{equation*}
	
\end{proof}

\vspace{7pt}
\subsection{Improving the $L^4$ bound for $Q_Lu$}\label{sec:middle}

Lemma~\ref{lem:positive-above-N} gives
$\norm{Q_Nu}{L^4(J\times\R^2)}\lesssim_M N^{-1}$. Since
\begin{equation*}
	Q_Lu=Q_LP_Nu+Q_LQ_Nu,
\end{equation*}
it remains to estimate $Q_LP_Nu$. Write $v=P_Nu$ and, for $r>0$,
let $v_{\ge r}$ denote the restriction of $v$ to frequencies
$|\xi|\ge r$:
\begin{equation*}
	\widehat{v_{\ge r}}(t,\xi)
	:=\mathbf1_{\{|\xi|\ge r\}}\widehat v(t,\xi).
\end{equation*}
Since $Q_Lv=Q_Lv_{\ge L}$, it suffices to bound
$\norm{v_{\ge L}}{L^4(J\times\R^2)}$. 
\begin{remark}
	The sharp cutoff $\mathbf1_{\{|\xi|\ge r\}}$ in $v_{\ge r}$ is introduced only to make the next frequency
	separation explicit.  The proof can also be adapted to use the smooth cutoff $Q_r$.
\end{remark}

\begin{lemma}\label{quartet}
	Assume \eqref{bootstrap}, and let $L/2\le\rho<r\le L$.
	For every finite partition of $J$ into intervals $I$ with
	left endpoints $c_I$, there exist $e_I\ge0$ such that
	\begin{equation}\label{eq:positive-projected-cubic}
		\left\|\int_{c_I}^t S(t-\tau)
		\mathbf1_{\{|D|\ge r\}}P_NF(v)(\tau)\dd\tau
		\right\|_{L^4(I\times\R^2)}
		\lesssim \norm{v_{\ge\rho}}{L^4(I\times\R^2)}^3+e_I
	\end{equation}
	and
	\begin{equation}\label{separated}
		\left(\sum_I e_I^2\right)^{1/2}
		\lesssim_M \frac{N^{5/2}}{L(r-\rho)^2}.
	\end{equation}
	The implicit constant in \eqref{eq:positive-projected-cubic} is
	independent of $r,\rho,N,J$ and the partition.
\end{lemma}

\begin{proof}
	To prove \eqref{eq:positive-projected-cubic}, fix $I$ and
	test the Duhamel integral against a function $\phi_I$ of unit
	norm in $L^{4/3}(I\times\R^2)$. Reversing the order of
	integration gives the pairing
	\begin{equation}\label{eq:positive-quartic-pairing}
		\int_I\int_{c_I}^t
		\left\langle S(t-\tau)\mathbf1_{\{|D|\ge r\}}P_N
		(|v|^2v)(\tau),\phi_I(t)\right\rangle_{L^2(\R^2)}
		\dd\tau\dd t
		=\int_I\int_{\R^2}
		v(t,x)\bar v(t,x)v(t,x)\overline{V_I(t,x)}\dd x\dd t ,
	\end{equation}
	where
	\begin{equation*}
		V_I(\tau):=\int_{I\cap\{t>\tau\}}S(\tau-t)
		\mathbf1_{\{|D|\ge r\}}P_N\phi_I(t)\dd t .
	\end{equation*}
	
	The symbol of $\mathbf1_{\{|D|\ge r\}}P_N$ has absolute value
	at most one and is supported in $\{r\le|\xi|\le2N\}$.
	The advanced forms of \eqref{eq:duhamel-multiplier-l4} and
	\eqref{eq:duhamel-gradient-product}, applied to $\phi_I$, give
	\begin{equation}\label{eq:positive-advanced-pair}
		\norm{V_I}{L^4(I\times\R^2)}\lesssim1,\qquad
		\sup_{y\in\R^2}
		\norm{\nabla_x[v(t,x)\overline{V_I(t,x-y)}]}
		{L^2(I\times\R^2)}
		\lesssim_M N^{3/2}L^{-1/2}.
	\end{equation}
	
	We first determine which frequencies in
	\eqref{eq:positive-quartic-pairing} require the gradient
	estimate in \eqref{eq:positive-advanced-pair}. Denote by $\xi_1,\xi_2,\xi_3,\xi_4$ the
	frequencies of $v,\bar v,v,\overline{V_I}$, respectively, measured
	before conjugation. For non-vanishing integral, we have
	\begin{equation*}
		\xi_1-\xi_2+\xi_3-\xi_4=0,\qquad r\le |\xi_4|\le 2N.
	\end{equation*}
	In particular,
	\begin{equation*}
		\xi_2-\xi_4=(\xi_1-\xi_4)+(\xi_3-\xi_4).
	\end{equation*}
	If any of $|\xi_1|,|\xi_2|,|\xi_3|$ is less than $\rho$, then
	for that $j$ we have $|\xi_j-\xi_4|>r-\rho$. Consequently,
	\begin{equation*}
		r-\rho
		<\max_{1\le j\le3}|\xi_j-\xi_4|
		\le |\xi_1-\xi_4|+|\xi_3-\xi_4| .
	\end{equation*}
	This suggests splitting first in $\xi_1-\xi_4$ and then, when
	this difference is small, in $\xi_3-\xi_4$. When both
	differences are small, only $v_{\ge\rho}$ occurs. When either
	difference is large, we keep the factors $v$ so that
	\eqref{vv} and \eqref{eq:positive-advanced-pair} apply directly.
	
	Set $\delta=(r-\rho)/8$. To split in $\xi_1-\xi_4$, group
	the first and fourth factors in \eqref{eq:positive-quartic-pairing}
	and write
	\begin{equation*}
		\int_I\int_{\R^2} v\bar vv\overline{V_I}\dd x\dd t
		=\int_I\int_{\R^2}P_\delta(v\overline{V_I})|v|^2\dd x\dd t
		+\int_I\int_{\R^2}Q_\delta(v\overline{V_I})|v|^2\dd x\dd t .
	\end{equation*}
	On the support
	of $Q_\delta$ we have $|\xi_1-\xi_4|\ge\delta$, and
	$Q_{\delta/2}$ equals one for $|\xi|\ge\delta$. Hence the
	second integral equals
	\begin{equation}\label{eq:positive-first-separated-pair}
		\int_I\int_{\R^2}
		Q_\delta(v\overline{V_I})\,Q_{\delta/2}(|v|^2)\dd x\dd t .
	\end{equation}
	
	In the first integral, write $P_\delta$ as convolution with
	$k_\delta$ and regroup the four factors as follows:
	\begin{equation*}
		\int_I\int_{\R^2}P_\delta(v\overline{V_I})|v|^2\dd x\dd t
		=\int_{\R^2}k_\delta(y)\int_I\int_{\R^2}
		\bigl[v(t,x)\overline{V_I(t,x-y)}\bigr]
		\bigl[v(t,x-y)\bar v(t,x)\bigr]\dd x\dd t\dd y .
	\end{equation*}
	The first bracket now pairs the third and fourth Fourier
	factors, so its frequency is $\xi_3-\xi_4$. Apply
	$\operatorname{Id}=Q_\delta+P_\delta$ to this bracket.
	For the $Q_\delta$ term, the frequency of the second bracket
	satisfies $|\xi_1-\xi_2|=|\xi_3-\xi_4|\ge\delta$.
	Inserting $Q_{\delta/2}$ on the second bracket therefore gives
	\begin{equation}\label{eq:positive-second-separated-pair}
		\int_{\R^2}k_\delta(y)\int_I\int_{\R^2}
		Q_\delta\bigl[v(t,\cdot)\overline{V_I(t,\cdot-y)}\bigr](x)
		\,Q_{\delta/2}\bigl[v(t,\cdot-y)\bar v(t,\cdot)\bigr](x)
		\dd x\dd t\dd y .
	\end{equation}
	For the $P_\delta$ term, the two projections applied above
	give the Fourier multiplier
	\begin{equation*}
		\chi\bigl((\xi_1-\xi_4)/\delta\bigr)
		\chi\bigl((\xi_3-\xi_4)/\delta\bigr).
	\end{equation*}
	On its support, $|\xi_1-\xi_4|,|\xi_3-\xi_4|\le2\delta$,
	whence $|\xi_2-\xi_4|\le4\delta$. It follows that
	\begin{equation*}
		|\xi_j|\ge r-4\delta=(r+\rho)/2>\rho,
		\qquad j=1,2,3 .
	\end{equation*}
	Thus inserting $\mathbf1_{\{|\xi_j|\ge\rho\}}$ on each of
	the three Fourier factors of $v$ leaves this term unchanged.
	Since $\mathbf1_{\{|\xi|\ge\rho\}}\widehat v
	=\widehat{v_{\ge\rho}}$, expanding the remaining $P_\delta$
	as convolution in $z$ gives exactly
	\begin{equation*}
		\iint_{\R^2\times\R^2}k_\delta(y)k_\delta(z)
		\int_I\int_{\R^2}
		v_{\ge\rho}(t,x-y)\overline{v_{\ge\rho}(t,x)}
		v_{\ge\rho}(t,x-z)\overline{V_I(t,x-y-z)}
		\dd x\dd t\dd y\dd z .
	\end{equation*}
	For each $y,z$, apply Hölder's inequality in $(t,x)$ with
	exponents $(4,4,4,4)$. Spatial translations preserve each
	$L^4$ norm, so the absolute value of the displayed integral
	is bounded by
	\begin{equation*}
		\norm{k_\delta}{L^1(\R^2)}^2
		\norm{v_{\ge\rho}}{L^4(I\times\R^2)}^3
		\norm{V_I}{L^4(I\times\R^2)}
		\lesssim \norm{v_{\ge\rho}}{L^4(I\times\R^2)}^3 .
	\end{equation*}
	Here we used \eqref{eq:positive-advanced-pair} and
	$\norm{k_\delta}{L^1(\R^2)}=\norm{k_1}{L^1(\R^2)}$.
	
	It remains to bound \eqref{eq:positive-first-separated-pair}
	and \eqref{eq:positive-second-separated-pair} by an error $e_I$
	whose square sum satisfies \eqref{separated}. For any two
	functions $f,g$ on $I\times\R^2$, Cauchy--Schwarz gives
	\begin{equation*}
		\left|\int_I\int_{\R^2}Q_\delta f\,Q_{\delta/2}g
		\dd x\dd t\right|
		\le \norm{Q_\delta f}{L^2(I\times\R^2)}
		\norm{Q_{\delta/2}g}{L^2(I\times\R^2)}
		\le 2\delta^{-2}
		\norm{\nabla_x f}{L^2(I\times\R^2)}
		\norm{\nabla_x g}{L^2(I\times\R^2)}.
	\end{equation*}
    For
	\eqref{eq:positive-first-separated-pair}, take
	$f=v\overline{V_I}$ and $g=|v|^2$. The absolute value is
	bounded by
	\begin{equation*}
		2\delta^{-2}
		\norm{\nabla_x(v\overline{V_I})}{L^2(I\times\R^2)}
		\norm{\nabla_x|v|^2}{L^2(I\times\R^2)}
		\lesssim_M \delta^{-2}N^{3/2}L^{-1/2}
		\norm{\nabla_x|v|^2}{L^2(I\times\R^2)}
	\end{equation*}
	by \eqref{eq:positive-advanced-pair} at $y=0$.
	For \eqref{eq:positive-second-separated-pair}, apply the same
	inequality at each $y$ to the two products inside the
	projections. This bounds the absolute value by
	\begin{equation*}
		2\delta^{-2}\int_{\R^2}|k_\delta(y)|
		\norm{\nabla_x[v(t,x)\overline{V_I(t,x-y)}]}{L^2(I\times\R^2)}
		\norm{\nabla_x[v(t,x-y)\bar v(t,x)]}{L^2(I\times\R^2)}
		\dd y .
	\end{equation*}
	The first norm under the integral is bounded by
	\eqref{eq:positive-advanced-pair} uniformly in $y$.
	
	We may therefore take
	\begin{equation*}
		e_I:=\delta^{-2}N^{3/2}L^{-1/2}
		\left(
		\norm{\nabla_x|v|^2}{L^2(I\times\R^2)}
		+\int_{\R^2}|k_\delta(y)|
		\norm{\nabla_x[v(t,x-y)\bar v(t,x)]}{L^2(I\times\R^2)}
		\dd y
		\right).
	\end{equation*}
	Combining the three integrals and taking the supremum over
	unit $\phi_I$ gives
	\begin{equation*}
		\left\|\int_{c_I}^t S(t-\tau)\mathbf1_{\{|D|\ge r\}}P_N
		(|v|^2v)(\tau)\dd\tau\right\|_{L^4(I\times\R^2)}
		\lesssim \norm{v_{\ge\rho}}{L^4(I\times\R^2)}^3+e_I ,
	\end{equation*}
	which proves \eqref{eq:positive-projected-cubic}.
	
	Finally, since $v=P_Nu$ and the intervals $I$ are disjoint
	in $J$, \eqref{vv} gives, for every $y\in\R^2$,
	\begin{equation*}
		\left(\sum_I
		\norm{\nabla_x[v(t,x-y)\bar v(t,x)]}{L^2(I\times\R^2)}^2
		\right)^{1/2}
		\le \norm{\nabla_x[v(t,x-y)\bar v(t,x)]}{L^2(J\times\R^2)}
		\lesssim_M NL^{-1/2}.
	\end{equation*}
	At $y=0$ this controls the first term in the definition of
	$e_I$. For the second term, Minkowski's inequality gives
	\begin{equation*}
		\left(\sum_I\left[
		\int_{\R^2}|k_\delta(y)|
		\norm{\nabla_x[v(t,x-y)\bar v(t,x)]}{L^2(I\times\R^2)}
		\dd y\right]^2\right)^{1/2}
	\end{equation*}
	\begin{equation*}
		\le \int_{\R^2}|k_\delta(y)|
		\left(\sum_I
		\norm{\nabla_x[v(t,x-y)\bar v(t,x)]}{L^2(I\times\R^2)}^2
		\right)^{1/2}\dd y \lesssim_M\norm{k_1}{L^1(\R^2)}NL^{-1/2}.
	\end{equation*}
	Substituting both bounds into the definition of $e_I$ and
	using $\delta=(r-\rho)/8$, we conclude that
	\begin{equation*}
		\left(\sum_I e_I^2\right)^{1/2}
		\lesssim_M \delta^{-2}N^{3/2}L^{-1/2}
		\bigl(1+\norm{k_1}{L^1(\R^2)}\bigr)NL^{-1/2}
		\lesssim \frac{N^{5/2}}{L(r-\rho)^2}.
	\end{equation*}
	This proves \eqref{separated}.
	
\end{proof}

By choosing suitable $r$ and $\rho$ in Lemma~\ref{quartet}, we can finally conclude the desired $L^4$-control of $Q_Lu$, which is the last piece of our bootstrap argument.

\begin{proposition}\label{tail}
	Under \eqref{bootstrap}, uniformly over finite $J$,
	\begin{equation}\label{tail-improved}
		\norm{Q_Lu}{L^4(J\times\R^2)}
		\lesssim_M N^{5/4}L^{-2}+N^{-1}
		\lesssim_M N^{-3/4+2\varepsilon}.
	\end{equation}
\end{proposition}

\begin{proof}
	By \eqref{w}, $Q_Lv=Q_Lv_{\ge L}$, and the uniform
	$L^4(\R^2)$ boundedness of $Q_L$, it suffices to prove
	\begin{equation*}
		\norm{v_{\ge L}}{L^4(J\times\R^2)}^4\lesssim_M N^5L^{-8}.
	\end{equation*}
	For this purpose, consider
	\begin{equation}\label{eq:positive-radius-supremum}
		\sup_{L/2<R\le L}
		(R-L/2)^4\norm{v_{\ge R}}{L^4(J\times\R^2)}^4 .
	\end{equation}
	The fourth power of $R-L/2$ is designed to compensate for the factor
	$(r-\rho)^{-4}$ obtained by squaring \eqref{separated}.
	The supremum is finite on each finite $J$. Indeed, for
	$R>L/2$, Plancherel and \eqref{bootstrap} give
	\begin{equation*}
		\norm{v_{\ge R}(t)}{L^2(\R^2)}
		\le 2L^{-1}\norm{\nabla v(t)}{L^2(\R^2)}
		\lesssim L^{-1},\qquad
		\norm{\nabla v_{\ge R}(t)}{L^2(\R^2)}
		\le \norm{\nabla v(t)}{L^2(\R^2)}\lesssim1 .
	\end{equation*}
	Hence \eqref{eq:sobolev-gn} yields
	$\norm{v_{\ge R}}{L^4(J\times\R^2)}^4\lesssim |J|L^{-2}$
	uniformly for $L/2<R\le L$.
	
	If \eqref{eq:positive-radius-supremum} is zero, the conclusion
	follows by taking $R=L$. Otherwise, choose $r\in(L/2,L]$ at
	which the expression is at least half its supremum.
	Hence, to obtain
	$\norm{v_{\ge L}}{L^4(J\times\R^2)}^4\lesssim_M N^{5}L^{-8}$,
	it suffices to prove
	\begin{equation}\label{asdf}
		(r-L/2)^4\norm{v_{\ge r}}{L^4(J\times\R^2)}^4
		\lesssim_M N^5L^{-4}.
	\end{equation}
	
	To estimate $v_{\ge r}$, we will apply
	\eqref{eq:positive-projected-cubic}, whose right-hand side
	involves $v_{\ge\rho}$ at a smaller radius $\rho$.
	Choose $\rho=(r+L/2)/2$, so that
	$r-\rho=\rho-L/2=(r-L/2)/2$. Then, from the choice
	of $r$, we can obtain
	\begin{equation}\label{eq:positive-two-radii}
		\norm{v_{\ge\rho}}{L^4(J\times\R^2)}^4
		\le 32\norm{v_{\ge r}}{L^4(J\times\R^2)}^4 .
	\end{equation}
	
	Fix a finite partition of $J$. On each interval $I$ with left endpoint $c_I$, we have the following Duhamel formula for $v_{\ge r}$
	\begin{equation*}
		v_{\ge r}(t)
		=S(t-c_I)v_{\ge r}(c_I)
		-\ii\int_{c_I}^t S(t-\tau)\mathbf1_{\{|D|\ge r\}}P_N
		F(u(\tau))\dd\tau,
		\qquad t\in I .
	\end{equation*}
	 By the Strichartz estimate \eqref{eq:strichartz},
	\begin{equation*}
		\norm{S(t-c_I)v_{\ge r}(c_I)}{L^4(I\times\R^2)}
		\lesssim \norm{v_{\ge r}(c_I)}{L^2(\R^2)}
		\lesssim L^{-1},
	\end{equation*}
	 For the Duhamel
	integral, Lemma~\ref{quartet} estimates the part containing
	$F(v)$. We next show that the part containing $F(u)-F(v)$
	also satisfies the error bound \eqref{separated}.
	
	Since $P_NP_L=P_L$, we have
	\begin{equation*}
		u=P_Lu+Q_Lu,\qquad
		v=P_Lu+P_NQ_Lu,\qquad Q_Lu-P_NQ_Lu=Q_Nu .
	\end{equation*}
	Expanding $F(u)-F(v)$ as the sum of
	\begin{equation}\label{eq:tail-no-low}
		F(Q_Lu)-F(P_NQ_Lu),
	\end{equation}
	\begin{equation}\label{eq:tail-one-low}
		2P_Lu\bigl(|Q_Lu|^2-|P_NQ_Lu|^2\bigr)
		+\overline{P_Lu}\bigl((Q_Lu)^2-(P_NQ_Lu)^2\bigr),
	\end{equation}
	and
	\begin{equation}\label{eq:tail-two-low}
		2|P_Lu|^2Q_Nu+(P_Lu)^2\overline{Q_Nu}.
	\end{equation}
	
	For \eqref{eq:tail-no-low} and \eqref{eq:tail-one-low}, it
	suffices to bound their $L^{4/3}(J\times\R^2)$ norms. Indeed,
	for either expression $G$, \eqref{eq:duhamel-multiplier-l4}
	and disjointness give
	\begin{equation*}
		\left(\sum_I
		\left\|\int_{c_I}^t S(t-\tau)
		\mathbf1_{\{|D|\ge r\}}P_NG(\tau)\dd\tau
		\right\|_{L^4(I\times\R^2)}^2\right)^{1/2}
		\lesssim
		\left(\sum_I\norm G{L^{4/3}(I\times\R^2)}^{4/3}\right)^{3/4}
		=\norm G{L^{4/3}(J\times\R^2)} .
	\end{equation*}
	For \eqref{eq:tail-no-low}, use
	$|F(z)-F(w)|\lesssim |z-w|(|z|^2+|w|^2)$ and
	$Q_Lu-P_NQ_Lu=Q_Nu$. H\"older, the $L^4$ boundedness of
	$P_N$, and \eqref{w} give
	\begin{equation*}
		\norm{F(Q_Lu)-F(P_NQ_Lu)}{L^{4/3}(J\times\R^2)}
		\lesssim \norm{Q_Nu}{L^4(J\times\R^2)}
		\norm{Q_Lu}{L^4(J\times\R^2)}^2
		\lesssim_M N^{-1}\norm{Q_Lu}{L^4(J\times\R^2)}^2 .
	\end{equation*}
	
	For \eqref{eq:tail-one-low}, factor the two differences as
	\begin{equation*}
		|Q_Lu|^2-|P_NQ_Lu|^2
		=(Q_Nu)\overline{Q_Lu}+(P_NQ_Lu)\overline{Q_Nu},
	\end{equation*}
	\begin{equation*}
		(Q_Lu)^2-(P_NQ_Lu)^2=(Q_Nu)(Q_Lu+P_NQ_Lu).
	\end{equation*}
	Since $Q_N=Q_NQ_{aN}=Q_{aN}-P_N Q_{aN}$,  \eqref{eq:translation-kernel} gives
	\begin{equation}\label{eq:tail-low-high}
		\norm{(P_Lu)(Q_Nu)}{L^2(J\times\R^2)}
		\le \bigl(1+\norm{k_N}{L^1(\R^2)}\bigr)Y(Q_{aN}u).
	\end{equation}
	H\"older with $3/4=1/2+1/4$ therefore bounds the
	$L^{4/3}(J\times\R^2)$ norm of \eqref{eq:tail-one-low} by
	\begin{equation*}
		\norm{(P_Lu)(Q_Nu)}{L^2(J\times\R^2)}
		\bigl(\norm{Q_Lu}{L^4(J\times\R^2)}
		+\norm{P_NQ_Lu}{L^4(J\times\R^2)}\bigr)
		\lesssim Y(Q_{aN}u)\norm{Q_Lu}{L^4(J\times\R^2)} .
	\end{equation*}
	
	For \eqref{eq:tail-two-low}, both products have frequencies
	$|\xi|\ge N-4L\ge N/2$ for large $N$. Since $P_N$ vanishes
	above $2N$, we may restrict the symbol of
	$\mathbf1_{\{|D|\ge r\}}P_N$ to $\{N/2\le|\xi|\le2N\}$
	when applying it to these products. Use \eqref{low-input}
	with $H=N$, $t_0=c_I$, and $\ell'=P_Lu$. The two choices of
	$G$ are $\overline{P_Lu}\,Q_Nu$ and
	$(P_Lu)\overline{Q_Nu}$. Hence the square sum of the
	corresponding Duhamel norms, after taking its square root,
	is bounded by
	\begin{equation*}
		(L/N)^{1/2}
		\left(\sum_I\norm{(P_Lu)(Q_Nu)}{L^2(I\times\R^2)}^2
		\right)^{1/2}
		\lesssim_M (L/N)^{1/2}Y(Q_{aN}u),
	\end{equation*}
	where the last inequality follows from \eqref{eq:tail-low-high}.
	
Increase $e_I$ from Lemma~\ref{quartet} by the sum of the three
Duhamel norms on $I$ estimated in \eqref{eq:tail-no-low},
\eqref{eq:tail-one-low}, and \eqref{eq:tail-two-low}. Combining
these estimates with \eqref{eq:positive-projected-cubic} and the
bound for the free solution gives
\begin{equation}\label{countinput}
	\norm{v_{\ge r}}{L^4(I\times\R^2)}
	\lesssim_M L^{-1}+\norm{v_{\ge\rho}}{L^4(I\times\R^2)}^3+e_I .
\end{equation}
By \eqref{bootstrap}, the three additional contributions to
$(\sum_I e_I^2)^{1/2}$ are bounded, up to constants depending
only on $M$, by
\begin{equation*}
	N^{-7/3},\qquad N^{-5/3-\varepsilon/4},\qquad
	N^{-1-3\varepsilon/4},
\end{equation*}
respectively. Each is at most $N^{5/2}L^{-3}=N^{-1/2+3\varepsilon}$.
Since $r-\rho\le L$, 
\eqref{separated} gives
\begin{equation}\label{eq:positive-count-errors}
	\sum_I e_I^2\lesssim_M \frac{N^5}{L^2(r-\rho)^4}.
\end{equation}

We now choose the partition so that the term $L^{-1}$ in
\eqref{countinput} can be absorbed on every interval. Starting at the left endpoint of $J$, take
successive intervals on which $\norm{v_{\ge r}}{L^4(I\times\R^2)}$
equals a fixed multiple of $L^{-1}$, with at most one final
interval of smaller norm. Choose the multiple to be at least
twice the implicit constant in \eqref{countinput}. 

On each interval that attains the chosen norm, subtracting the
contribution of $L^{-1}$ in \eqref{countinput} gives
\begin{equation*}
	\norm{v_{\ge\rho}}{L^4(I\times\R^2)}^3\gtrsim_M L^{-1}
	\qquad\text{or}\qquad
	e_I\gtrsim_M L^{-1}.
\end{equation*}
In the first case,
$1\lesssim_M L^{4/3}\norm{v_{\ge\rho}}{L^4(I\times\R^2)}^4$.
In the second case, $1\lesssim_M L^2e_I^2$. Summing over the
disjoint intervals and including the possible final interval gives
\begin{equation*}
	\#\{I\}\lesssim_M
	1+L^{4/3}\norm{v_{\ge\rho}}{L^4(J\times\R^2)}^4
	+L^2\sum_I e_I^2 .
\end{equation*}
Since $\norm{v_{\ge r}}{L^4(I\times\R^2)}^4\lesssim_M L^{-4}$
on every interval, we have
\begin{equation*}
	\norm{v_{\ge r}}{L^4(J\times\R^2)}^4
	=\sum_I\norm{v_{\ge r}}{L^4(I\times\R^2)}^4
	\lesssim_M L^{-4}\#\{I\}.
\end{equation*}
Substituting the bound for $\#\{I\}$ and then
\eqref{eq:positive-count-errors} yields
\begin{equation}\label{count}
	\norm{v_{\ge r}}{L^4(J\times\R^2)}^4
	\lesssim_M L^{-4}
	+L^{-8/3}\norm{v_{\ge\rho}}{L^4(J\times\R^2)}^4
	+\frac{N^5}{L^4(r-\rho)^4}.
\end{equation}

Recall that \eqref{eq:positive-two-radii} gives
\begin{equation*}
	L^{-8/3}\norm{v_{\ge\rho}}{L^4(J\times\R^2)}^4
	\le 32L^{-8/3}\norm{v_{\ge r}}{L^4(J\times\R^2)}^4 .
\end{equation*}
Choose $N$ sufficiently large that $32L^{-8/3}$ times the
implicit constant in \eqref{count} is at most $1/2$.  Absorbing the middle term in
\eqref{count}, multiplying by $(r-L/2)^4$, and using
$r-\rho=(r-L/2)/2$, we obtain 
\begin{equation*}
	(r-L/2)^4\norm{v_{\ge r}}{L^4(J\times\R^2)}^4
	\lesssim_M (r-L/2)^4L^{-4}+N^5L^{-4}
	\lesssim_M N^5L^{-4},
\end{equation*}
which is exactly the desired (\ref{asdf}). Thus, we finally conclude
\eqref{tail-improved}.

\end{proof}

\vspace{7pt}
\subsection{Scaling and continuation in positive regularity}

We use \eqref{variation}, \eqref{tail-improved}, and \eqref{Ximprove}
to close \eqref{bootstrap} after rescaling. 

\begin{proof}[Proof of Theorem~\ref{thm:positive}]
	For smooth data, conservation of mass and of the defocusing
	energy bounds $\norm{u(a)}{H^1(\R^2)}$ uniformly in $a$
	throughout the lifespan. By \eqref{eq:sobolev-gn} and unitarity,
	\begin{equation*}
		\norm{S(t-a)u(a)}{L^4([a,a+\delta]\times\R^2)}^4
		\lesssim \delta\,\norm{u(a)}{L^2(\R^2)}^2
		\norm{\nabla u(a)}{L^2(\R^2)}^2 .
	\end{equation*}
	Thus \eqref{eq:local-small} holds for some $\delta>0$ independent
	of $a$. Iterating Lemma~\ref{lem:local}, with persistence of
	regularity, gives a smooth solution
	on $\R$. 
	
	Fix any $R,T>0$. It suffices to prove that, for $0<s<1$, every smooth
	solution with $\norm{u(0)}{H^s(\R^2)}\le R$ satisfies
	\begin{equation}\label{eq:positive-uniform-spacetime}
		\norm u{L^4([0,T]\times\R^2)}\lesssim_{s,R,T}1 .
	\end{equation}
   Indeed, if \eqref{eq:positive-uniform-spacetime} holds,
   Lemma~\ref{lem:local} gives existence, uniqueness, and continuous
   dependence of $H^s$ solutions on initial data by smooth approximation, which yields Theorem \ref{thm:positive}.
	
	To prove \eqref{eq:positive-uniform-spacetime}, we write $f=u(0)$ and set
	\begin{equation*}
		u_\lambda(t,x):=\lambda^{-1}u(t/\lambda^2,x/\lambda),
		\qquad f_\lambda(x):=\lambda^{-1}f(x/\lambda).
	\end{equation*}
	Then $u_\lambda$ solves the same equation with initial datum
	$f_\lambda$, and
	\begin{equation*}
		\norm{f_\lambda}{L^2(\R^2)}=\norm f{L^2(\R^2)},\qquad
		\norm{u_\lambda}{L^4([0,\lambda^2T]\times\R^2)}
		=\norm u{L^4([0,T]\times\R^2)} .
	\end{equation*}
	To start \eqref{bootstrap}, we will choose $\lambda$ so that
	$\mathcal E_N(f_\lambda)\le1/4$. The scaling $f_\lambda(x)=\lambda^{-1}f(x/\lambda)$ gives
	\begin{equation*}
		P_Nf_\lambda(x)=\lambda^{-1}(P_{\lambda N}f)(x/\lambda),
		\qquad
		Q_Nf_\lambda(x)=\lambda^{-1}(Q_{\lambda N}f)(x/\lambda).
	\end{equation*}
	Differentiating the first identity and changing variables
	$x=\lambda y$ in both $L^2(\R^2)$ norms, we obtain
	\begin{equation*}
		\norm{\nabla P_Nf_\lambda}{L^2(\R^2)}^2
		+N^2\norm{Q_Nf_\lambda}{L^2(\R^2)}^2
		=\lambda^{-2}\left(
		\norm{\nabla P_{\lambda N}f}{L^2(\R^2)}^2
		+(\lambda N)^2\norm{Q_{\lambda N}f}{L^2(\R^2)}^2
		\right)\lesssim_s N^{2(1-s)}\lambda^{-2s}
		\norm f{\dot H^s(\R^2)}^2.
	\end{equation*}

	For the quartic part, \eqref{eq:sobolev-gn} and mass conservation give
	\begin{equation*}
		\norm{P_Nf_\lambda}{L^4(\R^2)}^4
		\lesssim \norm{P_Nf_\lambda}{L^2(\R^2)}^2
		\norm{\nabla P_Nf_\lambda}{L^2(\R^2)}^2
		\le R^2\norm{\nabla P_Nf_\lambda}{L^2(\R^2)}^2 .
	\end{equation*}
	Consequently,
	\begin{equation*}
		\mathcal E_N(f_\lambda)
		\lesssim_s (1+R^2)R^2N^{2(1-s)}\lambda^{-2s}.
	\end{equation*}
	Take $\lambda$ to be a sufficiently large multiple of
	$N^{(1-s)/s}$, depending only on $s,R$, so that the last
	bound is at most $1/4$. Fix $p=1+2(1-s)/s$. Since
	$\lambda^2T/N^p\lesssim_{s,R}T/N$, we have
	$\lambda^2T\le N^p$ for $N$ sufficiently large. From the discussion in Remark \ref{J}, the assumption \eqref{operator-small} therefore holds
	on every $J=[0,t_1]\subset[0,\lambda^2T]$.
	
	Next, we show that \eqref{bootstrap} holds for $u_\lambda$ on
	$[0,\lambda^2T]$ by using standard bootstrap argument:
	
	 For $J=[0,t_1]$, the three quantities in
	\eqref{bootstrap} are continuous in $t_1$, with $Y$ taken on
	$J$ and formed with $P_Lu_\lambda$. At $t_1=0$, the energy is
	$\mathcal E_N(f_\lambda)\le1/4$ and both spacetime norms vanish.
	Hence \eqref{bootstrap} holds on $[0,t_1]$ for all sufficiently
	small $t_1>0$.
	
    Let $t^*$ be the supremum of all $t_1\in[0,\lambda^2T]$ for which
    \eqref{bootstrap} holds on $[0,t_1]$. The preceding argument gives
    $t^*>0$, and continuity gives \eqref{bootstrap} on $J=[0,t^*]$.
    
    On this interval, \eqref{variation} gives
    \begin{equation*}
    	\sup_{t\in J}
    	|\mathcal E_N(u_\lambda(t))-\mathcal E_N(f_\lambda)|
    	\lesssim_R N^{-\varepsilon/2},
    \end{equation*}
    while \eqref{tail-improved} and \eqref{Ximprove} give, respectively,
    \begin{equation*}
    	N^{2/3}\norm{Q_Lu_\lambda}{L^4(J\times\R^2)}
    	\lesssim_R N^{-(1/12-2\varepsilon)},\qquad
    	N^{1+\varepsilon/4}Y(Q_{aN}u_\lambda)
    	\lesssim_R N^{-\varepsilon/4}.
    \end{equation*}
    Since $0<\varepsilon<1/24$, for all sufficiently large $N$ the
    energy variation is at most $1/4$ and the two rescaled norms are
    each at most $1/2$, uniformly in $t^*$ and in the initial datum
    with $\norm f{H^s(\R^2)}\le R$. Using
    $\mathcal E_N(f_\lambda)\le1/4$, we obtain
    \begin{equation*}
    	\sup_{t\in J}\mathcal E_N(u_\lambda(t))\le\tfrac12,\qquad
    	\norm{Q_Lu_\lambda}{L^4(J\times\R^2)}
    	\le\tfrac12N^{-2/3},\qquad
    	Y(Q_{aN}u_\lambda)\le\tfrac12N^{-1-\varepsilon/4}.
    \end{equation*}
    If $t^*<\lambda^2T$, these strict improvements and continuity
    would give \eqref{bootstrap} on $[0,t^*+\delta]$ for some
    $0<\delta\le\lambda^2T-t^*$, contradicting the definition of
    $t^*$. Hence $t^*=\lambda^2T$, and \eqref{bootstrap} holds on
    $[0,\lambda^2T]$.

	On $[0,\lambda^2T]$, the quartic term in \eqref{energy} gives
	\begin{equation*}
		\sup_{0\le t\le\lambda^2T}\norm{P_Nu_\lambda(t)}{L^4(\R^2)}^4
		\lesssim \sup_{0\le t\le\lambda^2T}\mathcal E_N(u_\lambda(t))
		\lesssim 1 .
	\end{equation*}
	Decomposing $u_\lambda=P_Nu_\lambda+Q_Nu_\lambda$ and applying
	\eqref{w} to the second term, we obtain
	\begin{equation*}
		\norm{u_\lambda}{L^4([0,\lambda^2T]\times\R^2)}
		\le (\lambda^2T)^{1/4}
		\sup_{0\le t\le\lambda^2T}\norm{P_Nu_\lambda(t)}{L^4(\R^2)}
		+\norm{Q_Nu_\lambda}{L^4([0,\lambda^2T]\times\R^2)}
		\lesssim_R (\lambda^2T)^{1/4}+N^{-1}.
	\end{equation*}
	Since $N$ and $\lambda$ were chosen in terms of $s,R,T$, scaling
	back gives
	\begin{equation*}
		\norm u{L^4([0,T]\times\R^2)}
		=\norm{u_\lambda}{L^4([0,\lambda^2T]\times\R^2)}
		\lesssim_{s,R,T}1 .
	\end{equation*}
	This proves \eqref{eq:positive-uniform-spacetime}.

\end{proof}

\vspace{10pt}
\section{Removing the time loss at the mass-critical endpoint}\label{sec:endpoint}

In the proof of Theorem~\ref{thm:positive}, the assumption $s>0$
was used to choose $\lambda\sim_{s,R}N^{(1-s)/s}$ so that
$\mathcal E_N(f_\lambda)\le1/4$. This choice gave
$\lambda^2T\le N^p$ for some fixed $p$, which allowed us to
deduce \eqref{operator-small} from \eqref{poslog}. For
$f\in L^2(\R^2)$, we can still make $\mathcal E_N(f_\lambda)$
small by fixing $N$ and taking $\lambda$ sufficiently large,
but we no longer have a polynomial bound on $\lambda^2T$.

It therefore suffices to prove \eqref{operator-small} independently
of $|J|$. Indeed, the bounds \eqref{variation} and
\eqref{tail-improved} are already uniform in $|J|$, and
Proposition~\ref{high} gives \eqref{Ximprove} once
\eqref{operator-small} holds. 

\begin{lemma}\label{forced-test}
	Assume \eqref{bootstrap} on a finite interval $J$, and fix
	$H\ge cN$ with $c>0$ fixed. Let $m(D)$ have smooth symbol
	supported in $\{H/2\le|\xi|\le2H\}$ and convolution kernel
	with uniformly bounded $L^1(\R^2)$ norm. Let $t_0$ be either
	endpoint of $J$, and suppose that
	\begin{equation*}
		\left\{
		\begin{array}{l}
			(\ii\partial_t+\Delta)h=m(D)\bigl((P_Lu)G\bigr),\\
			h(t_0)=f,
		\end{array}
		\right.
		\qquad
		\supp\widehat f\subset\{H/2\le|\xi|\le2H\},
		\qquad G\in L^2(J\times\R^2).
	\end{equation*}
	Then, for all sufficiently large $N$,
	\begin{equation}\label{forced}
		Y(h)\lesssim_M
		(L/H)^{1/2}\norm f{L^2(\R^2)}
		+(L/H)\norm G{L^2(J\times\R^2)} .
	\end{equation}
	The same estimate holds if $P_Lu$ in the equation is replaced
	by any spatial translate of $P_Lu$ or $\overline{P_Lu}$,
	uniformly in the translation, $J$, and $t_0$.
\end{lemma}

\begin{proof}
	Denote $U=P_Nu$, $V=h$ and write
	\begin{equation*}
		B(t,x,y)=\overline{h(t,y)}\,\nabla U(t,x)
		+U(t,x)\nabla\overline{h(t,y)} .
	\end{equation*}
 For each $d\in\R^2$,
the product $(P_Lu)(t,y+d)\overline{h(t,y)}$ has Fourier
support in $\{|\xi|\ge H/2-2L\}\subset\{|\xi|\ge H/4\}$.
Hence Plancherel gives
\begin{equation*}
	\norm{(P_Lu)(t,y+d)h(t,y)}{L^2(J\times\R^2)}^2
	\le \frac{16}{H^2}
	\norm{\nabla_y\bigl[(P_Lu)(t,y+d)\overline{h(t,y)}\bigr]}
	{L^2(J\times\R^2)}^2 .
\end{equation*}

Since $2L\le N$, we have $P_LU=P_LP_Nu=P_Lu$. Writing $k_L$
for the kernel of $P_L$, the product rule and the definition of
$B$ give
\begin{equation*}
	\nabla_y\bigl[(P_Lu)(t,y+d)\overline{h(t,y)}\bigr]
	=\int_{\R^2}k_L(y-x)B(t,x+d,y)\dd x .
\end{equation*}
By Cauchy--Schwarz,
\begin{equation*}
	\left|\int_{\R^2}k_L(y-x)B(t,x+d,y)\dd x\right|^2
	\le \norm{k_L}{L^1(\R^2)}
	\int_{\R^2}|k_L(y-x)|\,|B(t,x+d,y)|^2\dd x .
\end{equation*}
Now $k_L(z)=L^2k_1(Lz)$ with $k_1$ Schwartz, so
$\norm{k_L}{L^1(\R^2)}=\norm{k_1}{L^1(\R^2)}$ and
$|k_L(y-x)|\lesssim L^2(1+L|x-y|)^{-3}$.
Integrating in $y,t$, substituting into the Plancherel bound,
and taking the supremum over $d$, we obtain
\begin{equation}\label{eq:endpoint-product-from-gradient}
	Y(h)^2\lesssim \frac{L}{H^2}\sup_{d\in\R^2}
	\int_J\iint L(1+L|x-y|)^{-3}
	|B(t,x+d,y)|^2\dd x\dd y\dd t .
\end{equation}
	
	It remains to estimate the integral in
	\eqref{eq:endpoint-product-from-gradient}. We write the
	calculation for $d=0$. Integrating
	$\frac{\dd\mathcal M}{\dd t}=I_1+I_2+I_3$ and using
	\eqref{eq:radial-positive} gives
	\begin{equation*}
		\int_J\iint L(1+L|x-y|)^{-3}|B(t,x,y)|^2\dd x\dd y\dd t
		\lesssim 2\sup_{t\in J}|\mathcal M(t)|
		+\int_J|I_2(t)|\dd t+\int_J|I_3(t)|\dd t .
	\end{equation*}
	
	For $\mathcal M$ and $I_2$, the estimates already proved in
	\eqref{eq:radial-action-mass-momentum} and
	\eqref{eq:radial-source-bound}, together with \eqref{variation},
	give
	\begin{equation*}
		\sup_{t\in J}|\mathcal M(t)|
		\lesssim_M \sup_{t\in J}
		\left(\norm{\rho_h(t)}{L^1(\R^2)}
		+\norm{p_h(t)}{L^1(\R^2)}\right)
	\end{equation*}
	and
	\begin{equation*}
		\int_J|I_2(t)|\dd t
		\lesssim_M N^{-2-\varepsilon/2}
		\left(N\sup_{t\in J}\norm{\rho_h(t)}{L^1(\R^2)}
		+\sup_{t\in J}\norm{p_h(t)}{L^1(\R^2)}\right).
	\end{equation*}
	Since $\rho_h=|h|^2$, $p_h=\im(\bar h\nabla h)$, and
	$\norm{\nabla h(t)}{L^2(\R^2)}
	\le2H\norm{h(t)}{L^2(\R^2)}$, it follows that
	\begin{equation*}
		\sup_{t\in J}|\mathcal M(t)|+\int_J|I_2(t)|\dd t
		\lesssim_M H\sup_{t\in J}\norm{h(t)}{L^2(\R^2)}^2 .
	\end{equation*}
	Thus the two quantities still to be estimated are
	$\sup_{t\in J}\norm{h(t)}{L^2(\R^2)}^2$ and
	$\int_J|I_3(t)|\dd t$.
	
	Let $k$ be the convolution kernel of $m(D)$. For every
	$y,z\in\R^2$, Cauchy--Schwarz and \eqref{Y} give
	\begin{equation*}
		\int_J\int_{\R^2}
		|h(t,x+y+z)|\,|(P_Lu)(t,x)|\,|G(t,x)|\dd x\dd t
		\le Y(h)\norm G{L^2(J\times\R^2)} .
	\end{equation*}
	Expanding $R_h=k*((P_Lu)G)$ and integrating this inequality
	against $|k(z)|\dd z$ yields
	\begin{equation}\label{eq:endpoint-forcing-product}
		\sup_{y\in\R^2}\int_J\int_{\R^2}
		|h(t,x+y)|\,|R_h(t,x)|\dd x\dd t
		\lesssim Y(h)\norm G{L^2(J\times\R^2)} .
	\end{equation}
	Integrating the identity \eqref{masslaw} gives
	\begin{equation*}
		\norm{h(t)}{L^2(\R^2)}^2-\norm f{L^2(\R^2)}^2
		=2\int_{t_0}^t\int_{\R^2}\im(\bar hR_h)\dd x\dd\tau .
	\end{equation*}
	By \eqref{eq:endpoint-forcing-product} at $y=0$, we obtain
	\begin{equation}\label{testmass}
		\sup_{t\in J}\norm{h(t)}{L^2(\R^2)}^2
		\lesssim \norm f{L^2(\R^2)}^2
		+Y(h)\norm G{L^2(J\times\R^2)} .
	\end{equation}
	
	For $I_3$, formula \eqref{testsource} reads
	\begin{equation*}
		I_3(t)=\iint\nabla A(x-y)\cdot
		\left[
		2p_U(t,x)\im(\bar hR_h)(t,y)
		-|U(t,x)|^2
		\re(\overline{R_h}\nabla h-\bar h\nabla R_h)(t,y)
		\right]\dd x\dd y .
	\end{equation*}
	Mass conservation and \eqref{bootstrap} give
	$\norm U{L^2(\R^2)}^2+\norm{p_U}{L^1(\R^2)}\lesssim_M1$.
	Using $\norm{\nabla A}{L^\infty(\R^2)}\le1$ and integrating
	first in $x$, we obtain
	\begin{equation*}
		|I_3(t)|\lesssim_M\int_{\R^2}
		\bigl(|h|\,|R_h|+|\nabla h|\,|R_h|
		+|h|\,|\nabla R_h|\bigr)(t,y)\dd y .
	\end{equation*}
	Both $h$ and $R_h$ have frequencies $|\xi|\le2H$, so
	$\nabla h=(\nabla k_{2H})*h$ and
	$\nabla R_h=(\nabla k_{2H})*R_h$. Expanding either convolution
	and using $\norm{\nabla k_{2H}}{L^1(\R^2)}\lesssim H$ gives
	\begin{equation*}
		\int_J\int_{\R^2}
		\bigl(|\nabla h|\,|R_h|+|h|\,|\nabla R_h|\bigr)\dd x\dd t
		\lesssim H\sup_{y\in\R^2}\int_J\int_{\R^2}
		|h(t,x+y)|\,|R_h(t,x)|\dd x\dd t .
	\end{equation*}
	Hence \eqref{eq:endpoint-forcing-product} and $H\ge1$ yield
	\begin{equation*}
		\int_J|I_3(t)|\dd t
		\lesssim_M H Y(h)\norm G{L^2(J\times\R^2)} .
	\end{equation*}
	
	Combining the bounds for $\mathcal M,I_2,I_3$ with
	\eqref{testmass} gives
	\begin{equation}\label{eq:endpoint-weighted-gradient}
		\int_J\iint L(1+L|x-y|)^{-3}|B(t,x,y)|^2\dd x\dd y\dd t
		\lesssim_M H\left(\norm f{L^2(\R^2)}^2
		+Y(h)\norm G{L^2(J\times\R^2)}\right).
	\end{equation}
	This bound is uniform under translation of $U$, by the
	translation uniformity of the estimates for $\mathcal M,I_2$
	and the $L^1(\R^2)$ bounds for $|U|^2,p_U$ used in $I_3$.
	Substituting \eqref{eq:endpoint-weighted-gradient} into
	\eqref{eq:endpoint-product-from-gradient} gives
	\begin{equation*}
		Y(h)^2
		\lesssim_M \frac LH\left(\norm f{L^2(\R^2)}^2
		+Y(h)\norm G{L^2(J\times\R^2)}\right),
	\end{equation*}
    which gives \eqref{forced}.

\end{proof}

To verify \eqref{operator-small}, take $f=0$ in \eqref{forced}
and apply it with
\begin{equation*}
	m(D)=(P_H-P_{H/2})P Q_b,\qquad
	P\in\{\operatorname{Id},Q_{aN}\},\qquad H\ge b/2 .
\end{equation*}
The symbol is supported in $\{H/2\le|\xi|\le2H\}$. Expanding
$Q_b=\operatorname{Id}-P_b$ and
$Q_{aN}=\operatorname{Id}-P_{aN}$ bounds the $L^1(\R^2)$ norm
of its kernel by
\begin{equation*}
	\norm{k_H-k_{H/2}}{L^1(\R^2)}
	\bigl(1+\norm{k_b}{L^1(\R^2)}\bigr)
	\bigl(1+\norm{k_{aN}}{L^1(\R^2)}\bigr)
	\lesssim 1 .
\end{equation*}
Thus \eqref{forced} gives $d_H\lesssim_M L/H$. Since $b=aN/4$ and the sum is
dyadic,
\begin{equation*}
	\sum_{H\ge b/2}d_H
	\lesssim_M\sum_{H\ge b/2}\frac LH
	\lesssim\frac LN
	=N^{-\varepsilon}.
\end{equation*}
 This proves
\eqref{operator-small}, and Proposition~\ref{high} gives
\eqref{Ximprove} on every finite interval satisfying
\eqref{bootstrap}.

\vspace{10pt}
\section{Global continuation and scattering}\label{sec:global}

\begin{proof}[Proof of Theorem~\ref{thm:endpoint}]
	
\textit{Scaling and global continuation.}
Write $f=u(0)\in L^2(\R^2)$ for the prescribed initial datum.
By Lemma~\ref{lem:local} at $s=0$, it suffices to find
$\delta>0$ such that, for every $T>0$, the smooth solutions
$\widetilde u$ of \eqref{eq:nls} satisfy
\begin{equation*}
	\sup_{\norm{\widetilde u(0)-f}{L^2(\R^2)}<\delta}
	\norm{\widetilde u}{L^4([-T,T]\times\R^2)}<\infty .
\end{equation*}

Since \eqref{operator-small}, \eqref{variation},
\eqref{tail-improved}, and \eqref{Ximprove} are uniform for data
with a common mass bound $M$, we may choose $N$ sufficiently large
uniformly for $\norm{\widetilde u(0)-f}{L^2(\R^2)}<1$.

With $N$ fixed, set $f_\lambda(x)=\lambda^{-1}f(x/\lambda)$.
By scaling,
\begin{equation*}
	\norm{\nabla P_Nf_\lambda}{L^2(\R^2)}^2
	=\lambda^{-2}\norm{\nabla P_{\lambda N}f}{L^2(\R^2)}^2,
	\qquad
	\norm{Q_Nf_\lambda}{L^2(\R^2)}^2
	=\norm{Q_{\lambda N}f}{L^2(\R^2)}^2 .
\end{equation*}
The multipliers $\lambda^{-1}|\xi|\chi(\xi/(\lambda N))$
and $1-\chi(\xi/(\lambda N))$ tend pointwise to zero as
$\lambda\to\infty$ and are uniformly bounded. Thus both terms tend to zero by Plancherel
and dominated convergence. For the quartic term,
\eqref{eq:sobolev-gn} gives
\begin{equation*}
	\norm{P_Nf_\lambda}{L^4(\R^2)}^4
	\lesssim \norm f{L^2(\R^2)}^2
	\norm{\nabla P_Nf_\lambda}{L^2(\R^2)}^2\to 0 .
\end{equation*}
Fix $\lambda$ so that $\mathcal E_N(f_\lambda)<1/8$.
For fixed $N$, the operators $\nabla P_N$ and $NQ_N$ are
bounded on $L^2(\R^2)$, and Bernstein gives
$P_N:L^2(\R^2)\to L^4(\R^2)$. Hence $\mathcal E_N$ is
continuous on $L^2(\R^2)$. Since scaling preserves the $L^2$
norm, there exists $0<\delta\le1$ such that
\begin{equation*}
	\mathcal E_N(\widetilde u_\lambda(0))<\tfrac14,
	\qquad
	\widetilde u_\lambda(t,x)
	=\lambda^{-1}\widetilde u(t/\lambda^2,x/\lambda).
\end{equation*}
Then by \eqref{variation},
\eqref{tail-improved}, \eqref{Ximprove} and bootstrap argument as before, \eqref{bootstrap} holds for $\widetilde{u}_\lambda$ on $[0,\lambda^2T]$. 

By scaling, the quartic term in \eqref{energy}, and \eqref{w},
\begin{equation*}
	\norm{\widetilde u_\lambda}{L^4([0,\lambda^2T]\times\R^2)}
	\le (\lambda^2T)^{1/4}
	\sup_{0\le t\le\lambda^2T}
	\norm{P_N\widetilde u_\lambda(t)}{L^4(\R^2)}
	+\norm{Q_N\widetilde u_\lambda}{L^4([0,\lambda^2T]\times\R^2)}
	\lesssim_M (\lambda^2T)^{1/4}+N^{-1}.
\end{equation*}
Time reversal and conjugation give the same bound on $[-T,0]$.
Since $N$ and $\lambda$ were fixed for the entire ball, we obtain
\begin{equation}\label{localglobal}
	\sup_{\norm{\widetilde u(0)-f}{L^2(\R^2)}<\delta}
	\norm{\widetilde u}{L^4([-T,T]\times\R^2)}
	\lesssim_M (\lambda^2T)^{1/4}+N^{-1}
	\lesssim_{f,T}1 .
\end{equation}
By \eqref{localglobal} and Lemma~\ref{lem:local} at $s=0$,
there is a unique solution with initial datum $f$ on $[-T,T]$,
with continuous dependence in
$C([-T,T];L^2(\R^2))\cap L^4([-T,T]\times\R^2)$.
As $T$ is arbitrary, uniqueness on overlapping intervals gives
\begin{equation*}
	u\in C(\R;L^2(\R^2))\cap L^4_{\mathrm{loc}}(\R\times\R^2),
	\qquad u(0)=f .
\end{equation*}
	
\textit{Scattering.}
To obtain the $L^4$ bound on the two infinite time intervals,
we use the pseudoconformal transformation in
\cite[Section~3]{TaoPC}. For $0<|\tau|\le1$, set
\begin{equation}\label{pseudo}
	U(\tau,y):=|\tau|^{-1}e^{\ii|y|^2/(4\tau)}u(-1/\tau,y/\tau).
\end{equation}
For smooth solutions, differentiation gives
\begin{equation*}
	(\ii\partial_\tau+\Delta_y)U(\tau,y)
	=|\tau|^{-3}e^{\ii|y|^2/(4\tau)}
	\bigl[(\ii\partial_t+\Delta_x)u\bigr](-1/\tau,y/\tau)
	=|U(\tau,y)|^2U(\tau,y).
\end{equation*}
By \eqref{pseudo}, the data at $\tau=-1$ and $\tau=1$ are
\begin{equation*}
	U(-1,y)=e^{-\ii|y|^2/4}u(1,-y),\qquad
	U(1,y)=e^{\ii|y|^2/4}u(-1,y).
\end{equation*}
Both belong to $L^2(\R^2)$. By global existence, we may solve
\eqref{eq:nls} forward from $U(-1)$ to $\tau=0$ and backward
from $U(1)$ to $\tau=0$. By local uniqueness, these solutions
agree with \eqref{pseudo} on $[-1,0)$ and $(0,1]$, respectively. Changing variables by $t=-1/\tau$, $x=y/\tau$ gives
\begin{equation*}
	\int_1^\infty\int_{\R^2}|u|^4\dd x\dd t
	=\int_{-1}^0\int_{\R^2}|U|^4\dd y\dd\tau<\infty,\qquad
	\int_{-\infty}^{-1}\int_{\R^2}|u|^4\dd x\dd t
	=\int_0^1\int_{\R^2}|U|^4\dd y\dd\tau<\infty .
\end{equation*}
Thus $u\in L^4(\R\times\R^2)$. By the Duhamel formula and
the dual of \eqref{eq:strichartz}, for $t_1<t_2$,
\begin{equation*}
	\norm{S(-t_2)u(t_2)-S(-t_1)u(t_1)}{L^2(\R^2)}
	=\left\|\int_{t_1}^{t_2}S(-t)(|u|^2u)(t)\dd t
	\right\|_{L^2(\R^2)}
	\lesssim \norm u{L^4([t_1,t_2]\times\R^2)}^3 .
\end{equation*}
The right-hand side tends to zero as both times tend to
$+\infty$ or both tend to $-\infty$. Hence the limits
$u_\pm=\lim_{t\to\pm\infty}S(-t)u(t)$ exist in $L^2(\R^2)$.
By unitarity, they are the unique states satisfying
\begin{equation*}
	\lim_{t\to\pm\infty}
	\norm{u(t)-S(t)u_\pm}{L^2(\R^2)}=0,
\end{equation*}
which is exactly the desired scattering.

\end{proof}

\vspace{10pt}
\section{The Hamiltonian system}\label{sec:vector}

To apply the scalar proof to \eqref{sys-eq:system}, we need to
check three points. First, the local mass and momentum identities
must give a nonnegative term in \eqref{principal} after summation
over the components. Second, the estimates for the modified
energy and the two spacetime norms in \eqref{bootstrap} must
hold for products of different components. Third, since $W$ may vanish at nonzero points, the bound for
$\norm{P_Nu(t)}{L^4(\R^2)}$ used to prove \eqref{localglobal} must follow
from mass conservation and the kinetic term
$\tfrac12\norm{\nabla P_Nu(t)}{L^2(\R^2)}^2$ in the modified
energy via \eqref{eq:sobolev-gn}. Once these
points are established, the scalar continuation and scattering
arguments apply.

We use the Euclidean norm on $\C^m$. All implicit constants may
depend on the fixed polynomial $W$ and the number of components.
Write $v=P_Nu$ and, as in \eqref{Y}, set
\begin{equation*}
	Y(h):=\sup_{d\in\R^2}
	\norm{|(P_Lu)(t,x+d)|\,|h(t,x)|}{L^2(J\times\R^2)} .
\end{equation*}

\begin{lemma}\label{lem:system-virial}
	Assume \eqref{sys-eq:nonnegative}. Let $u=(u_1,\ldots,u_m)$
	be a smooth solution of \eqref{sys-eq:system} on $J\ni0$, and set
	\begin{equation*}
		E(u(t)):=\frac12\norm{\nabla u(t)}{L^2(\R^2)}^2
		+\frac12\int_{\R^2}W(u(t,x))\dd x .
	\end{equation*}
	Then $\norm{u(t)}{L^2(\R^2)}=\norm{u(0)}{L^2(\R^2)}$ and
	$E(u(t))=E(u(0))$ for every $t\in J$.
	
	Let $v=P_Nu$, so that $v_j=P_Nu_j$. For scalar functions
	$f,g$, write the vector in \eqref{principal} as
	\begin{equation*}
		B(f,g)(t,x,y)
		:=\overline{g(t,y)}\,\nabla f(t,x)
		+f(t,x)\nabla\overline{g(t,y)} .
	\end{equation*}
	Let $I_1$ denote the term in \eqref{principal}. 	If $U=V=v$ in \eqref{eq:radial-action}, then
	\begin{equation}\label{eq:system-virial-positive}
		I_1(t)\gtrsim
		\sum_{j,k=1}^m\iint L(1+L|x-y|)^{-3}
		|B(v_j,v_k)(t,x,y)|^2\dd x\dd y .
	\end{equation}
	If $U=v$ and $V=h$, where
	$h\colon J\times\R^2\to\C$ is smooth and has sufficient
	spatial decay, then
	\begin{equation*}
		I_1(t)\gtrsim
		\sum_{j=1}^m\iint L(1+L|x-y|)^{-3}
		|B(v_j,h)(t,x,y)|^2\dd x\dd y .
	\end{equation*}
	Both estimates are uniform under spatial
	translations of $U$ and $V$.
\end{lemma}

\begin{proof}
	Since $F_j=\partial_{\bar z_j}W$ and $W$ is real and of
	bidegree $(2,2)$, Euler's identity and the chain rule give
	\begin{equation}\label{euler}
		\sum_{j=1}^m\bar z_jF_j(z)=2W(z),\qquad
		\partial W(v)
		=2\re\sum_{j=1}^m\overline{F_j(v)}\,\partial v_j,
		\qquad \partial\in\{\partial_t,\partial_1,\partial_2\}.
	\end{equation}
	Conservation of mass and energy follows from \eqref{euler}
	and integration by parts in \eqref{sys-eq:system}.
	
	To compute $I_1$, write the equation for $v_j=P_Nu_j$ as
	\begin{equation*}
		(\ii\partial_t+\Delta)v_j=F_j(v)+R_j,\qquad
		R_j:=P_NF_j(u)-F_j(v).
	\end{equation*}
	Sum the mass and momentum identities \eqref{masslaw} and
	\eqref{momlaw} over $j$, with
	\begin{equation*}
		\rho_v:=\sum_{j=1}^m|v_j|^2,\qquad
		p_{v,q}:=\im\sum_{j=1}^m\bar v_j\partial_qv_j .
	\end{equation*}
	By \eqref{euler}, the contribution of $F(v)$ to
	$\partial_t\rho_v$ is $2\im\sum_j\bar v_jF_j(v)=0$.
	Its contribution to $\partial_t p_{v,q}$ is
	\begin{equation*}
		\re\sum_{j=1}^m
		\bigl(\overline{F_j(v)}\,\partial_qv_j
		-\bar v_j\partial_qF_j(v)\bigr)
		=2\re\sum_{j=1}^m\overline{F_j(v)}\,\partial_qv_j
		-\partial_q(2W(v))
		=-\partial_qW(v).
	\end{equation*}
	
	For $U=V=v$, the terms containing $R_j$ enter
	\eqref{actualsource} and \eqref{testsource}. Substituting the component sums for $\rho_v$ and $p_v$ into
	\eqref{eq:radial-action} gives a sum over pairs $(v_j,v_k)$.
	For each pair, the Laplacian terms in the time derivative give
	the Hessian term in \eqref{principal} with $B=B(v_j,v_k)$. Integrating
	$-\nabla W(v)$ by parts in $x$ and $y$ gives
	\begin{equation*}
		I_1(t)
		=2\sum_{j,k=1}^m\iint
		\re\bigl(\overline{B(v_j,v_k)}\cdot
		D^2A(x-y)B(v_j,v_k)\bigr)\dd x\dd y
	\end{equation*}
	\begin{equation*}
		+\iint\Delta A(x-y)
		\bigl[W(v(t,x))\rho_v(t,y)
		+\rho_v(t,x)W(v(t,y))\bigr]\dd x\dd y .
	\end{equation*}
	For $U=v$ and $V=h$, the terms containing
	$(\ii\partial_t+\Delta)h$ enter \eqref{testsource}, so
	\begin{equation*}
		I_1(t)
		=2\sum_{j=1}^m\iint
		\re\bigl(\overline{B(v_j,h)}\cdot
		D^2A(x-y)B(v_j,h)\bigr)\dd x\dd y
		+\iint\Delta A(x-y)W(v(t,x))|h(t,y)|^2\dd x\dd y .
	\end{equation*}
	In both identities, $\Delta A\ge0$ by \eqref{hessian} and
	$W\ge0$ by \eqref{sys-eq:nonnegative}. Applying
	\eqref{hessian} to the real and imaginary parts of each $B$
	gives the asserted bounds. 
	
\end{proof}

To close \eqref{bootstrap} for the system, we use
Lemma~\ref{lem:system-virial} with $U=V=P_Nu$ and with
$U=P_Nu$, $V=h$, where $h$ solves the equation in
Lemma~\ref{forced-test}. We illustrate how the scalar estimates
extend to the system by bounding the variation of $\mathcal E_N(u)$
and the $L^2(J\times\R^2)$ norm of
$\nabla\bigl((P_Nu_j)\overline{P_Nu_k}\bigr)$.
\begin{lemma}\label{lem:system-bootstrap}
	Let $u$ be a smooth solution of \eqref{sys-eq:system} on a
	finite interval $J$, let $L=N^{1-\varepsilon}$, and set
	\begin{equation*}
		\mathcal E_N(u)
		=E(P_Nu)+\frac{N^2}{2}\norm{Q_Nu}{L^2(\R^2)}^2 .
	\end{equation*}
	Under \eqref{bootstrap}, for all sufficiently large $N$,
	\begin{equation*}
		\int_J\left|\frac{\dd}{\dd t}\mathcal E_N(u(t))\right|\dd t
		\lesssim_M N^{-\varepsilon/2},\qquad
		\norm{Q_Nu}{L^4(J\times\R^2)}\lesssim_M N^{-1},
	\end{equation*}
	\begin{equation*}
		Y(Q_{aN}u)\lesssim_M N^{-1-\varepsilon/2},\qquad
		\norm{Q_Lu}{L^4(J\times\R^2)}
		\lesssim_M N^{-3/4+2\varepsilon}.
	\end{equation*}
	The bounds are uniform in $J$.
\end{lemma}

\begin{proof}
	We show how to estimate the variation of $\mathcal E_N$ and
	the product of two different components. Write $v=P_Nu$ and
	$R_j=P_NF_j(u)-F_j(v)$. Differentiating the modified energy gives
	\begin{equation*}
		\frac{\dd}{\dd t}\mathcal E_N(u)
		=\im\sum_{j=1}^m\int_{\R^2}
		\bigl[\overline{-\Delta v_j}\,R_j
		+\overline{F_j(v)}\,R_j
		+N^2\overline{Q_Nu_j}\,Q_NF_j(u)\bigr]\dd x .
	\end{equation*}
	Each $F_j(u)$ is a finite sum of monomials
	$c^j_{\alpha\beta\gamma}u_\alpha\bar u_\beta u_\gamma$.
	For each monomial in $R_j$, the multiplier vanishes when
	all three frequencies have magnitude at most $N/3$. Thus one
	frequency is greater than $N/3$, and spatial integration
	requires a second frequency of magnitude $\gtrsim N$.
	In the last integral, $Q_Nu_j$ already supplies the first
	high frequency. The proof of Lemma~\ref{source} therefore
	applies to each choice of component indices, since
	\begin{equation*}
		\sup_{d\in\R^2}
		\norm{(P_Lu_\alpha)(t,x+d)(Q_{aN}u_\beta)(t,x)}
		{L^2(J\times\R^2)}
		\le Y(Q_{aN}u),\qquad
		\norm{Q_Lu_\gamma}{L^4(J\times\R^2)}
		\le\norm{Q_Lu}{L^4(J\times\R^2)} .
	\end{equation*}
	Summing over the finitely many monomials yields
	$\int_J|\frac{\dd}{\dd t}\mathcal E_N(u(t))|\dd t
	\lesssim_M N^{-\varepsilon/2}$ and gives \eqref{qs1}--\eqref{qs2}
	for $\bar v_jR_j$ with either derivative placement.
	
	For the product estimate, fix $j,k$. Since every summand on
	the right of \eqref{eq:system-virial-positive} is nonnegative,
	we may retain the term containing $B(v_j,v_k)$. Integrating $\frac{d\mathcal M}{dt}=I_1+I_2+I_3$ over $J$,
	as in the proof of \eqref{eq:radial-two-projected}, and using
	\eqref{eq:system-virial-positive} together with the preceding
	bounds for the terms containing $R_j$, we obtain
	\begin{equation*}
		\int_J\iint L(1+L|x-y|)^{-3}
		|B(v_j,v_k)(t,x,y)|^2\dd x\dd y\dd t\lesssim_M1 .
	\end{equation*}
	For fixed $t,y$, the $x$-Fourier support of $B(v_j,v_k)$ lies
	in $\{|\xi|\le2N\}$. On the diagonal,
	$B(v_j,v_k)(t,x,x)=\nabla_x(v_j\bar v_k)(t,x)$.
	Thus \eqref{diag} with $K=2N$ gives
	\begin{equation*}
		\norm{\nabla_x(v_j\bar v_k)}{L^2(J\times\R^2)}^2
		\lesssim \frac{N^2}{L}
		\int_J\iint L(1+L|x-y|)^{-3}
		|B(v_j,v_k)(t,x,y)|^2\dd x\dd y\dd t
		\lesssim_M \frac{N^2}{L}.
	\end{equation*}
	The same estimate holds under translations of
	$v_j$ and $v_k$.
	
	Applying the scalar proofs of \eqref{Ximprove}, \eqref{w},
	and \eqref{tail-improved} to each cubic monomial, using
	Lemma~\ref{lem:system-virial} also for the pairs $(v_j,h)$
	with scalar $h$, and summing over the component indices
	gives the remaining three estimates.
	
\end{proof}

\begin{proof}[Proof of Theorem~\ref{sys-thm:main}]
	For simplicity, we prove the result for $s=0$, since the case $s>0$ follows from the global $L^4(\R\times\R^2)$ bound by applying the Strichartz and fractional Leibniz estimates on intervals where this norm is small. 
	
	We first obtain the local $L^4$ bound used in
	Section~\ref{sec:global}. Since $W\ge0$, Lemma~\ref{lem:system-virial}
	gives global smooth solutions by the local theory and
	persistence of regularity. For the fixed $N$,
	the two quadratic terms scale exactly as in the scalar case.
	The quartic term satisfies, by \eqref{eq:sobolev-gn} applied
	to $|P_Nf_\lambda|$,
	\begin{equation*}
		0\le\int_{\R^2}W(P_Nf_\lambda)\dd x
		\lesssim \norm {P_Nf_{\lambda}}{L^2(\R^2)}^2
		\norm{\nabla P_Nf_\lambda}{L^2(\R^2)}^2
		\longrightarrow0
		\qquad(\lambda\to\infty).
	\end{equation*}
	Thus the same $N,\lambda$ give initial modified energy below
	$1/4$ on a neighborhood of $f$ in $L^2(\R^2;\C^m)$.
	Lemma~\ref{lem:system-bootstrap} gives the three strict
	improvements in \eqref{bootstrap}, so the bootstrap
	argument of Section~\ref{sec:global} applies on each finite
	interval in either time direction. For the rescaled smooth
	solutions it gives
	\begin{equation*}
		\sup_{t\in J}\norm{P_Nu(t)}{L^4(\R^2)}^4
		\lesssim M^2\sup_{t\in J}\norm{\nabla P_Nu(t)}{L^2(\R^2)}^2
		\lesssim_M1,\qquad
		\norm u{L^4(J\times\R^2)}
		\lesssim_M |J|^{1/4}+N^{-1}.
	\end{equation*}
    The second inequality follows from the first and the bound for
    $\norm{Q_Nu}{L^4(J\times\R^2)}$ in Lemma~\ref{lem:system-bootstrap}.

	For scattering, $F(cz)=|c|^2cF(z)$ for every $c\in\C$.
	Hence the transform $U$ in \eqref{pseudo} still satisfies
	\begin{equation*}
		(\ii\partial_\tau+\Delta_y)U(\tau,y)
		=|\tau|^{-3}e^{\ii|y|^2/(4\tau)}
		F\bigl(u(-1/\tau,y/\tau)\bigr)=F(U(\tau,y)).
	\end{equation*}
	Solving forward from $U(-1)$ and backward from $U(1)$ to
	$\tau=0$ gives, exactly as in Section~\ref{sec:global},
	\begin{equation*}
		\int_{|t|\ge1}\int_{\R^2}|u(t,x)|^4\dd x\dd t
		=\int_{0<|\tau|\le1}\int_{\R^2}|U(\tau,y)|^4\dd y\dd\tau
		<\infty .
	\end{equation*}
    Thus it recovers the desired global $L^4$ bound, which directly gives scattering.

\end{proof}

\begin{proof}[Proof of Corollary~\ref{sys-cor:matrix}]
	Set $W(z)=\tfrac12\sum_{j,k=1}^m\beta_{jk}|z_j|^2|z_k|^2$.
	By symmetry,
	$\partial_{\bar z_j}W(z)=z_j\sum_{k=1}^m\beta_{jk}|z_k|^2$.
	The copositivity condition \eqref{sys-eq:copositive}
	gives $W(z)\ge0$, so Theorem~\ref{sys-thm:main} applies.
	
\end{proof}

\vspace{10pt}
\section{Invariant complex lines and necessity of the sign condition}
\label{sec:sharpness}

Theorem~\ref{sys-thm:main} gives global existence for every
$L^2(\R^2;\C^m)$ datum when $W\ge0$. We show that if $W$ takes
a negative value, then \eqref{sys-eq:system} admits smooth compactly
supported data that blow up in finite time. We do this by restricting
the system to a complex line on which it reduces to the focusing
scalar equation. We shall use the same reduction in
Section~\ref{sec:negative-systems} to transfer the scalar examples
of Theorem~\ref{neg-thm:scalar} to every nonzero $W$.

\begin{lemma}\label{neg-lem:direction}
	Let $W:\C^m\to\R$ be a polynomial of bidegree $(2,2)$.
	At every critical point $a$ of the restriction of $W$ to
	the unit sphere,
	\begin{equation}\label{neg-eq:direction}
		\partial_{\bar z}W(a)=\lambda a,\qquad \lambda=2W(a).
	\end{equation}
	Consequently, for every $\zeta\in\C$,
	\begin{equation}\label{neg-eq:homogeneity}
		\partial_{\bar z}W(\zeta a)=\lambda|\zeta|^2\zeta a .
	\end{equation}
	If $W\not\equiv0$, one can choose $a$ with $\lambda\ne0$.
	If $W$ takes a negative value, one can choose $\lambda<0$.
\end{lemma}

\begin{proof}
	Identity \eqref{neg-eq:direction} is the quartic case of
	\cite[Section~3]{Masaki}. Indeed, write $z_j=x_j+\ii y_j$. Since $a$ is a critical point
	of $W$ under the constraint $\sum_{j=1}^m(x_j^2+y_j^2)=1$,
	the real Lagrange multiplier rule gives $\lambda\in\R$ such that,
	for every $j$,
	\begin{equation*}
		\partial_{x_j}W(a)=2\lambda\re a_j,\qquad
		\partial_{y_j}W(a)=2\lambda\im a_j .
	\end{equation*}
	Hence
	\begin{equation*}
		\partial_{\bar z_j}W(a)
		=\tfrac12\bigl(\partial_{x_j}W(a)+\ii\partial_{y_j}W(a)\bigr)
		=\lambda a_j .
	\end{equation*}
	
	To determine $\lambda$, recall that every monomial of $W$ has
	the form $z_pz_q\bar z_k\bar z_l$. Differentiating the two
	barred factors gives
	\begin{equation*}
		\sum_{j=1}^m\bar z_j\,\partial_{\bar z_j}
		\bigl(z_pz_q\bar z_k\bar z_l\bigr)
		=z_pz_q\bigl(\bar z_k\bar z_l+\bar z_l\bar z_k\bigr)
		=2z_pz_q\bar z_k\bar z_l .
	\end{equation*}
	Summing over the monomials and evaluating at $a$, we obtain
	\begin{equation*}
		2W(a)=\sum_{j=1}^m\bar a_j\,\partial_{\bar z_j}W(a)
		=\lambda\sum_{j=1}^m|a_j|^2=\lambda .
	\end{equation*}
	Each monomial of $\partial_{\bar z_j}W$ contains two unbarred
	factors and one barred factor. Thus, for every $\zeta\in\C$,
	\begin{equation*}
		\partial_{\bar z_j}W(\zeta a)
		=\zeta^2\bar\zeta\,\partial_{\bar z_j}W(a)
		=\lambda|\zeta|^2\zeta a_j ,
	\end{equation*}
	which proves \eqref{neg-eq:homogeneity}.
	
	Finally, if $W\not\equiv0$, either the maximum of $W$ on
	the unit sphere is positive or its minimum is negative. At the
	corresponding extremum $a$, we have $\lambda=2W(a)\ne0$.
	In particular, if $W(z)<0$, choosing a minimizer gives
	$\lambda=2W(a)\le2|z|^{-4}W(z)<0$.
	
\end{proof}

For $a$ as in Lemma~\ref{neg-lem:direction} with $\lambda\ne0$,
\eqref{neg-eq:homogeneity} gives
\begin{equation*}
	(\ii\partial_t+\Delta)\left(\frac{a}{\sqrt{|\lambda|}}v\right)
	-\partial_{\bar z}W\left(\frac{a}{\sqrt{|\lambda|}}v\right)
	=\frac{a}{\sqrt{|\lambda|}}
	\bigl((\ii\partial_t+\Delta)v
	-\operatorname{sgn}(\lambda)|v|^2v\bigr).
\end{equation*}
Thus a negative minimum of $W$ gives a copy of the focusing cubic
equation inside \eqref{sys-eq:system}.

\begin{proposition}\label{sys-prop:negative}
	If $W$ takes a negative value, there exists
	$u_0\in C_c^\infty(\R^2;\C^m)$ for which the solution of
	\eqref{sys-eq:system} blows up in finite positive time.
\end{proposition}

\begin{proof}
	By the preceding identity, it suffices to construct a solution that blows up in finite time for the scalar focusing equation
	\begin{equation*}
		\left\{
		\begin{array}{l}
			(\ii\partial_t+\Delta)v=-|v|^2v,\\
			v(0)=b\phi .
		\end{array}
		\right.
	\end{equation*}
	Fix a nonzero $\phi\in C_c^\infty(\R^2)$ and choose $b>0$
	large enough that
	\begin{equation*}
		E(v(0))
		=\frac{b^2}{2}\norm{\nabla\phi}{L^2(\R^2)}^2
		-\frac{b^4}{4}\norm{\phi}{L^4(\R^2)}^4<0 .
	\end{equation*}
	Since $|x|\phi\in L^2(\R^2)$, Glassey's virial identity
	\cite{Glassey} gives
	\begin{equation*}
		\frac{\dd^2}{\dd t^2}\int_{\R^2}|x|^2|v(t,x)|^2\dd x
		=8\norm{\nabla v(t)}{L^2(\R^2)}^2
		-4\norm{v(t)}{L^4(\R^2)}^4
		=16E(v(0))<0,
	\end{equation*}
	which implies the maximal existence time $T_+$ is finite.
	
\end{proof}

\vspace{10pt}
\section{An approximate profile}\label{sec:negative-profile}

Write $F(z):=|z|^2z$. To prove
Theorem~\ref{neg-thm:scalar}, we first construct one approximate
solution whose $H^s$ norm is small before a time close to $1$
and whose value at $t=1$ has a large pairing with a function
bounded in $H^{-s}$. In Section~\ref{sec:negative-correction}
we shall add this approximate solution to a given smooth exact
solution.

\begin{proposition}\label{neg-prop:profile}
	Fix $\sigma>0$ and $\mu\in\{-1,1\}$. There exist
	$t_n\uparrow1$, $Z_n\in C([0,1];L^2(\R^2))$, and
	$g_n\in\Sclass(\R^2)$ such that
	\begin{equation*}
		Z_n|_{[0,t_n]}\in C^\infty([0,t_n];\Sclass(\R^2)),
		\qquad
		Z_n|_{[t_n,1]}\in C^\infty([t_n,1];\Sclass(\R^2)).
	\end{equation*}
	As $n\to\infty$,
	\begin{equation}\label{neg-eq:profile-small}
		\sup_{0\le t<t_n}\norm{Z_n(t)}{H^{-\sigma}(\R^2)}
		+\norm{Z_n}{L^2([0,1];L^\infty(\R^2))}
		\longrightarrow0 ,
	\end{equation}
	\begin{equation}\label{neg-eq:profile-error}
		\norm{(\ii\partial_t+\Delta)Z_n-\mu F(Z_n)}
		{L^1([0,1];L^2(\R^2))}
		\longrightarrow0 .
	\end{equation}
	Moreover,
	\begin{equation*}
		\norm{Z_n(0)}{L^2(\R^2)}\longrightarrow\infty,\qquad
		|\langle Z_n(1),g_n\rangle|\longrightarrow\infty,\qquad
		\norm{g_n}{L^2(\R^2)}\le1,\qquad
		\sup_n\norm{g_n}{H^\sigma(\R^2)}<\infty .
	\end{equation*}
\end{proposition}

\vspace{7pt}
\subsection{Construction of the approximate solution}

We prescribe three packets at $t=t_n$. On $[0,t_n]$ we solve
the cubic equation backwards for each packet and add the three
solutions. On $[t_n,1]$ we take a finite Duhamel iterate with
their sum as its initial value.

For each sufficiently large integer $n$, set
\begin{equation*}
	\ell:=n^{-4/\sigma^2},\qquad N:=n^{20+8/\sigma},\qquad
	t_n:=1-\frac{n\ell^2}{N}.
\end{equation*}
Write
\begin{equation*}
	Z_n(t,x):=\frac{\sqrt N}{n^3\ell}\,U(\tau,y),
	\qquad \tau:=\frac{N(t-1)}{\ell^2},\qquad y:=\frac{x}{\ell},
\end{equation*}
where $U$ will be chosen later.
The two intervals in the rescaled time variable are
\begin{equation*}
	0\le t\le t_n
	\quad\Longleftrightarrow\quad -N/\ell^2\le\tau\le-n,
	\qquad
	t_n\le t\le1
	\quad\Longleftrightarrow\quad -n\le\tau\le0 .
\end{equation*}

We first write the estimates required by
Proposition~\ref{neg-prop:profile} in these variables.
Since $\dd x=\ell^2\dd y$ and $\dd t=(\ell^2/N)\dd\tau$,
\begin{equation*}
	\norm{Z_n(t)}{H^{-\sigma}(\R^2)}
	=n^{7+4/\sigma}
	\left\|\langle\xi/\ell\rangle^{-\sigma}
	\widehat U(\tau,\xi)\right\|_{L^2(\R^2)},\qquad
	\norm{Z_n(t)}{L^2(\R^2)}
	=n^{7+4/\sigma}\norm{U(\tau)}{L^2(\R^2)} .
\end{equation*}
The square of the $L^2([0,1];L^\infty(\R^2))$ norm acquires
the factor $n^{-6}$. For the error in the equation, the same
change of variables gives
\begin{equation*}
	\norm{(\ii\partial_t+\Delta)Z_n-\mu F(Z_n)}
	{L^1([0,1];L^2(\R^2))}
	=n^{7+4/\sigma}
	\norm{(\ii\partial_\tau+N^{-1}\Delta_y)U-\mu n^{-6}F(U)}
	{L^1([-N/\ell^2,0];L^2(\R^2))}.
\end{equation*}
Thus, for \eqref{neg-eq:profile-small} and
\eqref{neg-eq:profile-error}, it suffices to construct $U$
continuous at $\tau=-n$ and prove
\begin{equation}\label{neg-eq:rescaled-small}
	\sup_{-N/\ell^2\le\tau<-n}
	\left\|\langle\xi/\ell\rangle^{-\sigma}
	\widehat U(\tau,\xi)\right\|_{L^2(\R^2)}
	=o(n^{-7-4/\sigma}),\qquad
	n^{-6}\int_{-N/\ell^2}^0\norm{U(\tau)}{L^\infty(\R^2)}^2
	\dd\tau=o(1),
\end{equation}
\begin{equation}\label{neg-eq:rescaled-error}
	\int_{-N/\ell^2}^0
	\norm{(\ii\partial_\tau+N^{-1}\Delta_y)U-\mu n^{-6}F(U)}
	{L^2(\R^2)}\dd\tau=o(n^{-7-4/\sigma}).
\end{equation}
To ensure that $\|Z_n(0)\|_{L^2}$ tends to infinity, it suffices that
$\norm{U(-N/\ell^2)}{L^2(\R^2)}\ge1$.
We give the corresponding condition for the terminal pairing
once the test function $g_n$ has been defined.

Let $\psi(x):=(2/\pi)^{1/4}e^{-x^2}$ and prescribe the
one-dimensional data
\begin{equation*}
	a(-n,x):=\psi(x),\qquad
	b(-n,x):=e^{\ii Nx-\ii N(x+2n)^2/n^4}\psi(x+2n).
\end{equation*}
For $-n\le\tau\le0$, define their free evolutions by
\begin{equation*}
	a(\tau):=e^{\ii(\tau+n)\partial_x^2/N}a(-n),\qquad
	b(\tau):=e^{\ii(\tau+n)\partial_x^2/N}b(-n).
\end{equation*}
For one-dimensional functions $f,g$, write
$(f\otimes g)(y_1,y_2):=f(y_1)g(y_2)$, with both factors
evaluated at the same time when they depend on $\tau$.
Since $\Delta_y=\partial_{y_1}^2+\partial_{y_2}^2$,
the products $b\otimes a$, $b\otimes b$, and $a\otimes b$
are two-dimensional free solutions.

Let $e_1=(1,0)$ and $e_2=(0,1)$. Since the factor
$e^{\ii Nx}$ in $b(-n,x)$ occurs in the first coordinate of
$b\otimes a$, in both coordinates of $b\otimes b$, and in the
second coordinate of $a\otimes b$, the three data contain the
respective factors $e^{\ii Ne_1\cdot y}$,
$e^{\ii N(e_1+e_2)\cdot y}$, and $e^{\ii Ne_2\cdot y}$.

For $\theta\in[0,2\pi]$, set
\begin{equation}\label{neg-eq:three-free-packets}
	P_0(\tau):=e^{\ii\theta}(b\otimes a)(\tau)
	+(b\otimes b)(\tau)+e^{3\ii\theta}(a\otimes b)(\tau).
\end{equation}

On $-n\le\tau\le0$, fix an integer $k\ge1$ with
$5(k+1)>7+4/\sigma$. Starting from
\eqref{neg-eq:three-free-packets}, define
\begin{equation}\label{neg-eq:collision-approx}
	P_{j+1}(\tau):=P_0(\tau)-\ii\mu n^{-6}
	\int_{-n}^\tau e^{\ii(\tau-r)\Delta_y/N}F(P_j(r))\dd r,
	\qquad 0\le j<k,
\end{equation}
and set $U(\tau):=P_k(\tau)$.

On $-N/\ell^2\le\tau\le-n$, let
$u_{e_1},u_{e_1+e_2},u_{e_2}$ solve
\begin{equation*}
	\left\{
	\begin{array}{l}
		(\ii\partial_\tau+N^{-1}\Delta_y)u_p
		=\mu n^{-6}F(u_p),
		\qquad p=(p_1,p_2)\in\{e_1,e_1+e_2,e_2\},\\
		u_{e_1}(-n)=(b\otimes a)(-n),\qquad
		u_{e_1+e_2}(-n)=(b\otimes b)(-n),\qquad
		u_{e_2}(-n)=(a\otimes b)(-n).
	\end{array}
	\right.
\end{equation*}
We also set
\begin{equation*}
	U(\tau):=e^{\ii\theta}u_{e_1}(\tau)
	+u_{e_1+e_2}(\tau)+e^{3\ii\theta}u_{e_2}(\tau).
\end{equation*}
The two definitions agree at $\tau=-n$ since $P_k(-n)=P_0(-n)$.

Finally, set
\begin{equation*}
	g_n(x):=\ell^{\sigma-1}(a(0)\otimes a(0))(x/\ell).
\end{equation*}
Choose $\theta$ to maximize
$|\langle P_k(0),a(0)\otimes a(0)\rangle|$ on $[0,2\pi]$.
The pairing is a trigonometric polynomial in $\theta$ by
\eqref{neg-eq:three-free-packets} and \eqref{neg-eq:collision-approx},
so the maximum is attained. By scaling and unitarity,
\begin{equation}\label{neg-eq:test}
	\norm{g_n}{L^2(\R^2)}=\ell^\sigma\le1,\qquad
	\norm{g_n}{H^\sigma(\R^2)}
	\le\norm{a(0)\otimes a(0)}{H^\sigma(\R^2)}
	=\norm{\psi\otimes\psi}{H^\sigma(\R^2)} .
\end{equation}
Since $(\sqrt N/n^3)\ell^\sigma=n^7$, the terminal pairing is
\begin{equation*}
	|\langle Z_n(1),g_n\rangle|
	=n^7|\langle P_k(0),a(0)\otimes a(0)\rangle|.
\end{equation*}
It therefore suffices to prove
\begin{equation}\label{neg-eq:rescaled-pairing}
	|\langle P_k(0),a(0)\otimes a(0)\rangle|\gtrsim n^{-6}.
\end{equation}

\vspace{7pt}
\subsection{Proof of Proposition~\ref{neg-prop:profile}}

We first prove \eqref{neg-eq:rescaled-small} and
\eqref{neg-eq:rescaled-error} on $[-N/\ell^2,-n]$. On this
interval,
\begin{equation*}
	U=e^{\ii\theta}u_{e_1}+u_{e_1+e_2}+e^{3\ii\theta}u_{e_2},
\end{equation*}
and substitution of the three equations gives
\begin{equation}\label{neg-eq:early-U-error}
	(\ii\partial_\tau+N^{-1}\Delta_y)U-\mu n^{-6}F(U)
	=\mu n^{-6}\bigl[
	e^{\ii\theta}F(u_{e_1})+F(u_{e_1+e_2})
	+e^{3\ii\theta}F(u_{e_2})-F(U)\bigr].
\end{equation}
Note that every product remaining on the right (after expansion and cancellation) contains two factors
$u_p,u_q$ with $p\ne q$, up to conjugation.

To estimate these mixed products, we first bound each $u_p$
and its spatial moments about $y=2p\tau$, so that the separation
between $2p\tau$ and $2q\tau$ yields a small bound for their product.
\begin{lemma}\label{neg-lem:single-packet-bounds}
	Fix an integer $K>8+4/\sigma$. For all sufficiently large $n$,
	the three solutions $u_p$, $p\in\{e_1,e_1+e_2,e_2\}$, exist
	on $[-N/\ell^2,-n]$ and are smooth with values in
	$\Sclass(\R^2)$. On this interval, with
	$d(\tau):=1+4|\tau+n|/n^4$, we have
	\begin{equation}\label{neg-eq:packet-decay}
		d(\tau)^{1/2}\norm{u_p(\tau)}{L^\infty(\R^2)}
		+\left\|
		\left\langle\frac{y-2p\tau}{d(\tau)}\right\rangle^K
		u_p(\tau)\right\|_{L^2(\R^2)}
		\lesssim_K1
	\end{equation}
	and
	\begin{equation}\label{neg-eq:packet-fourier-tail}
		\left\|\mathbf1_{\{|\xi-Np|>N/4\}}
		\widehat{u_p}(\tau,\xi)\right\|_{L^2(\R^2)}
		\lesssim_K n^{-4K}.
	\end{equation}
\end{lemma}

\begin{proof}
	To prove \eqref{neg-eq:packet-decay} and
	\eqref{neg-eq:packet-fourier-tail}, we should first study the multiplication
	by $y_j-2p_j\tau$  and the operator
	$-\ii\partial_{y_j}-Np_j$, whose Fourier multiplier is just
	$\xi_j-Np_j$. For convenience, we call them \emph{centered operators}.
	
	For $-N/\ell^2\le\tau\le-n$, the Duhamel formula is
	\begin{equation}\label{neg-eq:single-packet-duhamel}
		u_p(\tau)=e^{\ii(\tau+n)\Delta_y/N}u_p(-n)
		-\ii\mu n^{-6}\int_{-n}^{\tau}
		e^{\ii(\tau-r)\Delta_y/N}F(u_p(r))\dd r .
	\end{equation}
	We begin with the multiplication
	by $y_j-2p_j\tau$.
	To see what it requires of the datum, apply $y_j-2p_j\tau$. Indeed, we can obtain
	\begin{equation}\label{qw}
		(y_j-2p_j\tau)e^{\ii(\tau+n)\Delta_y/N}u_p(-n)
		=e^{\ii(\tau+n)\Delta_y/N}
		\left[(y_j+2np_j)u_p(-n)
		+\frac{2(\tau+n)}{N}(-\ii\partial_{y_j}-Np_j)u_p(-n)\right].
	\end{equation}

	We compute the two terms in brackets in \eqref{qw}.
	Set $z=y+2np$. The first term is $z_j u_p(-n,y)$.
	For the second, the datum is
	\begin{equation*}
		u_p(-n,y)
		=e^{\ii Np\cdot y-\ii(N/n^4)\sum_{k=1}^2p_kz_k^2}
		\psi(z_1)\psi(z_2).
	\end{equation*}
	Since $\psi'(z_j)=-2z_j\psi(z_j)$, differentiation gives
	\begin{equation*}
		(-\ii\partial_{y_j}-Np_j)u_p(-n,y)
		=-\frac{2Np_j}{n^4}z_j u_p(-n,y)+2\ii z_j u_p(-n,y).
	\end{equation*}
	 To keep the $L^2$ bound uniform in $n$ under repeated applications,
	 we move the term containing $N/n^4$ to the left so that each
	 application changes the polynomial factor by the $n$-independent
	 rule. More precisely, for any polynomial $q$ on $\R^2$, we can conclude
	\begin{equation}\label{fir}
		(y_j+2np_j)\bigl[q(z)u_p(-n,y)\bigr]
		=z_jq(z)u_p(-n,y),
	\end{equation}
	\begin{equation}\label{sec}
		\left(-\ii\partial_{y_j}-Np_j
		+\frac{2Np_j}{n^4}(y_j+2np_j)\right)
		\bigl[q(z)u_p(-n,y)\bigr]
		=\bigl(-\ii\partial_{z_j}q(z)+2\ii z_jq(z)\bigr)u_p(-n,y).
	\end{equation}
	Starting from $q=1$, each application of \eqref{fir} or
	\eqref{sec} therefore gives a new polynomial times $u_p(-n)$: $q(z)\mapsto z_jq(z)$ or $q(z)\mapsto -\ii\partial_{z_j}q(z)+2\ii z_jq(z)$.
	Thus, for any fixed sequence of applications, the resulting
	polynomial $Q$ has coefficients independent of $n$, and
	\begin{equation*}
		\norm{Q(y+2np)u_p(-n,y)}{L^2(\R^2)}
		=\norm{Q(z)\psi(z_1)\psi(z_2)}{L^2(\R^2)} .
	\end{equation*}
	
 To use these bounds at time $\tau$, we ask
which operator acting on the free solution
$e^{\ii(\tau+n)\Delta_y/N}u_p(-n)$ gives the free evolution
of the corresponding operator applied to $u_p(-n)$.

For multiplication by $y_j+2np_j$ in \eqref{fir}, we can find
\begin{equation*}
	e^{\ii(\tau+n)\Delta_y/N}(y_j+2np_j)
	=\underbrace{\left(y_j+2np_j+
		\frac{2\ii(\tau+n)}{N}\partial_{y_j}\right)}_{X_j(\tau)}
	e^{\ii(\tau+n)\Delta_y/N}.
\end{equation*}
For the operator in \eqref{sec}, we also have
\begin{equation*}
	e^{\ii(\tau+n)\Delta_y/N}
	\left(-\ii\partial_{y_j}-Np_j+
	\frac{2Np_j}{n^4}(y_j+2np_j)\right)
	=\underbrace{\left(-\ii\partial_{y_j}-Np_j+
		\frac{2Np_j}{n^4}X_j(\tau)\right)}_{Y_j(\tau)}
	e^{\ii(\tau+n)\Delta_y/N}.
\end{equation*}
	In particular, $X_j(-n)$ and $Y_j(-n)$ are exactly the operators in \eqref{fir} and \eqref{sec}, respectively. For $T=X_j$ or $Y_j$, applying these formulas at $\tau$
	and $r$ gives
	\begin{equation}\label{aqw}
		T(\tau)e^{\ii(\tau-r)\Delta_y/N}
		=e^{\ii(\tau+n)\Delta_y/N}T(-n)e^{-\ii(r+n)\Delta_y/N}
		=e^{\ii(\tau-r)\Delta_y/N}T(r).
	\end{equation}
	It also holds for each ordered composition of
	$X_1,X_2,Y_1,Y_2$.
	
	Note that we can actually express the \emph{centered operators}
	in terms of $X_j$ and $Y_j$. In fact, we have
	\begin{equation}\label{asdfgh}
		y_j-2p_j\tau
		=X_j(\tau)+\frac{2(\tau+n)}N(-\ii\partial_{y_j}-Np_j),
		\qquad
		-\ii\partial_{y_j}-Np_j
		=Y_j(\tau)-\frac{2Np_j}{n^4}X_j(\tau).
	\end{equation}
	Substituting the second identity into the first gives
	\begin{equation}\label{asdfg}
		\frac{y_j-2p_j\tau}{d(\tau)}
		=\frac{1-4p_j(\tau+n)/n^4}{d(\tau)}X_j(\tau)
		+\frac{2(\tau+n)}{Nd(\tau)}Y_j(\tau).
	\end{equation}
It therefore suffices to control the compositions
of $X_j$ and $Y_j$, which we collect in the following norm.
	
   For an integer $m\ge0$, define
   \begin{equation*}
   	\norm f{m,r}
   	=\norm f{L^2(\R^2)}
   	+\sum_{k=1}^{m}
   	\sum_{T_1,\ldots,T_k\in\{X_1(r),X_2(r),Y_1(r),Y_2(r)\}}
   	\norm{T_1\cdots T_kf}{L^2(\R^2)} .
   \end{equation*}
   The inner sum includes every ordered choice of the $k$ operators.
   We claim that
   \begin{equation}\label{neg-eq:moments}
   	\sup_{-N/\ell^2\le\tau\le-n}
   	\norm{u_p(\tau)}{K,\tau}\lesssim_K1 .
   \end{equation}
   Assuming (\ref{neg-eq:moments}), we first derive the conclusions of the lemma.
   
   For the $L^\infty$ estimate in \eqref{neg-eq:packet-decay},
   substitute the formula for $X_j(r)$ into $Y_j(r)$:
   \begin{equation*}
   	Y_j(r)
   	=\left(1-\frac{4p_j(r+n)}{n^4}\right)
   	(-\ii\partial_{y_j}-Np_j)
   	+\frac{2Np_j}{n^4}(y_j-2p_jr).
   \end{equation*}
   The multiplication coefficients are real, so
   \begin{equation*}
   	\ii\bigl(\bar fY_j(r)g-\overline{Y_j(r)f}\,g\bigr)
   	=\left(1-\frac{4p_j(r+n)}{n^4}\right)\partial_j(\bar f g).
   \end{equation*}
   Applying this identity in $y_2$ and then in $y_1$ gives
   \begin{equation*}
   	\prod_{j=1}^2\left(1-\frac{4p_j(r+n)}{n^4}\right)
   	\partial_1\partial_2|f|^2
   	=2\re\bigl(\overline{Y_1(r)f}\,Y_2(r)f
   	-\bar fY_1(r)Y_2(r)f\bigr).
   \end{equation*}
   Since $r\le-n$ and $p\in\{(1,0),(1,1),(0,1)\}$,
   \begin{equation*}
   	\prod_{j=1}^2\left(1-\frac{4p_j(r+n)}{n^4}\right)
   	=d(r)^{p_1+p_2}\ge d(r).
   \end{equation*}
   Integrating over $(-\infty,y_1]\times(-\infty,y_2]$ and
   applying Cauchy--Schwarz yields
   \begin{equation}\label{neg-eq:packet-sup}
   	d(r)\norm f{L^\infty(\R^2)}^2
   	\le 2\norm{Y_1(r)f}{L^2(\R^2)}\norm{Y_2(r)f}{L^2(\R^2)}
   	+2\norm f{L^2(\R^2)}\norm{Y_1(r)Y_2(r)f}{L^2(\R^2)} .
   \end{equation}
   Applying it to $u_p(\tau)$ and using
   \eqref{neg-eq:moments}, we obtain
   \begin{equation*}
   	\norm{u_p(\tau)}{L^\infty(\R^2)}
   	\lesssim_K d(\tau)^{-1/2},
   \end{equation*}
   which is exactly the $L^\infty$ bound in \eqref{neg-eq:packet-decay}.
  
   To prove the $L^2$ bound, we note that the coefficients in \eqref{asdfg} depend only on time and satisfy
   \begin{equation*}
   	\frac{|1-4p_j(\tau+n)/n^4|}{d(\tau)}\le1,
   	\qquad
   	\frac{2|\tau+n|}{Nd(\tau)}\le\frac{n^4}{2N}\le1.
   \end{equation*}
   Expanding the $K$-th power as a sum of ordered compositions
   and using \eqref{neg-eq:moments}, we obtain
   \begin{equation*}
   	\left\|\left\langle\frac{y-2p\tau}{d(\tau)}\right\rangle^K
   	u_p(\tau)\right\|_{L^2(\R^2)}\lesssim_K1 .
   \end{equation*}

   For \eqref{neg-eq:packet-fourier-tail}, the second identity in \eqref{asdfgh} gives
   \begin{equation*}
   	\frac{-\ii\partial_{y_j}-Np_j}{N}
   	=\frac{Y_j(\tau)}N-\frac{2p_j}{n^4}X_j(\tau).
   \end{equation*}
   Both coefficients are $O(n^{-4})$, so each ordered composition
   in the $K$th power has coefficient $O_K(n^{-4K})$.
   On $|\xi-Np|>N/4$, at least one coordinate satisfies
   $|\xi_j-Np_j|>N/(4\sqrt2)$. Plancherel and
   \eqref{neg-eq:moments} therefore give
   \begin{equation*}
   	\left\|\mathbf1_{\{|\xi-Np|>N/4\}}\widehat{u_p}(\tau,\xi)
   	\right\|_{L^2(\R^2)}
   	\lesssim_K\sum_{j=1}^2
   	\left\|\left(\frac{-\ii\partial_{y_j}-Np_j}{N}\right)^K
   	u_p(\tau)\right\|_{L^2(\R^2)}
   	\lesssim_K n^{-4K}.
   \end{equation*}
   
  We next check that $u_p$ is smooth in time with values in
  $\Sclass(\R^2)$. It suffices to show 
  \begin{equation*}
  	\sup_{-N/\ell^2\le\tau\le-n}
  	\norm{u_p(\tau)}{m,\tau}\lesssim_{m,n} 1 , \quad \forall m\ge K.
  \end{equation*}
  Indeed, we can express $y_j$ and $\partial_{y_j}$
  in terms of $X_j(\tau)$, $Y_j(\tau)$.
  
  Note that each choice of $T\in\{X_1(\tau),X_2(\tau),Y_1(\tau),Y_2(\tau)\}$ satisfies
  \begin{equation*}
  	T(f_1\bar f_2f_3)
  	=(Tf_1)\bar f_2f_3-f_1\overline{Tf_2}f_3
  	+f_1\bar f_2(Tf_3).
  \end{equation*}
   After at most $m$ applications, put
  the factor with the most operators in $L^2$. Each other factor
  carries at most $\lfloor m/2\rfloor$ operators. Since
  $\lfloor m/2\rfloor+2\le m-1$, \eqref{neg-eq:packet-sup} gives
  \begin{equation}\label{az}
  	\norm{F(u_p(r))}{m,r}
  	\lesssim_m d(r)^{-1}
  	\norm{u_p(r)}{m-1,r}^2\norm{u_p(r)}{m,r}.
  \end{equation}
  Commuting the compositions through the propagator by
  \eqref{aqw} and applying unitarity to
  \eqref{neg-eq:single-packet-duhamel}, we obtain
  \begin{equation*}
  	\norm{u_p(\tau)}{m,\tau}
  	\le \norm{u_p(-n)}{m,-n}
  	+C_mn^{-6}\int_\tau^{-n}d(r)^{-1}
  	\norm{u_p(r)}{m-1,r}^2\norm{u_p(r)}{m,r}\dd r .
  \end{equation*}
  Thus,  Gronwall's inequality and induction on $m$ with base case given by \eqref{neg-eq:moments} give spatial
  Schwartz regularity.
  And temporal regularity follows directly from the equation
  \begin{equation*}
  	\partial_\tau u_p
  	=\ii N^{-1}\Delta_yu_p-\ii\mu n^{-6}F(u_p).
  \end{equation*}

   We now complete the proof by proving \eqref{neg-eq:moments} through contraction for the
   solution map in \eqref{neg-eq:single-packet-duhamel}.
   By \eqref{aqw} and unitarity,
   \begin{equation*}
   	\norm{e^{\ii(\tau-r)\Delta_y/N}f}{K,\tau}
   	=\norm f{K,r}.
   \end{equation*}
   The free term therefore satisfies
   \begin{equation*}
   	\sup_{-N/\ell^2\le\tau\le-n}
   	\norm{e^{\ii(\tau+n)\Delta_y/N}u_p(-n)}{K,\tau}
   	=\norm{u_p(-n)}{K,-n}\lesssim_K1
   \end{equation*}
   by \eqref{fir} and \eqref{sec}.
   
   For the Duhamel integral, the argument in \eqref{az} actually gives the following trilinear estimate
   \begin{equation*}
   	\norm{f_1\bar f_2f_3}{K,r}
   	\lesssim_K d(r)^{-1}\prod_{j=1}^3\norm{f_j}{K,r}.
   \end{equation*}

   The choices of $N$ and $\ell$ give
   \begin{equation}\label{neg-eq:eta}
   	n^{-6}\int_{-N/\ell^2}^{-n}\frac{\dd r}{d(r)}
   	=\frac{1}{4n^2}
   	\log\left(1+\frac{4(N/\ell^2-n)}{n^4}\right)
   	\lesssim_\sigma n^{-2}\log n .
   \end{equation}
   Choose $B$, independent of $n$, such that
   $\norm{u_p(-n)}{K,-n}\le B/2$. On the ball
   $\lbrace f: \sup_r\norm{f(r)}{K,r}\le B\rbrace$, the trilinear estimate yields
   \begin{equation*}
   	\sup_\tau
   	\left\|n^{-6}\int_{-n}^{\tau}
   	e^{\ii(\tau-r)\Delta_y/N}F(f(r))\dd r\right\|_{K,\tau}
   	\le C_{\sigma,K}B^3n^{-2}\log n .
   \end{equation*}
   For two functions in this ball, expand
   \begin{equation*}
   	F(f)-F(g)
   	=(f-g)\bar f f+g\overline{(f-g)}f+g\bar g(f-g).
   \end{equation*}
   The same estimates give
   \begin{equation*}
   	\sup_\tau
   	\left\|n^{-6}\int_{-n}^{\tau}
   	e^{\ii(\tau-r)\Delta_y/N}
   	\bigl(F(f(r))-F(g(r))\bigr)\dd r\right\|_{K,\tau}
   	\le C_{\sigma,K}B^2n^{-2}\log n
   	\sup_r\norm{f(r)-g(r)}{K,r}.
   \end{equation*}
   Here $C_{\sigma,K}$ may be chosen to cover both bounds.
   Take $n$ sufficiently large that
   $C_{\sigma,K}B^2n^{-2}\log n\le1/2$.
   The nonlinear term then has norm at most $B/2$, and the
   Duhamel map has Lipschitz constant at most $1/2$.
   Together with the bound for the free term, this gives a
   contraction on the ball. Its fixed point satisfies
   $\sup_\tau\norm{u_p(\tau)}{K,\tau}\le B$, which implies the desired 
   \eqref{neg-eq:moments}.
   
\end{proof}

We now use Lemma~\ref{neg-lem:single-packet-bounds} to estimate
$U$ on $[-N/\ell^2,-n]$. The spatial moments control the mixed
products in \eqref{neg-eq:early-U-error}, while the Fourier
estimate controls the negative Sobolev norm. The following conclusions hold uniformly in $\theta$.

\begin{corollary}\label{neg-cor:early-profile}
	On $[-N/\ell^2,-n]$, $U$ satisfies both estimates in
	\eqref{neg-eq:rescaled-small} and the estimate in
	\eqref{neg-eq:rescaled-error}, with the time integrals
	restricted to this interval. Moreover, if $n$ tends to infinity,
	\begin{equation*}
		\norm{U(-N/\ell^2)}{L^2(\R^2)}\longrightarrow\sqrt3,
	\end{equation*}
	which ensures $\norm{Z_n(0)}{L^2(\R^2)}\longrightarrow\infty$.
	
\end{corollary}

\begin{proof}
	To prove \eqref{neg-eq:rescaled-error}, we estimate the
	products containing two distinct indices in
	\eqref{neg-eq:early-U-error}. Fix $p\ne q$. Since
	$|p-q|\ge1$,
	\begin{equation*}
		2|\tau|
		\le |2p\tau-2q\tau|
		\le |y-2p\tau|+|y-2q\tau|.
	\end{equation*}
	For $\tau\le-n$ and sufficiently large $n$,
	\begin{equation*}
		|\tau|-nd(\tau)
		=|\tau+n|\left(1-\frac4{n^3}\right)\ge0.
	\end{equation*}
	Thus every $y$ belongs to at least one of
	$\{|y-2p\tau|\ge nd(\tau)\}$ and
	$\{|y-2q\tau|\ge nd(\tau)\}$. By
	\eqref{neg-eq:packet-decay},
	\begin{equation*}
		\norm{\mathbf1_{\{|y-2p\tau|\ge nd(\tau)\}}u_p(\tau)}
		{L^2(\R^2)}
		\le n^{-K}
		\left\|\left\langle\frac{y-2p\tau}{d(\tau)}\right\rangle^K
		u_p(\tau)\right\|_{L^2(\R^2)}
		\lesssim_K n^{-K}.
	\end{equation*}
	The same estimate holds for $q$. Taking the other two
	factors in $L^\infty$, we obtain, for every
	$p'\in\{e_1,e_1+e_2,e_2\}$,
	\begin{equation*}
		\norm{u_p(\tau)u_q(\tau)u_{p'}(\tau)}{L^2(\R^2)}
		\lesssim_K n^{-K}d(\tau)^{-1}.
	\end{equation*}
	Conjugating any factor preserves this bound. Hence
	\eqref{neg-eq:early-U-error} and \eqref{neg-eq:eta} give
	\begin{equation*}
		\int_{-N/\ell^2}^{-n}
		\norm{(\ii\partial_\tau+N^{-1}\Delta_y)U-\mu n^{-6}F(U)}
		{L^2(\R^2)}\dd\tau
		\lesssim_{\sigma,K}n^{-K-2}\log n
		=o(n^{-7-4/\sigma}),
	\end{equation*}
	where $K>8+4/\sigma$. This proves
	\eqref{neg-eq:rescaled-error} on $[-N/\ell^2,-n]$.
	
	For the spacetime estimate, the triangle inequality,
	\eqref{neg-eq:packet-decay}, and \eqref{neg-eq:eta} give
	\begin{equation*}
		n^{-6}\int_{-N/\ell^2}^{-n}
		\norm{U(\tau)}{L^\infty(\R^2)}^2\dd\tau
		\lesssim_K n^{-6}\int_{-N/\ell^2}^{-n}
		d(\tau)^{-1}\dd\tau
		\lesssim_{\sigma,K}n^{-2}\log n
		\longrightarrow0.
	\end{equation*}
	This proves the second estimate in
	\eqref{neg-eq:rescaled-small} on this interval.
	
	For the negative Sobolev estimate, split each
	$\widehat{u_p}(\tau)$ into $\lbrace |\xi-Np|\le N/4\rbrace$ and $\lbrace|\xi-Np|> N/4\rbrace$. 
	Note that mass conservation and the normalization of $\psi$ give
	\begin{equation*}
		\norm{u_p(\tau)}{L^2(\R^2)}
		=\norm{u_p(-n)}{L^2(\R^2)}=1.
	\end{equation*}
	On the region $\lbrace |\xi-Np|\le N/4\rbrace$, we have $|\xi|\ge3N/4$, so 
	\begin{equation*}
		\left\|\mathbf1_{\{|\xi-Np|\le N/4\}}
		\langle\xi/\ell\rangle^{-\sigma}\widehat{u_p}(\tau,\xi)
		\right\|_{L^2(\R^2)}
		\le \left(\frac{4\ell}{3N}\right)^\sigma
		\norm{u_p(\tau)}{L^2(\R^2)}
		\lesssim_\sigma(\ell/N)^\sigma .
	\end{equation*}
	On the region $\lbrace |\xi-Np|> N/4\rbrace$, we apply \eqref{neg-eq:packet-fourier-tail} and obtain
	\begin{equation*}
		\left\|\langle\xi/\ell\rangle^{-\sigma}
		\widehat{u_p}(\tau,\xi)\right\|_{L^2(\R^2)}
		\lesssim_\sigma(\ell/N)^\sigma
		+\left\|\mathbf1_{\{|\xi-Np|>N/4\}}
		\widehat{u_p}(\tau,\xi)\right\|_{L^2(\R^2)}
		\lesssim_{\sigma,K}(\ell/N)^\sigma+n^{-4K}.
	\end{equation*}
	Consequently,
	\begin{equation*}
		\sup_{-N/\ell^2\le\tau\le-n}
		\left\|\langle\xi/\ell\rangle^{-\sigma}
		\widehat U(\tau,\xi)\right\|_{L^2(\R^2)}
		\lesssim_{\sigma,K}(\ell/N)^\sigma+n^{-4K}
		=o(n^{-7-4/\sigma}),
	\end{equation*}
	since
	\begin{equation*}
		n^{7+4/\sigma}(\ell/N)^\sigma
		=n^7N^{-\sigma}=n^{-1-20\sigma}\longrightarrow0,
		\qquad
		n^{7+4/\sigma-4K}\longrightarrow0.
	\end{equation*}
	This proves the first estimate in
	\eqref{neg-eq:rescaled-small}.
	
	Finally, $\{|\xi-Np|\le N/4\}$ are disjoint for different $p$.
	Plancherel, \eqref{neg-eq:packet-fourier-tail}, and
	$\norm{u_p(\tau)}{L^2(\R^2)}=1$ give
	\begin{equation*}
		\left|
		\norm{U(\tau)}{L^2(\R^2)}
		-\left(\sum_{p\in\{e_1,e_1+e_2,e_2\}}
		\left\|\mathbf1_{\{|\xi-Np|\le N/4\}}
		\widehat{u_p}(\tau)\right\|_{L^2(\R^2)}^2\right)^{1/2}
		\right|
		\lesssim_K n^{-4K},
	\end{equation*}
	\begin{equation*}
		\sum_{p\in\{e_1,e_1+e_2,e_2\}}
		\left\|\mathbf1_{\{|\xi-Np|\le N/4\}}
		\widehat{u_p}(\tau)\right\|_{L^2(\R^2)}^2
		=3+O_K(n^{-8K}).
	\end{equation*}
	Thus $\norm{U(\tau)}{L^2(\R^2)}
	=\sqrt3+O_K(n^{-4K})$ uniformly on the interval.
	
\end{proof}

On $[-n,0]$, the definition $U=P_k$ and
\eqref{neg-eq:collision-approx} give
\begin{equation}\label{neg-eq:iterate-residual}
	(\ii\partial_\tau+N^{-1}\Delta_y)U-\mu n^{-6}F(U)
	=\mu n^{-6}\bigl(F(P_{k-1})-F(P_k)\bigr).
\end{equation}
We therefore estimate $P_k$ in $L^\infty$ and
$P_k-P_{k-1}$ in $L^2$, respectively.

\begin{lemma}\label{neg-lem:finite-iterates}
	Uniformly in $\theta$, the iterates in
	\eqref{neg-eq:collision-approx} satisfy
	\begin{equation}\label{neg-eq:picard-differences}
		\sup_{-n\le\tau\le0}
		\norm{P_j(\tau)}{L^\infty(\R^2)}
		\lesssim_k1
		\qquad(0\le j\le k),
	\end{equation}
	\begin{equation}\label{azx}
		\sup_{-n\le\tau\le0}
		\norm{P_j(\tau)-P_{j-1}(\tau)}{L^2(\R^2)}
		\lesssim_k n^{-5j}
		\qquad(1\le j\le k).
	\end{equation}
\end{lemma}

\begin{proof}
	To estimate $P_j$ in $L^\infty$, repeatedly substitute
	\eqref{neg-eq:collision-approx}, replacing each $P_i$ with
	$i\ge1$ by its defining formula until only $P_0$ remains.
	Then substitute
	\begin{equation*}
		P_0=e^{\ii\theta}b\otimes a+b\otimes b
		+e^{3\ii\theta}a\otimes b .
	\end{equation*}
	Expanding each cubic and distributing the time integrals
	reduces the calculation to products and free evolutions
	of the one-dimensional Gaussians $a,b$, since
	\begin{equation*}
		(f_1\otimes g_1)\overline{(f_2\otimes g_2)}(f_3\otimes g_3)
		=(f_1\bar f_2f_3)\otimes(g_1\bar g_2g_3),
	\end{equation*}
	\begin{equation*}
		e^{\ii s\Delta_y/N}(f\otimes g)
		=\left(e^{\ii s\partial_x^2/N}f\right)
		\otimes\left(e^{\ii s\partial_x^2/N}g\right),
		\qquad
		\norm{f\otimes g}{L^\infty(\R^2)}
		=\norm f{L^\infty(\R)}\norm g{L^\infty(\R)}.
	\end{equation*}
	We estimate these Gaussians before performing the time
	integrations.
	
	For $\re z>0$ and real $\xi,x_0$, the free evolution is
	\begin{equation}\label{Gau}
		e^{\ii s\partial_x^2/N}
		\bigl[e^{\ii\xi x}e^{-z(x-x_0)^2}\bigr]
		=\frac{e^{\ii\xi x-\ii s\xi^2/N}}
		{\sqrt{1+4\ii sz/N}}
		\exp\left(
		-\frac{z(x-x_0-2s\xi/N)^2}{1+4\ii sz/N}
		\right).
	\end{equation}
	For the data defining $a,b$, respectively,
	$z=1$ and $z=1+\ii N/n^4$.
	If $\re z\sim_k1$, $|z|\lesssim_k N/n^4$, and
	$0\le s\le n$, then, for sufficiently large $n$,
	\begin{equation*}
		|4sz/N|\lesssim_k n^{-3},\qquad
		\re\left(\frac{z}{1+4\ii sz/N}\right)
		=\frac{\re z}{|1+4\ii sz/N|^2}\sim_k1,
	\end{equation*}
	\begin{equation*}
		\left|\frac{z}{1+4\ii sz/N}\right|
		\lesssim_k N/n^4,\qquad
		|1+4\ii sz/N|^{-1/2}\lesssim_k1.
	\end{equation*}
	Thus free evolution preserves these bounds on the quadratic
	coefficient and increases the supremum by at most a factor
	depending on $k$.
	
   Hence, for $G_i(x)=A_i e^{\ii\xi_i x}e^{-z_i(x-x_i)^2}$, the cubic product is (for some $A\in\C$ and $\xi\in\R$)
   \begin{equation*}
   	G_1\overline{G_2}G_3
   	=Ae^{\ii\xi x}e^{-z(x-x_0)^2},
   	\qquad
   	z=z_1+\bar z_2+z_3,\qquad
   	x_0=\frac{\sum_{i=1}^3(\re z_i)x_i}{\sum_{i=1}^3\re z_i},
   \end{equation*}
   where
   \begin{equation*}
   	\re z=\sum_{i=1}^3\re z_i\sim_k1,\qquad
   	|z|\le\sum_{i=1}^3|z_i|\lesssim_k N/n^4.
   \end{equation*}
   Then we can obtain
   \begin{equation*}
   	\norm{G_1\overline{G_2}G_3}{L^\infty(\R)}
   	=|A|
   	=\prod_{i=1}^3|A_i|\,
   	e^{-\sum_{i=1}^3(\re z_i)(x_0-x_i)^2}
   	\le\prod_{i=1}^3\norm{G_i}{L^\infty(\R)}.
   \end{equation*}
    Thus,
   these formulas bound every Gaussian integrand by $C_k$,
   uniformly in its time variables.
   
   A term with $m$ time integrations is therefore bounded in
   $L^\infty(\R^2)$ by
   \begin{equation*}
   	C_kn^{-6m}\int_{[-n,0]^m}1\dd r_1\cdots\dd r_m
   	=C_kn^{-5m}.
   \end{equation*}
   The number of terms and integrations depends only on $k$.
   Summing gives
   \begin{equation*}
   	\sup_{-n\le\tau\le0}
   	\norm{P_j(\tau)}{L^\infty(\R^2)}
   	\lesssim_k1,\qquad 0\le j\le k,
   \end{equation*}
   uniformly in $\theta$. This proves the $L^\infty$ estimate (\ref{neg-eq:picard-differences}).
	
	For the $L^2$ differences, unitarity gives
	$\sup_{-n\le\tau\le0}\norm{P_0(\tau)}{L^2(\R^2)}\le3$.
	Consequently,
	\begin{equation*}
		\sup_{-n\le\tau\le0}
		\norm{P_1(\tau)-P_0(\tau)}{L^2(\R^2)}
		\le n^{-6}\int_{-n}^0
		\norm{P_0(r)}{L^\infty(\R^2)}^2
		\norm{P_0(r)}{L^2(\R^2)}\dd r
		\lesssim_k n^{-5}.
	\end{equation*}
	For $1\le j<k$, subtracting the two Duhamel formulas gives
	\begin{equation*}
		P_{j+1}(\tau)-P_j(\tau)
		=-\ii\mu n^{-6}\int_{-n}^{\tau}
		e^{\ii(\tau-r)\Delta_y/N}
		\bigl(F(P_j(r))-F(P_{j-1}(r))\bigr)\dd r .
	\end{equation*}
	Using
	$|F(z)-F(w)|\lesssim(|z|^2+|w|^2)|z-w|$
	and the $L^\infty$ bounds just proved, we obtain
	\begin{equation*}
		\sup_{-n\le\tau\le0}
		\norm{P_{j+1}(\tau)-P_j(\tau)}{L^2(\R^2)}
		\lesssim_k n^{-5}
		\sup_{-n\le\tau\le0}
		\norm{P_j(\tau)-P_{j-1}(\tau)}{L^2(\R^2)} .
	\end{equation*}
	Induction proves the asserted $L^2$ estimate (\ref{azx}).
	
\end{proof}

To prove \eqref{neg-eq:rescaled-pairing}, we first compute the
pairing for $P_1(0)$. Lemma~\ref{neg-lem:finite-iterates} gives
\begin{equation*}
	\sup_\theta\norm{P_k(0)-P_1(0)}{L^2(\R^2)}
	\le\sum_{j=2}^k
	\sup_\theta\norm{P_j(0)-P_{j-1}(0)}{L^2(\R^2)}
	\lesssim_k n^{-10}.
\end{equation*}
Since $\norm{a(0)\otimes a(0)}{L^2(\R^2)}=1$, it suffices to consider the pairing with $P_1(0)$.

\begin{lemma}\label{neg-lem:terminal-pairing}
	For the maximizing $\theta$ chosen in the construction,
	\eqref{neg-eq:rescaled-pairing} holds.
\end{lemma}

\begin{proof}
	We first compute the Fourier coefficient of $e^{4i\theta}$.
	Note that the powers of $e^{\ii\theta}$ in
	$P_0$ are $1,0,3$, so
	\begin{equation*}
		\frac1{2\pi}\int_0^{2\pi}
		e^{-4\ii\theta}F(P_0(r))\dd\theta
		=2(b\otimes a)(r)\,
		\overline{(b\otimes b)(r)}\,(a\otimes b)(r),
	\end{equation*}

\begin{equation*}
	\frac1{2\pi}\int_0^{2\pi}e^{-4\ii\theta}P_0(0)\dd\theta=0.
\end{equation*}
Moreover, the free evolution of $a$ gives
\begin{equation*}
	e^{\ii r\Delta_y/N}\bigl(a(0)\otimes a(0)\bigr)
	=a(r)\otimes a(r).
\end{equation*}
Substituting the preceding formula for
$\frac1{2\pi}\int_0^{2\pi}e^{-4\ii\theta}F(P_0(r))\dd\theta$
into the iteration (\ref{neg-eq:collision-approx}), and moving the
propagator $e^{\ii r\Delta_y/N}$ to the second entry of the pairing, gives
\begin{equation*}
	\frac1{2\pi}\int_0^{2\pi}e^{-4\ii\theta}
	\langle P_1(0),a(0)\otimes a(0)\rangle\dd\theta
	=-2\ii\mu n^{-6}\int_{-n}^0
	\left\langle
	(b\otimes a)\overline{(b\otimes b)}(a\otimes b),
	a\otimes a\right\rangle(r)\dd r .
\end{equation*}
\begin{equation*}
	=-2\ii\mu n^{-6}\int_{-n}^0
	\left(\int_{\R}|a(r,x)|^2|b(r,x)|^2\dd x\right)^2\dd r .
\end{equation*}

To obtain a lower bound, restrict the last integral to the rectangle 
$\lbrace-1/4\le r\le0, \ |x|\le1/4\rbrace$. Applying \eqref{Gau} to the Gaussian initial data
$a(-n)$ and $b(-n)$ gives
\begin{equation*}
	|a(r,x)|
	=\frac{(2/\pi)^{1/4}}
	{|1+4\ii(r+n)/N|^{1/2}}
	\exp\left(-\frac{x^2}{|1+4\ii(r+n)/N|^2}\right),
\end{equation*}
\begin{equation*}
	|b(r,x)|
	=\frac{(2/\pi)^{1/4}}
	{|1-4(r+n)/n^4+4\ii(r+n)/N|^{1/2}}
	\exp\left(-\frac{(x-2r)^2}
	{|1-4(r+n)/n^4+4\ii(r+n)/N|^2}\right).
\end{equation*}
For sufficiently large $n$, uniformly on the rectangle,
\begin{equation*}
	\frac12\le |1+4\ii(r+n)/N|\le2,\qquad
	\frac12\le |1-4(r+n)/n^4+4\ii(r+n)/N|\le2,
	\qquad |x-2r|\le\frac34 .
\end{equation*}
Thus $|a(r,x)|,|b(r,x)|\gtrsim1$ there, and
\begin{equation*}
	\left|\frac1{2\pi}\int_0^{2\pi}e^{-4\ii\theta}
	\langle P_1(0),a(0)\otimes a(0)\rangle\dd\theta\right|
	\ge cn^{-6}
\end{equation*}
for some fixed $c>0$.
Lemma~\ref{neg-lem:finite-iterates} and
$\norm{a(0)\otimes a(0)}{L^2(\R^2)}=1$ give
\begin{equation*}
	\left|\frac1{2\pi}\int_0^{2\pi}e^{-4\ii\theta}
	\langle P_k(0)-P_1(0),a(0)\otimes a(0)\rangle\dd\theta\right|
	\le\sup_\theta\norm{P_k(0)-P_1(0)}{L^2(\R^2)}
	\le C_kn^{-10}.
\end{equation*}
Therefore
\begin{equation*}
	\max_{\theta\in[0,2\pi]}
	|\langle P_k(0),a(0)\otimes a(0)\rangle|
	\ge
	\left|\frac1{2\pi}\int_0^{2\pi}e^{-4\ii\theta}
	\langle P_k(0),a(0)\otimes a(0)\rangle\dd\theta\right|
	\ge cn^{-6}-C_kn^{-10}\gtrsim n^{-6}.
\end{equation*}
This proves \eqref{neg-eq:rescaled-pairing}.

\end{proof}

\begin{proof}[Proof of Proposition~\ref{neg-prop:profile}]
	
	On $[-n,0]$, we have $U=P_k$. By
	Lemma~\ref{neg-lem:finite-iterates},
	\begin{equation*}
		n^{-6}\int_{-n}^0
		\norm{U(\tau)}{L^\infty(\R^2)}^2\dd\tau
		\lesssim_k n^{-5}.
	\end{equation*}
	Combining this with Corollary~\ref{neg-cor:early-profile}
	gives
	\begin{equation*}
		n^{-6}\int_{-N/\ell^2}^0
		\norm{U(\tau)}{L^\infty(\R^2)}^2\dd\tau
		\lesssim_{\sigma,K,k}n^{-2}\log n+n^{-5}
		\longrightarrow0.
	\end{equation*}
	This completes \eqref{neg-eq:rescaled-small}.
	Applying
	\eqref{neg-eq:iterate-residual} and
	Lemma~\ref{neg-lem:finite-iterates}, we can also obtain
	\begin{equation*}
		n^{-6}\int_{-n}^0
		\norm{F(P_{k-1}(r))-F(P_k(r))}{L^2(\R^2)}\dd r
		\lesssim_k n^{-5(k+1)}
		=o(n^{-7-4/\sigma}),
	\end{equation*}
	where $5(k+1)>7+4/\sigma$.  This completes \eqref{neg-eq:rescaled-error}.
	Lemma~\ref{neg-lem:single-packet-bounds} gives smoothness
	with values in $\Sclass(\R^2)$ on $[-N/\ell^2,-n]$.
	The iterations (\ref{neg-eq:collision-approx}) give the same regularity on
	$[-n,0]$.
	
	Finally, the scaling in the construction gives
	\begin{equation*}
		\norm{Z_n(0)}{L^2(\R^2)}
		=\frac{\sqrt N}{n^3}
		\norm{U(-N/\ell^2)}{L^2(\R^2)}
		=n^{7+4/\sigma}\bigl(\sqrt3+o(1)\bigr)
		\longrightarrow\infty
	\end{equation*}
	by Corollary~\ref{neg-cor:early-profile}. Similarly,
	Lemma~\ref{neg-lem:terminal-pairing} yields
	\begin{equation*}
		|\langle Z_n(1),g_n\rangle|
		=\frac{\sqrt N}{n^3}\ell^\sigma
		|\langle P_k(0),a(0)\otimes a(0)\rangle|
		=n^7|\langle P_k(0),a(0)\otimes a(0)\rangle|
		\gtrsim n .
	\end{equation*}
	This concludes Proposition~\ref{neg-prop:profile}.
	
\end{proof}

\vspace{10pt}
\section{Adding one approximate solution to a smooth solution}
\label{sec:negative-correction}

We now add a translate of $Z_n$ to a fixed smooth solution.
Proposition~\ref{neg-prop:profile} makes both the error in its
equation and its interactions with that solution small.
The following lemma corrects their sum to an exact solution.

\begin{lemma}\label{neg-lem:correction}
	Let $v\in C([0,1];H^2(\R^2))$ solve
	\eqref{neg-eq:scalar}, and let $\delta_0>0$.
	For all sufficiently large $n$, uniformly in $d\in\R^2$,
	set $P(t,x):=Z_n(t,x-d)$. Then the solution $\widetilde v$
	of \eqref{neg-eq:scalar} with
	$\widetilde v(0)=v(0)+P(0)$ exists in
	$C([0,1];L^2(\R^2))\cap L^4([0,1]\times\R^2)$ and satisfies
	\begin{equation}\label{neg-eq:correction-bound}
		\norm{\widetilde v-v-P}{L^\infty([0,1];L^2(\R^2))}
		+\norm{\widetilde v-v-P}{L^4([0,1]\times\R^2)}
		\le\delta_0 .
	\end{equation}
\end{lemma}

\begin{proof}
	We seek $w$ with $w(0)=0$ such that
	$\widetilde v=v+P+w$ solves \eqref{neg-eq:scalar}.
	Since $v$ is an exact solution, the error in the equation
	for $v+P$ is
	\begin{equation*}
		E:=(\ii\partial_t+\Delta)P-\mu F(P)
		-\mu\bigl(2|v|^2P+v^2\bar P+2v|P|^2+\bar vP^2\bigr).
	\end{equation*}
	H\"older's inequality and
	$v\in C([0,1];H^2(\R^2))$ give
	\begin{equation*}
		\norm E{L^1([0,1];L^2(\R^2))}
		\lesssim_v\norm{(\ii\partial_t+\Delta)P-\mu F(P)}
		{L^1([0,1];L^2(\R^2))}
		+\left(
		\norm P{L^2([0,1];L^\infty(\R^2))}
		+\norm P{L^2([0,1];L^\infty(\R^2))}^2\right).
	\end{equation*}
	By \eqref{neg-eq:profile-small} and
	\eqref{neg-eq:profile-error},
	\begin{equation}\label{neg-eq:residual}
		\sup_{d\in\R^2}\norm E{L^1([0,1];L^2(\R^2))}
		\longrightarrow0.
	\end{equation}
	Uniformity in $d$ follows from translation invariance of
	the norms of $P$.
	
	The standard contraction argument, applied successively
	on a fixed partition of $[0,1]$ on which
	$\norm v{L^2(K;L^\infty(\R^2))}$ is small, now gives
	$w\in C([0,1];L^2(\R^2))\cap L^4([0,1]\times\R^2)$
	with $w(0)=0$ and
	\begin{equation*}
		\norm w{L^\infty([0,1];L^2(\R^2))}
		+\norm w{L^4([0,1]\times\R^2)}
		\lesssim_v \norm E{L^1([0,1];L^2(\R^2))}.
	\end{equation*}
	For sufficiently large $n$, \eqref{neg-eq:residual}
	makes the last bound at most $\delta_0$, uniformly in $d$.
	Thus $\widetilde v=v+P+w$ has the prescribed initial value
	and satisfies \eqref{neg-eq:correction-bound}.
	
\end{proof}

\vspace{10pt}
\section{Finite-time $H^s$-norm blow-up}
\label{sec:negative-limit}

\begin{proof}[Proof of Theorem~\ref{neg-thm:scalar}]
	Fix $s<0$ and $\varepsilon>0$, and take $\sigma=-s$ in
	Proposition~\ref{neg-prop:profile}. We construct exact
	solutions by adding translated profiles successively.
	The corrections will be summable in $L^2$, while the
	pairings with the translated functions $g_n$ will force
	the $H^s$ norm to blow up as $t\uparrow1$.
	
	Start with $u^{(0)}=0$. Suppose $u^{(n-1)}$ has been
	constructed on $[0,1]$. Apply Lemma~\ref{neg-lem:correction}
	with $v=u^{(n-1)}$ and $\delta_0=2^{-n}$.
	Choose an increasing sequence of profile indices $m_n$,
	large enough for that lemma to apply and such that
	$t_{m_n}\uparrow1$ and
	\begin{equation}\label{neg-eq:sequence}
		\sup_{0\le t<t_{m_n}}
		\norm{Z_{m_n}(t)}{H^s(\R^2)}
		\le\varepsilon2^{-n-1},\qquad
		\norm{Z_{m_n}(0)}{L^2(\R^2)}\ge n,\qquad
		|\langle Z_{m_n}(1),g_{m_n}\rangle|\ge n^2.
	\end{equation}
	These choices are possible by
	Proposition~\ref{neg-prop:profile}.
	
	We next choose translations to keep the pairings between
	different profiles small and to make their sum converge
	locally. Write
	\begin{equation*}
		P_n(t,x):=Z_{m_n}(t,x-X_n),\qquad
		G_n(x):=g_{m_n}(x-X_n),\qquad
		\chi_n:=\frac{P_n(0)}{\norm{P_n(0)}{L^2(\R^2)}}.
	\end{equation*}
	Choose $|X_n|$ sufficiently large so that, for every $j<n$,
	\begin{equation}\label{neg-eq:cross-pairings}
		|\langle P_j(1),G_n\rangle|
		+|\langle P_n(1),G_j\rangle|
		+|\langle P_j(0),\chi_n\rangle|
		+|\langle P_n(0),\chi_j\rangle|
		\le2^{-j-n},
	\end{equation}
	and, for $|\alpha|\le n$,
	\begin{equation}\label{neg-eq:weighted-translations}
		\sup_{0\le t\le1}
		\norm{\langle x\rangle^{-1}\partial_x^\alpha P_n(t,x)}
		{L^\infty(\R^2)}\le2^{-n}.
	\end{equation}
	Lemma~\ref{neg-lem:correction} now gives an exact solution
	$u^{(n)}$ on $[0,1]$ with
	\begin{equation*}
		u^{(n)}(0)=u^{(n-1)}(0)+P_n(0)
		=\sum_{j=1}^nP_j(0)\in\Sclass(\R^2).
	\end{equation*}
	Define $w_n:=u^{(n)}-u^{(n-1)}-P_n$. Then
	\begin{equation}\label{neg-eq:small-correction}
		w_n(0)=0,\qquad
		\norm{w_n}{L^\infty([0,1];L^2(\R^2))}
		+\norm{w_n}{L^4([0,1]\times\R^2)}
		\le2^{-n}.
	\end{equation}
	Persistence of regularity gives
	$u^{(n)}\in C([0,1];H^2(\R^2))$, as required at the next
	step. Since $P_n$ is continuous in $L^2$ across $t_{m_n}$,
	we also have $w_n\in C([0,1];L^2(\R^2))$.
	
	We first prove convergence before $t=1$. For each $T<1$,
	all sufficiently large $n$ satisfy $T<t_{m_n}$.
	Using $s<0$, \eqref{neg-eq:sequence}, and
	\eqref{neg-eq:small-correction}, we obtain
	\begin{equation*}
		\sup_{0\le t\le T}
		\norm{u^{(n)}(t)-u^{(n-1)}(t)}{H^s(\R^2)}
		\le
		\sup_{0\le t\le T}\norm{P_n(t)}{H^s(\R^2)}
		+\sup_{0\le t\le T}\norm{w_n(t)}{L^2(\R^2)}
		\le\varepsilon2^{-n-1}+2^{-n}.
	\end{equation*}
	Thus there exists a limit $u\in C([0,T];H^s(\R^2))$ such that
	\begin{equation*}
		u^{(m)}=\sum_{n=1}^m(P_n+w_n)
		\longrightarrow u
		\quad\text{in }C([0,T];H^s(\R^2))
		\qquad\text{for every }T<1.
	\end{equation*}
	
	To extend the limit to $t=1$ and pass to the equation,
	we use \eqref{neg-eq:weighted-translations}.
	For every bounded ball $B\subset\R^2$,
	\begin{equation*}
		\sum_{n=1}^\infty
		\norm{P_n}{L^\infty([0,1]\times B)}
		\le \sup_{x\in B}\langle x\rangle
		\sum_{n=1}^\infty2^{-n}<\infty.
	\end{equation*}
	For every $\varphi\in\Sclass(\R^2)$, we also have
	\begin{equation*}
		\sup_{0\le t\le1}
		\left|\left\langle\sum_{n>m}P_n(t),\varphi\right\rangle\right|
		\le \sum_{n>m}2^{-n}
		\norm{\langle x\rangle\varphi}{L^1(\R^2)}
		\longrightarrow0.
	\end{equation*}
	Together with \eqref{neg-eq:small-correction}, these estimates
	define the same limit by
	\begin{equation*}
		u(t)=\sum_{n=1}^\infty P_n(t)+\sum_{n=1}^\infty w_n(t),
		\qquad 0\le t\le1,
	\end{equation*}
	and give
	\begin{equation*}
		u\in C([0,1];L^2_{\mathrm{loc}}(\R^2))
		\cap C([0,1];\Sclass'(\R^2)),
		\qquad
		u^{(m)}\longrightarrow u
		\quad\text{in }L^4([0,1]\times B).
	\end{equation*}
	H\"older's inequality yields
	\begin{equation*}
		\norm{F(u^{(m)})-F(u)}{L^{4/3}([0,1]\times B)}
		\lesssim
		\bigl(\norm{u^{(m)}}{L^4([0,1]\times B)}^2
		+\norm u{L^4([0,1]\times B)}^2\bigr)
		\norm{u^{(m)}-u}{L^4([0,1]\times B)}
		\longrightarrow0.
	\end{equation*}
	Passing to the limit in the equations for $u^{(m)}$
	therefore gives \eqref{neg-eq:scalar} in the distributional sense.
	
	We next verify the assertions about the initial datum.
	Since $w_n(0)=0$,
	\begin{equation*}
		u(0)=\sum_{n=1}^\infty P_n(0),\qquad
		\norm{u(0)}{H^s(\R^2)}
		\le\sum_{n=1}^\infty\varepsilon2^{-n-1}
		=\varepsilon/2.
	\end{equation*}
	For every multi-index $\alpha$ and bounded ball $B$,
	\eqref{neg-eq:weighted-translations} gives
	\begin{equation*}
		\sum_{n\ge\max\{1,|\alpha|\}}
		\norm{\partial_x^\alpha P_n(0)}{L^\infty(B)}
		\le\sup_{x\in B}\langle x\rangle
		\sum_{n\ge\max\{1,|\alpha|\}}2^{-n}<\infty .
	\end{equation*}
	Thus $u(0)\in C^\infty(\R^2)$.
	On the other hand, \eqref{neg-eq:cross-pairings} gives
	\begin{equation*}
		|\langle u(0),\chi_n\rangle|
		\ge\norm{P_n(0)}{L^2(\R^2)}
		-\sum_{j\ne n}|\langle P_j(0),\chi_n\rangle|
		\ge n-1,
		\qquad
		\norm{\chi_n}{L^2(\R^2)}=1.
	\end{equation*}
	Cauchy--Schwarz would bound these pairings uniformly if
	$u(0)\in L^2(\R^2)$. Hence
	\begin{equation*}
		u(0)\in C^\infty(\R^2)\cap
		\bigl(H^s(\R^2)\setminus L^2(\R^2)\bigr).
	\end{equation*}
	
	Finally, we prove blow-up of the $H^s$ norm.
	By \eqref{neg-eq:small-correction},
	\begin{equation*}
		\left\|\sum_{j=1}^\infty w_j(1)\right\|_{L^2(\R^2)}
		\le1.
	\end{equation*}
	Since $\norm{G_n}{L^2(\R^2)}\le1$,
	\eqref{neg-eq:sequence} and \eqref{neg-eq:cross-pairings}
	give
	\begin{equation*}
		|\langle u(1),G_n\rangle|
		\ge|\langle P_n(1),G_n\rangle|
		-\sum_{j\ne n}|\langle P_j(1),G_n\rangle|
		-\left|\left\langle\sum_{j=1}^\infty w_j(1),G_n\right\rangle\right|
		\ge n^2-2.
	\end{equation*}
	Translation invariance and \eqref{neg-eq:test} imply
	$\sup_n\norm{G_n}{H^{-s}(\R^2)}\lesssim 1.$
	Thus distributional continuity at $t=1$
	 yields
	\begin{equation*}
		\liminf_{t\uparrow1}\norm{u(t)}{H^s(\R^2)}
		\gtrsim \lim_{t\uparrow1}|\langle u(t),G_n\rangle|
		=|\langle u(1),G_n\rangle|
		\ge n^2-2.
	\end{equation*}
	Letting $n\to\infty$ proves \eqref{neg-eq:blowup}.
	
\end{proof}

\vspace{10pt}
\section{Transfer of the blow-up examples to Hamiltonian systems}
\label{sec:negative-systems}

The scalar construction applies to both signs of the cubic nonlinearity.
We can therefore choose the sign determined by the invariant
line in Lemma~\ref{neg-lem:direction}.

\begin{proof}[Proof of Theorem~\ref{neg-thm:system}]
	Since $W\not\equiv0$, Lemma~\ref{neg-lem:direction} gives
	a unit vector $a$ and $\lambda=2W(a)\ne0$ satisfying
	\eqref{neg-eq:homogeneity}. Apply Theorem~\ref{neg-thm:scalar}
	with $\mu=\operatorname{sgn}\lambda$ and initial norm less
	than $\varepsilon\sqrt{|\lambda|}$. Set
	\begin{equation}\label{neg-eq:embedding}
		U:=\frac{a}{\sqrt{|\lambda|}}u .
	\end{equation}
	Then
	\begin{equation*}
		\partial_{\bar z}W(U)
		=\frac{\lambda a}{|\lambda|^{3/2}}|u|^2u
		=\frac{a}{\sqrt{|\lambda|}}(\ii\partial_t+\Delta)u
		=(\ii\partial_t+\Delta)U .
	\end{equation*}
	Multiplication by $a/\sqrt{|\lambda|}$ only scales every spatial
	Sobolev norm by $|\lambda|^{-1/2}$. Thus \eqref{neg-eq:embedding}
	preserves all the properties stated in Theorem~\ref{neg-thm:scalar}.
	
\end{proof}

\clearpage
\section*{Acknowledgement}
The author thanks Prof. Alex Cohen for providing financial support to cover the costs of ChatGPT usage during the preparation of this work.

\section*{Conflict of interest statement}
The author does not have any possible conflict of interest.

\section*{Data availability statement}
The manuscript has no associated data.

\bibliographystyle{plain}
\bibliography{Planar_NLS_Hamiltonian}

\begin{thebibliography}{10}

\bibitem{CDS2012}
R{\'e}mi Carles, Eric Dumas, and Christof Sparber.
\newblock Geometric optics and instability for {NLS} and {Davey--Stewartson}
  models.
\newblock {\em J. Eur. Math. Soc. (JEMS)}, 14(6):1885--1921, 2012.

\bibitem{CCT}
Michael Christ, James Colliander, and Terence Tao.
\newblock Ill-posedness for nonlinear {Schr\"odinger} and wave equations.
\newblock \href{https://arxiv.org/abs/math/0311048}{arXiv:math/0311048}, 2003.

\bibitem{CK}
Michael Christ and Alexander Kiselev.
\newblock Maximal functions associated to filtrations.
\newblock {\em J. Funct. Anal.}, 179(2):409--425, 2001.

\bibitem{CKSTT}
James Colliander, Markus Keel, Gigliola Staffilani, Hideo Takaoka, and Terence
  Tao.
\newblock Almost conservation laws and global rough solutions to a nonlinear
  {Schr\"odinger} equation.
\newblock {\em Math. Res. Lett.}, 9(5--6):659--682, 2002.

\bibitem{CKSTT2}
James Colliander, Markus Keel, Gigliola Staffilani, Hideo Takaoka, and Terence
  Tao.
\newblock Resonant decompositions and the {$I$}-method for the cubic nonlinear
  {Schr\"odinger} equation on {$\mathbb{R}^2$}.
\newblock {\em Discrete Contin. Dyn. Syst.}, 21(3):665--686, 2008.

\bibitem{Dodson}
Benjamin Dodson.
\newblock Global well-posedness and scattering for the defocusing,
  {$L^2$}-critical, nonlinear {Schr\"odinger} equation when {$d=2$}.
\newblock {\em Duke Math. J.}, 165(18):3435--3516, 2016.

\bibitem{evans2022partial}
Lawrence~C Evans.
\newblock {\em Partial differential equations}, volume~19.
\newblock American mathematical society, 2022.

\bibitem{Glassey}
Robert~T. Glassey.
\newblock On the blowing up of solutions to the {Cauchy} problem for nonlinear
  {Schr\"odinger} equations.
\newblock {\em J. Math. Phys.}, 18(9):1794--1797, 1977.

\bibitem{GO}
Loukas Grafakos and Seungly Oh.
\newblock The {Kato--Ponce} inequality.
\newblock {\em Comm. Partial Differential Equations}, 39(6):1128--1157, 2014.

\bibitem{KT}
Markus Keel and Terence Tao.
\newblock Endpoint {Strichartz} estimates.
\newblock {\em Amer. J. Math.}, 120(5):955--980, 1998.

\bibitem{MSZ}
Li~Ma, Xianfa Song, and Lin Zhao.
\newblock On global rough solutions to a non-linear {Schr\"odinger} system.
\newblock {\em Glasg. Math. J.}, 51(3):499--511, 2009.

\bibitem{Masaki}
Satoshi Masaki.
\newblock On scalar-type standing-wave solutions to systems of nonlinear
  {Schr\"odinger} equations.
\newblock \href{https://arxiv.org/abs/2212.00754v2}{arXiv:2212.00754v2}, 2022.
\newblock Revised 10 February 2023.

\bibitem{Oh2017}
Tadahiro Oh.
\newblock A remark on norm inflation with general initial data for the cubic
  nonlinear {Schr\"odinger} equations in negative {Sobolev} spaces.
\newblock {\em Funkcial. Ekvac.}, 60(2):259--277, 2017.

\bibitem{PV}
Fabrice Planchon and Luis Vega.
\newblock Bilinear virial identities and applications.
\newblock {\em Ann. Sci. {\'E}c. Norm. Sup{\'e}r. (4)}, 42(2):261--290, 2009.

\bibitem{Tao}
Terence Tao.
\newblock {\em Nonlinear dispersive equations: Local and global analysis},
  volume 106 of {\em CBMS Regional Conference Series in Mathematics}.
\newblock American Mathematical Society, Providence, RI, 2006.

\bibitem{TaoPC}
Terence Tao.
\newblock A pseudoconformal compactification of the nonlinear {Schr\"odinger}
  equation and applications.
\newblock {\em New York J. Math.}, 15:265--282, 2009.

\bibitem{WangZakharov}
Jiajun Wang.
\newblock Almost sharp global well-posedness and low-regularity blow-up for the
  two-dimensional periodic {Zakharov} system.
\newblock \href{https://arxiv.org/abs/2611.11048}{arXiv:2611.11048}, 2026.

\bibitem{YZ}
Kailong Yang and Lifeng Zhao.
\newblock Global well-posedness and scattering for mass-critical, defocusing,
  infinite dimensional vector-valued resonant nonlinear {Schr\"odinger} system.
\newblock {\em SIAM J. Math. Anal.}, 50(2):1593--1655, 2018.

\end{thebibliography}
\noindent
{\small
	Courant Institute of Mathematical Sciences,
	New York University, New York, NY 10012, USA\par
	\noindent
	\textit{Email:} \texttt{jw9409@nyu.edu}\par
}
\end{document}